\pdfoutput=1 
\RequirePackage{ifpdf}
\ifpdf 
\documentclass[pdftex]{sigma}
\else
\documentclass{sigma}
\fi

\numberwithin{equation}{section}

\usepackage{lscape}
\usepackage[all]{xy}

\usepackage{DLdef}

\makeatletter

\renewenvironment{lmatrix}{%
 \arraycolsep=4pt \left(%
  \matrix@check\lmatrix\env@matrix
}{
  \endmatrix\right)%
}

\newcommand {\fbj}{{\mathfrak{bj}}}
\newcommand {\fbr}{{\mathfrak{br}}}
\newcommand {\fbrj}{{\mathfrak{brj}}}
\newcommand {\fbgl}{{\mathfrak{bgl}}}

\newcommand {\fB}{{\mathfrak{B}}}
\newcommand {\fel}{{\mathfrak{el}}}

\newcommand {\fwk}{{\mathfrak{wk}}}

\newcommand {\foo}{{\mathfrak{oo}}}
\newcommand {\fooc}{{\mathfrak{ooc}}}
\newcommand {\foc}{{\mathfrak{oc}}}
\newcommand {\fpec}{{\mathfrak{pec}}}

\newcommand {\size}{{\text{size}}}

\rmnameii{Charr}{char}

\renewcommand{\mbullet}{{\fontsize{16}{9pt}\selectfont
\text{\raisebox{-1.3pt}[0pt]{$\bullet$}}}}
\renewcommand{\mcirc}{{\fontsize{16}{9pt}\selectfont
\text{\raisebox{-1.3pt}[0pt]{$\circ$}}}}
\renewcommand{\motimes}{{\fontsize{10}{9pt}\selectfont
\text{\raisebox{0.3pt}[0pt]{$\otimes$}}}}

\begin{document}
\allowdisplaybreaks

\renewcommand{\thefootnote}{$\star$}

\renewcommand{\PaperNumber}{060}

\FirstPageHeading

\ShortArticleName{Classif\/ication of Finite Dimensional Modular Lie Superalgebras}

\ArticleName{Classif\/ication of Finite Dimensional Modular Lie\\ Superalgebras with Indecomposable Cartan Matrix\footnote{This paper is a
contribution to the Special Issue on Kac--Moody Algebras and Applications. The
full collection is available at
\href{http://www.emis.de/journals/SIGMA/Kac-Moody_algebras.html}{http://www.emis.de/journals/SIGMA/Kac-Moody{\_}algebras.html}}}

\Author{Sof\/iane BOUARROUDJ ${}^\dag$, Pavel GROZMAN ${}^\ddag$ and Dimitry
LEITES ${}^\S$}

\AuthorNameForHeading{S.~Bouarroudj, P.~Grozman and D.~Leites}

\Address{${}^\dag$~Department of Mathematics, United Arab Emirates
University,\\
\hphantom{${}^\dag$}~Al Ain, PO. Box: 17551, United Arab Emirates}
\EmailD{\href{mailto:Bouarroudj.sofiane@uaeu.ac.ae}{Bouarroudj.sofiane@uaeu.ac.ae}}

\Address{${}^\ddag$~Equa Simulation AB, Stockholm, Sweden}
\EmailD{\href{mailto:pavel.grozman@bredband.net}{pavel.grozman@bredband.net}}

\Address{${}^\S$~Department of Mathematics, University of Stockholm,\\
\hphantom{${}^\S$}~Roslagsv. 101, Kr\"aft\-riket hus 6, SE-106 91 Stockholm, Sweden}
\EmailD{\href{mailto:mleites@math.su.se}{mleites@math.su.se}}

\ArticleDates{Received September 17, 2008, in f\/inal form May 25,
2009; Published online June 11, 2009}

\Abstract{Finite dimensional modular Lie superalgebras over
algebraically closed f\/ields with indecomposable Cartan matrices are
classif\/ied under some technical, most probably inessential,
hypotheses. If the Cartan matrix is invertible, the corresponding
Lie superal\-geb\-ra is simple otherwise the quotient of the derived Lie
super\-al\-geb\-ra modulo center is simple (if its rank is greater than~1). Eleven new exceptional simple modular Lie superalgebras are
discovered. Several features of classic notions, or notions
themselves, are clarif\/ied or introduced, e.g., Cartan matrix,
several versions of restrictedness in characteristic~2, Dynkin
diagram, Chevalley generators, and even the notion of Lie
superalgebra if the characteristic is equal to~2.
Interesting phenomena in characteristic~2: (1)~all simple Lie
superalgebras with Cartan matrix  are obtained from simple Lie
algebras with Cartan matrix by declaring several (any) of its
Chevalley generators odd; (2)~there exist simple Lie superalgebras
whose even parts are solvable.
The Lie superalgebras of f\/ixed points of automorphisms corresponding
to the symmetries of Dynkin diagrams are also listed and their
simple subquotients described.}


\Classification{17B50; 70F25}

\setcounter{tocdepth}{1}\tableofcontents

\renewcommand{\thefootnote}{\arabic{footnote}}
\setcounter{footnote}{0}

\section{Introduction}\label{Sintr}

The ground f\/ield $\Kee$ is algebraically closed of characteristic
$p>0$. (Algebraic closedness of $\Kee$ is only needed in the quest
for parametric families.) The term ``Cartan matrix'' will often
be abbreviated to CM. Let the {\it size} of the square $n\times n$
matrix be equal to $n$.

Except the last section, we will only consider {\it indecomposable}
Cartan matrices $A$ with entries in $\Kee$ (so the integer entries
are considered as elements of $\Zee/p\Zee\subset\Kee$),
indecomposability being the two conditions:
\begin{enumerate}\itemsep=0pt
  \item $A_{ij}=0\Longleftrightarrow A_{ji}=0$.
  \item By a
reshuf\/f\/le of its rows and columns $A$ can not be reduced to a
block-diagonal form.
\end{enumerate}

\subsection{Main results}
\begin{enumerate}\itemsep=0pt
\item  We clarify several key notions~-- of Lie
superalgebra in characteristic 2, of Lie superalgebra with Cartan
matrix, of weights and roots.

\item  We introduce several versions of restrictedness for Lie
(super)algebras in characteristic 2.

These clarif\/ications are obtained by/with A.~Lebedev.

\item  We give an algorithm that, under certain (conjecturally
immaterial) hypotheses, see Section~\ref{conj}.5), produces the
complete list of all f\/inite dimensional Lie superalgebras possessing
indecomposable Cartan matrices $A$, i.e., of the form $\fg(A)$. Our
proof follows the same lines Weisfeiler and Kac outlined for the Lie
algebra case in \cite{WK}. The result of application of the
algorithm~-- the classif\/ication -- is summarized in the following theorem:
\end{enumerate}

\begin{Theorem}[Sections~\ref{Sans>5}--\ref{Sans2}] \qquad{}
\renewcommand{\theenumi}{$\arabic{enumi}$}
\begin{enumerate}\itemsep=0pt
  \item There are listed all finite dimensional Lie superalgebras $($which are not Lie
algebras$)$ of the form $\fg(A)$ with indecomposable $A$ over $\Kee$.
  \item There are listed all inequivalent systems of simple roots $($inequivalent Cartan
matrices$)$ of each of the above listed Lie superalgebras $\fg(A)$.
  \item For the above listed Lie superalgebras $\fg=\fg(A)$, their even
  parts $\fg_\ev$ and the
$\fg_\ev$-modules $\fg_\od$ are explicitly described in terms of
simple and solvable Lie algebras and modules over them.
\item For the above listed Lie superalgebras $\fg=\fg(A)$, the
existence of restrictedness is explicitly established; for $p=2$,
various cases of restrictedness are considered and explicit formulas
given in each case.
\end{enumerate}
The results for $p>5$, $p=5$, $3$ and $2$ are summarized in Sections~{\rm \ref{Sans>5}}, {\rm \ref{Sans5}}, {\rm \ref{Sans3}} and {\rm \ref{Sans2}},
respectively.
\end{Theorem}

The following simple Lie superalgebras (where $\fg$ is the Lie
superalgebra with Cartan matrix and $\fg^{(1)}/\fc$ is the quotient
of its f\/irst derived algebra modulo the center) are new:
\begin{enumerate}\itemsep=0pt
  \item[1)] $\fe(6,1)$, $\fe(6,6)$, $\fe(7,1)^{(1)}/\fc$, $\fe(7,6)^{(1)}/\fc$,
  $\fe(7,7)^{(1)}/\fc$,
$\fe(8,1)$, $\fe(8, 8)$; $\fbgl(4; \alpha)$ and $\fbgl(3;
\alpha)^{(1)}\!/\fc$ for $p=2$;
  \item[2)] $\fel(5; 3)$ for $p=3$;
  \item[3)] $\fbrj(2; 5)$ for $p=5$.
\end{enumerate}
Observe that although several of the exceptional examples were known
for $p>2$, together with one indecomposable Cartan matrix per each
Lie superalgebra \cite{El1, El2, CE, CE2}, the complete
description of all inequivalent Cartan matrices for all the
exceptional Lie superalgebras of the form $\fg(A)$ and for ALL cases
for $p=2$ is new.

{\it A posteriori} we see that for each f\/inite dimensional Lie
superalgebra $\fg(A)$ with indecomposable Cartan matrix, the module
$\fg_\od$ is a completely reducible $\fg_\ev$-module\footnote{For
the {\it simple} subquotient $\fg=\fg(A)^{(i)}/\fc$ of $\fg(A)$,
where $i=\text{corank}\;A$ is equal to $1$ or 2, complete reducibility
of the $\fg_\ev$-module $\fg_\od$ is sometimes violated.}.

\subsection*{Fixed points of the symmetries of Dynkin diagrams}

We also listed the Lie superalgebras of f\/ixed points of
automorphisms corresponding to the symmetries of Dynkin diagrams and
described their simple subquotients. In characteristic 0, this is
the way all Lie algebras whose Dynkin diagrams has multiple bonds
(roots of dif\/ferent lengths) are obtained. Since, for $p=2$, there
are no multiple bonds or roots of dif\/ferent length (at least, this
notion is not invariant), it is clear that this is the way to obtain
something new, although, perhaps, not simple. Lemma 2.2 in \cite{FG}
implicitly describes the ideal in the Lie algebra of f\/ixed points of
an automorphism of a Lie algebra, but one still has to describe the
Lie algebra of f\/ixed points explicitly. This explicit answer is
given in the last section. No new simple Lie (super)algebras are
obtained.

\subsubsection[On simple subquotients of Lie (super)algebras of the form
$\fg(A)$  and a terminological problem]{On simple subquotients of Lie (super)algebras of the form
$\boldsymbol{\fg(A)}$\\ and a terminological problem}

 Observe that if a given
indecomposable Cartan matrix $A$ is invertible, the Lie
(super)algebra $\fg(A)$ is simple, otherwise $\fg(A)^{(i)}/\fc$
-- the quotient of its f\/irst (if $i=1$) or second (if $i=2$)
derived algebra modulo the center -- is simple (if $\size A >1$).
\textit{Hereafter in this situation $i=1$ or~$2$}; meaning that the chain
of derived algebras stabilizes ($\fg(A)^{(j)}\simeq \fg(A)^{(i)}$
for any $j\geq i$). We will see {\it a posteriori} that $\dim
\fc=i=\text{corank}\;A$, the maximal possible value of $i$ above.

This simple Lie algebra $\fg(A)^{(i)}/\fc$ does NOT possess any
Cartan matrix. Except for Lie algeb\-ras over f\/ields of characteristic
$p=0$, this subtlety is never mentioned causing confusion: The
conventional sloppy practice is to refer to the simple Lie
(super)algebra $\fg(A)^{(i)}/\fc$ as ``possessing a Cartan
matrix'' (although it does not possess any) and at the same time to
say ``$\fg(A)$ is simple'' whereas it is not. However, it is
indeed extremely inconvenient to be completely correct, especially
in the cases where $\fg(A)$, as well as its simple subquotient
$\fg(A)^{(i)}/\fc$, and a central extension of the latter, the
algebra $\fg(A)^{(i)}$, or all three at the same time might be
needed. So we suggest to refer to either of these three algebras as
``simple'' ones, and extend the same convention to Lie
superalgebras.

Thus, for $i=1$, there are three distinct Lie (super)algebras:
$\fg(A)$ with Cartan matrix, $\fg(A)^{(1)}$ its derived, and
$\fg(A)^{(1)}/\fc$ which is simple. It sometimes happens that the
three versions are needed at the same time and to appropriately
designate them is a non-trivial task dealt with in Section~\ref{Sg(A)}. When the
second derived of $\fg(A)$ is not isomorphic to the f\/irst one, i.e.,
when $i=2$, there are even more intermediate derived objects, but
fortunately only the three of them~-- $\fg(A)$ with Cartan matrix,
$\fg(A)^{(2)}$, and $\fg(A)^{(2)}/\fc$, the quotient modulo the
whole center~-- appear (so far) in applications.

\subsubsection{The Elduque Supermagic Square}

For $p>2$, Elduque interpreted
most of the exceptional (when their exceptional nature was only
conjectured; now this is proved) simple Lie superalgebras (of the
form $\fg(A)$) in charac\-te\-ris\-tic~3~\cite{CE2} in terms of super
analogs of division algebras and collected them into a Supermagic
Square (an analog of Freudenthal's Magic Square). The rest of the
exceptional examples for $p=3$ and $p=5$, not entering the Elduque
Supermagic Square (the ones described here for the f\/irst time) are,
nevertheless, somehow af\/f\/iliated to the Elduque Supermagic Square~\cite{El3}.

\subsection{Characteristic 2}

Very interesting, we think, is the situation in
characteristic 2. {\it A posteriori} we see that the list of Lie
{\bf super}algebras in characteristic 2 of the form $\fg(A)$ or
$\fg(A)^{(i)}/\fc$, where $i=\text{corank}\;A$ can be equal to either
$1$ or 2, with an indecomposable matrix $A$ is as follows:
\begin{enumerate}\itemsep=0pt
  \item Take any f\/inite dimensional Lie algebra of the form
$\fg(A)$ with indecomposable Cartan matrix~\cite{WK} and declare
some of its Chevalley generators (simultaneously a positive one and
the respective negative one) odd (the corresponding diagonal
elements of $A$ should be changed accordingly: $\ev$ to $0$ and
$\od$ to $1$, see Section~\ref{normA}). Let $I$ be the vector of
parities of the generators; the parity of each positive generator
should equal to that of the corresponding negative one.
  \item Do this for each of the inequivalent Cartan matrices of $\fg(A)$ and any
distribution $I$ of parities.
  \item Construct Lie superalgebra $\fg(A, I)$ from these Chevalley generators
by factorizing a certain $\widetilde\fg(A, I)$ modulo an ideal
(explicitly described in~\cite{LCh, BGL1, BGL2}; for a summary, see~\cite{BGLL}).
\item For the Lie superalgebra $
\fg(A, I)$, list all its inequivalent Cartan matrices.
\item If $A$ is not
invertible, pass to $\fg^{(i)}(A, I)/\fc$.
\end{enumerate}

Such superization turns
\begin{enumerate}\itemsep=0pt
  \item[1)] a given orthogonal Lie algebra into either
an ortho-orthogonal or a periplectic Lie superalgebra;
  \item[2)] the three
exceptional Lie algebras of $\fe$ type turn into seven
non-isomorphic Lie superalgebras of $\fe$ type;
 \item[3)] the $\fwk$
type Lie algebras discovered in \cite{WK} turn into $\fbgl$ type Lie
superalgebras.
\end{enumerate}

\ssbegin{Remarks} {\samepage
1) Observe that the Lie (super)algebra uniformly
def\/ined for any characteristic~as
\begin{itemize}\itemsep=0pt
  \item preserving a tensor, e.g., (super)trace-less; (ortho-)orthogonal
  and periplectic
ones,
  \item given by an integer Cartan matrix whose entries are considered
modulo $p$
\end{itemize}}
may have dif\/ferent (super)dimensions  as $p$ varies (not only from 2
to ``not 2''): The ``same'' Cartan matrices might def\/ine
algebras of similar type but with dif\/ferent properties and names as
the characteristic changes: for example $\fsl(np)$ has CM in all
characteristics, except $p$, in which case it is $\fgl(np)$ that has
a CM.

2) Although the number of inequivalent Cartan matrices grows with
the size of $A$, it is easy to list all possibilities for serial Lie
(super)algebras. Certain exceptional Lie superalgebras have dozens
of inequivalent Cartan matrices; nevertheless, there are at least
the following reasons to list all of them:
\begin{enumerate}\itemsep=0pt
  \item To classify all $\Zee$-gradings of a given $\fg(A)$ (in
particular, inequivalent Cartan matrices) is a very natural problem.
Besides, sometimes the knowledge of the best, for the occasion,
$\Zee$-grading is important. Examples of dif\/ferent cases: \cite{RU}
(all simple roots non-isotropic), \cite{LSS} (all simple roots odd);
for computations ``by hand'' the cases where only one simple
root is odd are useful.

In particular, the def\/ining relations between the natural
(Chevalley) generators of $\fg(A)$ are of completely dif\/ferent form
for inequivalent $\Zee$-gradings, and this is used in \cite{RU}.
  \item Distinct $\Zee$-gradings yield distinct
Cartan--Tanaka--Shchepochkina (CTS) prolongs (vectorial Lie
(super)algebras), cf.~\cite{Shch, L}. So to classify them is vital,
for example, in the quest for simple vectorial Lie (super)algebras.
  \item Certain properties of Cartan matrices may vary under the
  passage from one inequivalent CM to the other (the Lie superalgebras that
  correspond to such matrices may have dif\/ferent rate of growth as $\Zee$-graded
  algebras); this is a novel, previously
unnoticed, feature
  of Lie superalgebras that had lead to false claims (rectif\/ied in~\cite{CCLL}).
\end{enumerate}
\end{Remarks}

\subsection{Related results} 1) For explicit presentations in terms of
(the analogs of) Chevalley generators of the Lie algebras and
superalgebras listed here, see \cite{BGLL}. In addition to
Serre-type relations there always are more complicated relations.

2) For deformations of the f\/inite dimensional Lie (super)algebras of
the form $\fg(A)$ and $\fg(A)^{(i)}/\fc$, see \cite{BGL4}. Observe
that whereas ``{\sl for $p>3$, the Lie algebras with Cartan
matrices of the same types that exist over $\Cee$ are either rigid
or have deforms which also possess Cartan matrices}'', this is not so
if $p=3$ or $2$.

3) For generalized Cartan--Tanaka--Shchepochkina (CTS) prolongs of the
simple Lie (super)algebras of the form $\fg(A)$, and the simple
subquotients of such prolongs, see \cite{BGL3, BGL6}.

4) With restricted Lie algebras one can associate algebraic groups;
analogously, with restricted Lie superalgebras one can associate
algebraic supergroups. For this and other results of Lebedev's Ph.D.
thesis pertaining to the classif\/ication of simple modular Lie
superalgebras, see~\cite{LCh}.

\section{Basics: Linear algebra in superspaces (from \cite{LSh})}\label{Sbasics}

\subsection{General notation} For further details on basics on Linear
Algebra in Superspaces, see~\cite{LSh}.

For the def\/inition of the term ``Lie superalgebra'',
especially, in characteristic~2, see Section~\ref{Ssalgin2}.

For the def\/inition of the Lie (super)algebras of the form $\fg(A)$,
i.e., with Cartan matrix $A$ (brief\/ly referred to {\it CM Lie
$($super$)$algebras}), and simple relatives of the Lie superalgebras of
the form $\fg(A)$ with indecomposable Cartan matrix $A$, see Section~\ref{Sg(A)}. For simplicity of typing, the $i$th incarnation of the
Lie (super)algebra $\fg(A)$ with the $i$th CM (according to the
lists given below) will be denoted by $i\fg(A)$.

For the split form of any simple Lie algebra $\fg$, we denote the
$\fg$-module with the $i$th fundamental weight $\pi_i$ by $R(\pi_i)$
(as in \cite{OV, Bou}; these modules are denoted by $\Gamma_i$ in~\cite{FH}). In particular, for $\fo(2k+1)$, the {\it spinor
representation} $\spin_{2k+1}$ is {\it defined} to be the $k$th
fundamental representation, whereas for $\fo_\Pi(2k)$, the {\it
spinor representations} are the $k$th and the $(k-1)$st fundamental
representations. The {\it realizations} of these representations and
the corresponding modules by means of quantization (as in~\cite{LSh1}) can also be def\/ined, since, as is easy to see, the same
quantization procedure is well-def\/ined for the restricted version of
the Poisson superalgebra. (Fortunately, we do not need irreducible
representations of $\fo_{\rm I}(2k)$; their description is unknown, except
the trivial, identity and adjoint ones.)

In what follows, $\ad$ denotes both the adjoint representation and
the module in which it acts. Similarly, $\id$ denotes both the {\it
identity} (a.k.a.\ {\it standard}) representation of the linear Lie
(super)algebra $\fg\subset \fgl(V)$ in the (super)space $V$ and $V$
itself. (We disfavor the adjective ``natural'' applied to $\id$
only, since it is no less appropriate to any of the tensor
(symmetric, exterior) powers of $\id$.) In particular, having f\/ixed
a basis in the $n$-dimensional space and having realized $\fsl(V)$
as $\fsl(n)$, we write $\id$ instead of $V$, so $V$ does not
explicitly appear. The {\it exterior powers}\index{power, exterior}
and {\it symmetric powers}\index{power, symmetric} of the vector
space $V$ are def\/ined as quotients of its tensor powers
\begin{gather*} T^0(V):=\bigwedge^0(V):=S^0(V):=\Kee,\\
T^1(V):=\bigwedge^1(V):=S^1(V):=V,\\
T^i(V):=\underbrace{V\otimes \dots \otimes V}_{i
\text{~factors}}\text{~~ for $i>0$}.
\end{gather*} We set:
\[\wedge{}^{\bcdot}(V):=T^{\bcdot}(V)/(x\otimes x\mid x\in V),\qquad
S^{\bcdot}(V):=T^{\bcdot}(V)/(x\otimes y+y\otimes x\mid x, y\in V),
\]
where $T^{\bcdot}(V):=\oplus T^{i}(V)$; let $\bigwedge^{i}(V)$ and
$S^{i}(V)$ be homogeneous components of degree $i$.

Describing the $\fg_\ev$-module structure of $\fg_\od$ for a Lie
superalgebra $\fg$, we write $\fg_\od\simeq R(\bcdot)$, though it
is, actually, $\Pi(R(\bcdot))$.

The symbol $A\subplus B$ denotes a {\it semi-direct sum} of modules
of which $A$ is a submodule; when dealing with algebras, $A$ is an
ideal in $A\subplus B$.

\subsection{Superspaces} A {\it superspace} is a~$\Zee /2$-graded
space; for any superspace $V=V_{\ev}\oplus V_{\od }$, where $\ev$
and $\od$ are residues modulo~2; we denote by $\Pi (V)$ another copy
of the same superspace: with the shifted parity, i.e.,
$(\Pi(V))_{\bar i}= V_{\bar i+\od }$. The parity function is denoted
by $p$. The {\it superdimension} of $V$ is $\sdim V=a+b\eps $, where
$\eps ^2=1$, and $a=\dim V_{\ev}$, $b=\dim V_{\od }$. (Usually,
$\sdim V$ is shorthanded as a~pair $(a,b)$ or $a|b$; this notation
obscures the fact that $\sdim V\otimes W=\sdim V\cdot \sdim W$.)

A superspace structure in $V$ induces the superspace structure in
the space $\End (V)$. A {\it basis of a~superspace} is always
a~basis consisting of {\it homogeneous} vectors; let $\Par=(p_1,
\dots, p_{\dim V})$ be an ordered collection of their parities. We
call $\Par$ the {\it format} of (the basis of) $V$. A square {\it
supermatrix} of format (size) $\Par$ is a~$\sdim V\times \sdim V$
matrix whose $i$th row and $i$th column are of the same parity
$p_i\in\Par$. The matrix unit $E_{i, j}$ is of parity $p_i+p_j$.

Whenever possible, we consider one of the simplest formats $\Par$,
e.g., the format $\Par_{st}$ of the form $(\ev , \dots, \ev ; \od ,
\dots, \od)$ is called {\it standard}. Systems of simple roots of
Lie superalgebras corresponding to distinct nonstandard formats of
supermatrix realizations of these superalgebras are related by
so-called {\it odd reflections}.

\subsection{The Sign Rule} The formulas of Linear Algebra are superized
by means of the {\bf Sign Rule}:
\begin{center}
\begin{minipage}[l]{12cm}
{\it if something of parity $p$ moves past something of parity $q$,
the~sign~$(-1)^{pq}$ accrues; the expressions defined on homogeneous
elements are extended to arbitrary ones via linearity}.
\end{minipage}
\end{center}

Examples of application of the Sign Rule: By setting \[ [X,
Y]=XY-(-1)^{p(X)p(Y)}YX\] we get the notion of the supercommutator
and the ensuing notions of {\it supercommutative} and {\it
superanti-commutative} superalgebras; {\it Lie superalgebra} is the
one which, in addition to superanti-commutativity, satisf\/ies Jacobi
identity amended with the Sign Rule; a {\it superderivation of
a~given superalgebra} $A$ is a~linear map $D\colon A\tto A$
satisfying the {\it super Leibniz rule} \[
D(ab)=D(a)b+(-1)^{p(D)p(a)}aD(b). \]

Let $V$ be a superspace, $\sdim V=m|n$. The {\it general linear} Lie
superalgebra of operators acting in $V$ is denoted by $\fgl(V)$ or
$\fgl(m|n)$ if an homogeneous basis of $V$ is f\/ixed. The Lie
subsuperalgebra of supertraceless operators (supermatrices) is
denoted $\fsl(V)\simeq \fsl(m|n)$. If a Lie superalgebra $\fg\subset
\fgl(V)$ contains the ideal of scalar matrices $\fs=\Kee 1_{m|n}$,
then $\fp\fg=\fg/\fs$ denotes the projective version of $\fg$.

Observe that sometimes the Sign Rule requires some dexterity in
application. For example, we have to distinguish between super-skew
and super-anti although both versions coincide in the non-super
case:
\begin{alignat*}{3}
& ba=(-1)^{p(b)p(a)}ab \qquad && (\text{supercommutativity}), &\\
& ba=-(-1)^{p(b)p(a)}ab \qquad &&(\text{superanti-commutativity}), &\\
& ba=(-1)^{(p(b)+1)(p(a)+1)}ab \qquad &&(\text{superskew-commutativity}), & \\
& ba=-(-1)^{(p(b)+1)(p(a)+1)}ab\qquad & &(\text{superantiskew-commutativity}).
\end{alignat*}
In other words, ``anti'' means the total change of the sign, whereas
any ``skew'' notion can be straightened by the parity change. In
what follows, the supersymmetric bilinear forms and supercommutative
superalgebras are named according to the above def\/initions.

The supertransposition of supermatrices is def\/ined so as to assign
the supertransposed supermatrix to the dual operator in the dual
bases. For details, see \cite{LSoS}; an explicit expression in the
standard format is as follows (we give a general formula for
matrices with entries in a~supercommutative superalgebra; in this
paper we only need the lower formula):
\[
X=\mat{a&b\\c&d}\mapsto
X^{st}:=\begin{cases}\mat{a^t&c^t\\-b^t&d^t}&\text{for $X$
even},\vspace{1mm}\\
\mat{a^t&-c^t\\b^t&d^t}&\text{for $X$ odd}.
 \end{cases}
\]

\subsection{What the~Lie superalgebra is}\label{SS:2.2}

Dealing with
superalgebras it sometimes becomes useful to know their def\/inition.
Lie superalgebras were distinguished in topology in 1930's, and the
Grassmann superalgebras half a century earlier. So it might look
strange when somebody of\/fers a~``better'' def\/inition of a~notion
which was established about 70 year ago. Nevertheless, the answer to
the question ``what is a~(Lie) superalgebra?'' is still not a common
knowledge.

So far we def\/ined Lie superalgebras naively: via the Sign Rule.
However, the naive def\/inition suggested above (``apply the Sign Rule
to the def\/inition of the Lie algebra'') is manifestly inadequate for
considering the supervarieties\footnote{A {\it
supervariety}\index{supervariety} is a ringed space such that the
collection of {\it functions} on it~-- the sections of its sheaf~-- constitute a supercommutative superring. {\it Morphisms} of
supervarieties are only the ring space morphisms that preserve
parity of the superrings of sections of the structure sheaves.} of
deformations and for applications of representation theory to
mathematical physics, for example, in the study of the coadjoint
representation of the Lie supergroup which can act on
a~supermanifold or supervariety but never on a~vector superspace~-- an object from another category. We were just lucky in the case
of f\/inite dimensional Lie algebras over $\Cee$ that the vector
spaces can be viewed as manifolds or varieties. In the case of
spaces over $\Kee$ and in the super setting, to be able to deform
Lie (super)algebras or to apply group-theoretical methods, we must
be able to recover a~supermanifold from a~vector superspace, and
{\it vice versa}.

A proper def\/inition of Lie superalgebras is as follows. The {\it Lie
superalgebra} in the category of supervarieties corresponding to the
``naive'' Lie superalgebra $L= L_{\ev} \oplus L_{\od}$ is a~linear
supermanifold $\cL=(L_{\ev}, \cO)$, where the sheaf of functions
$\cO$ consists of functions on $L_{\ev}$ with values in the
Grassmann superalgebra on $L_{\od}^*$; this supermanifold should be
such that for ``any'' (say, f\/initely generated, or from some other
appropriate category) supercommutative superalgebra $C$, the space
$\cL(C)=\Hom (\Spec C, \cL)$, called {\it the space of $C$-points of} $\cL$, is a~Lie algebra and the correspondence $C\tto
\cL(C)$ is a~functor in $C$. (A.~Weil introduced this approach in
algebraic geometry in 1953; in super setting it is called {\it the
language of points} or {\it families}.) This def\/inition might look
terribly complicated, but fortunately one can show that the
correspondence $\cL\longleftrightarrow L$ is one-to-one and the Lie
algebra $\cL(C)$, also denoted $L(C)$, admits a~very simple
description: $L(C)=(L\otimes C)_{\ev}$.

A {\it Lie superalgebra homomorphism} $\rho\colon L_1 \tto L_2$ in
these terms is a~functor morphism, i.e., a~collection of Lie algebra
homomorphisms $\rho_C\colon L_1 (C)\tto L_2(C)$ such that any
homomorphism of supercommutative superalgebras $\varphi\colon C\tto
C_1$ induces a~Lie algebra homomorphism $\varphi\colon L(C)\tto
L(C_1)$ and products of such homomorphisms are naturally compatible.
In particular, a~{\it representation} of a~Lie superalgebra $L$ in
a~superspace $V$ is a~homomorphism $\rho\colon L\tto \fgl (V)$,
i.e., a~collection of Lie algebra homomorphisms $\rho_C\colon L(C)
\tto ( \fgl (V )\otimes C)_{\ev}$.

\sssbegin{Example} Consider a~representation $\rho\colon
\fg\tto\fgl(V)$. The space of inf\/initesimal deformations of $\rho$
is isomorphic to $H^1(\fg; V\otimes V^*)$. For example, if $\fg$ is
the $0|n$-dimensional (i.e., purely odd) Lie superalgebra
(with the only bracket possible: identically equal to zero), its
only irreducible modules are the trivial one, ${\bf 1}$, and
$\Pi({\bf 1})$. Clearly, \[{\bf 1}\otimes {\bf 1}^*\simeq \Pi({\bf
1})\otimes \Pi({\bf 1})^*\simeq {\bf 1},\] and, because the Lie
superalgebra $\fg$ is commutative, the dif\/ferential in the cochain
complex is zero. Therefore $H^1(\fg; {\bf
1})=\bigwedge^1(\fg^*)\simeq\fg^*$, so there are $\dim\,\fg$ odd
parameters of deformations of the trivial representation. If we
consider $\fg$ ``naively'', all of these odd parameters will be
lost.
\end{Example}

Examples that lucidly illustrate why one should always remember that
a~Lie superalgebra is not a~mere linear superspace but a~linear
supermanifold are, e.g., the deforms with odd parameters. In the
category of supervarieties, these deforms, listed in \cite{BGL4},
are simple Lie superalgebras.

\subsection{Examples of simple Lie superalgebras over
$\Cee$}\label{SS:2.2.1}

Recall that the Lie superalgebra $\fg$
without proper ideals and of dimension $>1$ is said to be {\it
simple}. Examples: Serial Lie superalgebras $\fsl(m|n)$ for $m>
n\geq 1$, $\fpsl(n|n):=\fsl(n|n)/\fc$ for $n>1$, $\fosp(m|2n)$ for
$mn\neq 0$, and $\fspe(n)$ for $n>2$; and the exceptional Lie
superalgebras: $\fag(2)$, $\fab(3)$, and $\fosp(4|2;\alpha)$ for
$\alpha\neq 0, -1$.

\section{What the Lie superalgebra in characteristic 2 is
(from \cite{Leb1})}\label{Ssalgin2}

Let us give a naive def\/inition of a Lie superalgebra for $p=2$. We
def\/ine it as a superspace $\fg=\fg_\ev\oplus\fg_\od$ such that
$\fg_\ev$ is a Lie algebra, $\fg_\od$ is an $\fg_\ev$-module (made
into the two-sided one by anti-symmetry, but if $p=2$, it is the
same) and on $\fg_\od$ a {\it squaring} (roughly speaking, the
halved bracket) is def\/ined as a map
\begin{gather}\begin{split}
& x\mapsto x^2\quad \text{such that $(ax)^2=a^2x^2$ for any $x\in
\fg_\od$ and $a\in \Kee$, and} \\
& {}(x+y)^2-x^2-y^2\text{~is a bilinear form on $\fg_\od$ with values
in $\fg_\ev$.}
\end{split}\label{squaring}
\end{gather}
(We use a minus sign, so the def\/inition also works for $p\neq 2$.)
The origin of this operation is as follows: If $\Char \Kee\neq 2$,
then for any Lie superalgebra $\fg$ and any odd element
$x\in\fg_\od$, we have $x^2=\frac12 [x,x]\in\fg_\ev$. If $p=2$, we
def\/ine $x^2$ f\/irst, and then def\/ine the bracket of odd elements to
be (this equation is valid for $p\neq 2$ as well):
\begin{equation}\label{bracket}
{}[x, y]:=(x+y)^2-x^2-y^2.
\end{equation}
We also assume, as usual, that
\begin{center}
\begin{minipage}[l]{12cm}
if $x,y\in\fg_\ev$, then $[x,y]$ is the bracket on the Lie algebra;

if $x\in\fg_\ev$ and $y\in\fg_\od$, then
$[x,y]:=l_x(y)=-[y,x]=-r_x(y)$, where $l$ and $r$ are the left and
right $\fg_\ev$-actions on $\fg_\od$, respectively.
\end{minipage}
\end{center}

The Jacobi identity involving odd elements has now the following
form:
\begin{equation}\label{JI}
~[x^2,y]=[x,[x,y]]\quad \text{for any~} x\in\fg_\od, y\in\fg.
\end{equation}
If $\Kee\neq \Zee/2\Zee$, we can replace the condition on two odd
elements by a simpler one:
\begin{equation}\label{JIbest}
[x,x^2]=0\quad \text{for any $x\in\fg_\od$.}
\end{equation}

Because of the squaring, the def\/inition of derived algebras should
be modif\/ied. For any Lie superalgebra $\fg$, set $\fg^{(0)}:=\fg$
and
\begin{equation}\label{deralg}
\fg^{(1)}:=[\fg,\fg]+\Span\{g^2\mid g\in\fg_\od\} ,\qquad
\fg^{(i+1)}:=[\fg^{(i)},\fg^{(i)}]+\Span\{g^2\mid
g\in(\fg^{(i)})_\od\}.
\end{equation}

\subsection{Examples: Lie superalgebras preserving non-degenerate forms
\cite{Leb1}}

Lebedev investigated various types of equivalence of
bilinear forms for $p=2$, see \cite{Leb1}; we just recall the
verdict and say that {\it two (anti)-symmetric bilinear forms $B$
and $B'$ on a superspace~$V$ are equivalent} if there is an even
non-degenerate linear map $M\colon V\to V$ such that
\begin{equation}\label{eqform}
B'(x,y)=B(Mx,My) \quad \text{for all~~} x,y\in V.
\end{equation} \textit{We fix some basis in $V$ and identify a given
bilinear form with its
Gram matrix in this basis; let us also identify any linear operator
on $V$ with its matrix}. Then two bilinear forms (rather
supermatrices) are {\it equivalent} if there is an even invertible
matrix $M$ such that
\begin{equation}\label{eqformM}
B'=MBM^T,\quad \text{where $T$ is for transposition}.
\end{equation}

We often use the following matrices
\begin{equation}\label{matrices}
J_{2n}=\mat{0&1_n\\
-1_n&0},\qquad \Pi_n=\begin{cases} \mat{0&1_k\\
1_k&0}&\text{if $n=2k$},\\[3mm]
\mat{0&0&1_k\\
0&1&0\\
1_k&0&0}&\text{if $n=2k+1$}.\end{cases}
\end{equation}
Let $J_{n|n}$ and $\Pi_{n|n}$ be the same as $J_{2n}$ and $\Pi_{2n}$
but considered as supermatrices.

Lebedev proved that, with respect to the above natural\footnote{It
is interesting and unexpected that for non-symmetric bilinear forms,
another equivalence is more natural.} equivalence of forms
(\ref{eqformM}), every even symmetric non-degenerate form on a
superspace of dimension $n_{\ev}|n_{\od}$ over a perfect f\/ield of
characteristic $2$ is equivalent to a form of the shape (here:
$i=\bar 0$ or $\bar 1$ and each $n_i$ may equal to~0)
\[
B=\mat{ B_{\ev}&0\\0&B_{\od}}, \qquad \text{where
$B_i=\begin{cases}1_{n_i}&\text{if $n_i$ is odd,}\\
\text{either $1_{n_i}$ or $\Pi_{n_i}$}&\text{if $n_i$ is
even.}\end{cases}$}
\]
In other words, the bilinear forms with matrices $1_{n}$ and
$\Pi_{n}$ are equivalent if $n$ is odd and non-equivalent if $n$ is
even; $\antidiag(1, \dots, 1)\sim \Pi_{n}$ for any $n$. The precise
statement is as follows:

\sssbegin{Theorem}\label{SymForm}\label{s2.2.1} Let $\Kee$ be a
perfect field of characteristic $2$. Let $V$ be an $n$-dimensio\-nal
space over $\Kee$.

\textup{1)} For $n$ odd, there is only one equivalence class of
non-degenerate symmetric bilinear forms on $V$.

\textup{2)} For $n$ even, there are two equivalence classes of
non-degenerate symmetric bilinear forms, one~-- with at least one
non-zero element on the main diagonal~-- contains $1_n$ and the
other one~-- with only $0$s on the main diagonal~--  contains
$S_n:=\antidiag(1, \dots, 1)$ and $\Pi_n$.
\end{Theorem}

The Lie superalgebra preserving $B$~-- by analogy with the orthosymplectic Lie superalgebras $\fosp$ in
characteristic $0$ we call it {\it ortho-orthogonal} and denote
$\fo\fo_B(n_\ev|n_\od)$~-- is spanned by the supermatrices which in the standard format are
of the form
\[
\mat{ A_{\ev}&B_{\ev}C^TB_{\od}^{-1}\\C&A_{\od}}, \quad
\begin{matrix}\text{where
$A_{\ev}\in\fo_{B_{\ev}}(n_\ev)$, $A_{\od}\in\fo_{B_{\od}}(n_\od)$, and}\\
\text{$C$ is arbitrary $n_{\od}\times n_{\ev}$ matrix.}\end{matrix}
\]
Since, as is easy to see, \[\fo\fo_{\Pi I}(n_\ev|n_\od)\simeq
\fo\fo_{\rm I\Pi}(n_\od|n_\ev),\] we do not have to consider the Lie
superalgebra $\fo\fo_{\Pi I}(n_\ev|n_\od)$ separately unless we
study Cartan prolongations~-- the case where the dif\/ference between
these two incarnations of the same algebra is vital.

For an odd symmetric form $B$ on a superspace of dimension
$(n_{\ev}|n_{\od})$ over a f\/ield of characteristic $2$ to be
non-degenerate, we need $n_{\ev}=n_{\od}$, and every such form $B$
is equivalent to $\Pi_{k|k}$, where $k=n_{\ev}=n_{\od}$. This form
is preserved by linear transformations with supermatrices in the
standard format of the shape
\begin{equation}
\label{pe} \mat{ A&C\\D&A^T}, \qquad \text{where $A\in\fgl(k)$, $C$
and $D$ are symmetric
$k\times k$ matrices}. 
\end{equation}
The Lie superalgebra $\fpe(k)$ of supermatrices $(\ref{pe})$ will be
referred to as {\it periplectic}, as A.~Weil suggested, and denoted
by $\fpe_B(k)$ or just $\fpe(k)$. (Notice that the matrix
realization of $\fpe_B(k)$ over $\Cee$ or $\Ree$ is dif\/ferent: its
set of roots is not symmetric relative the change of roots ``positive $\longleftrightarrow$ negative''.)
\begin{equation}\label{nt}
\begin{tabular}{c}
{\bf The fact that two bilinear forms are inequivalent does}\\
{\bf not, generally, imply that the Lie (super)algebras that}\\
{\bf preserve them are not isomorphic}.
\end{tabular}
\end{equation}
In \cite{Leb1}, Lebedev proved that for the {\it non-degenerate
symmetric} forms, this implication (\ref{nt}) is, however, true
(except for $\foo_{\rm I\Pi}(n_\ev|n_\od)\simeq\foo_{\Pi
I}(n_\od|n_\ev)$ and $\foo_{\Pi\Pi}^{(1)}(6|2)\simeq\fpe^{(1)}(4)$)
and described the distinct types of Lie (super)algebras preserving
non-degenerate forms. In what follows, we describe which of these
Lie (super)algebras (or their derived ones) are simple, and which of
them (or their central extensions) ``have Cartan matrix''. But
f\/irst, we recall what does the term in quotation marks mean.

\section[What $\fg(A)$ is]{What $\boldsymbol{\fg(A)}$ is}\label{Sg(A)}

\subsection[Warning: certain of $\fsl$'s and all $\fpsl$'s have no Cartan
matrix.
Which of their relatives have Cartan matrices]{Warning: certain of $\boldsymbol{\fsl}$'s and all $\boldsymbol{\fpsl}$'s have no Cartan
matrix.\\
Which of their relatives have Cartan matrices}\label{warn}

For the most reasonable def\/inition of Lie algebra with Cartan matrix
over $\Cee$, see \cite{K}. The same def\/inition applies, practically
literally, to Lie superalgebras and to modular Lie algebras and to
modular Lie superalgebras. However, the usual sloppy practice is to
attribute Cartan matrices to many of those (usually simple) modular
Lie algebras and (modular or not) Lie superalgebras which, strictly
speaking, have no Cartan matrix!

Although it may look strange for the reader with non-super
experience over $\Cee$, neither the simple modular Lie algebra
$\fpsl(pk)$, nor the simple modular Lie superalgebra
$\fpsl(a|pk+a)$, nor~-- in characteristic~$0$~-- the simple Lie
superalgebra $\fpsl(a|a)$ possesses a Cartan matrix. Their central
extensions~-- $\fsl(pk)$, the modular Lie superalgebra
$\fsl(a|pk+a)$, and~-- in characteristic~$0$~-- the Lie
superalgebra $\fsl(a|a)$~-- do not have Cartan matrix, either.

Their relatives possessing a Cartan matrix are, respectively,
$\fgl(pk)$, $\fgl(a|pk+a)$, and $\fgl(a|a)$, and for the ``extra'' (from the point of view of $\fsl$ or $\fpsl$) grading
operator (such operators are denoted in what follows $d_i$, to
distinguish them from the ``inner'' grading operators $h_i$) we
take~$E_{1,1}$.

Since often all the Lie (super)algebras involved (the simple one,
its central extension, the derivation algebras thereof) are needed
(and only representatives of one of the latter types of Lie
(super)algebras are of the form $\fg(A)$), it is important to have
(preferably short and easy to remember) notation for each of them.

For the Lie algebras that preserve a tensor (bilinear form or a
volume element) we retain the same notation in all characteristic;
but the (super)dimension of various incarnations of the algebras
with the same name may dif\/fer as characteristic changes; simplicity
may also be somewhat spoilt.

For the Lie (super)algebras which are easiest to be determined by
their Cartan matrix (the Elduque Supermagic Square is of little help
here), we have:
\begin{enumerate}\itemsep=0pt
\item[]\underline{for $p=3$}: $\fe(6)$ is of dimension 79, its derived
$\fe(6)^{(1)}$ is of dimension 78, whereas its ``simple core''
is $\fe(6)^{(1)}/\fc$ of dimension 77;

\item[]$\fg(2)$ is not simple (moreover, its CM is decomposable); its
``simple core'' is isomorphic to $\fpsl(3)$;

\item[]\underline{for $p=2$}: $\fe(7)$ is of dimension 134, its derived
$\fe(7)^{(1)}$ is of dimension 133, whereas its ``simple core''
is $\fe(7)^{(1)}/\fc$ of dimension 132;

\item[]$\fg(2)$ constructed from its CM reduced modulo 2 is now isomorphic
to $\fgl(4)$; it is not simple, its ``simple core'' is
isomorphic to $\fpsl(4)$;

\item[]\underline{the orthogonal Lie algebras and their super analogs} are
considered in detail later.
\end{enumerate}

In our examples, the notation $D/d|B$ means that $\sdim \fg(A)=D|B$
whereas $\sdim \fg(A)^{(i)}/\fc=d|B$, where
\begin{equation}\label{dims}
d=D-2(\size(A)-\rk(A))\qquad \text{and} \qquad i=\size(A)-\rk(A)=\dim\fc.
\end{equation}

\subsection{Generalities} Let us start with the construction of a CM Lie
(super)algebra. Let $A=(A_{ij})$ be an $n\times n$-matrix. Let $\rk
A=n-l$. It means that there exists an $l\times n$-matrix
$T=(T_{ij})$ such that
\begin{gather}
\mbox{a) the rows of $T$ are linearly independent;}\nonumber\\
\mbox{b) $TA=0$ (or, more precisely, ``zero $l\times n$-matrix'')}.\label{rankCM}
\end{gather}
Indeed, if $\rk A^T=\rk A=n-l$, then there exist $l$ linearly
independent vectors $v_i$ such that $A^Tv_i=0$; set
\[
T_{ij}=(v_i)_j.
\]

Let the elements $e_i^\pm$ and $h_i$, where $i=1,\dots,n$, generate
a Lie superalgebra denoted $\tilde \fg(A, I)$, where $I=(p_1, \dots
p_n)\in(\Zee/2)^n$ is a collection of parities ($p(e_i^\pm)=p_i$),
free except for the relations
\begin{equation}\label{gArel_0}
{}[e_{i}^+, e_{j}^-] = \delta_{ij}h_i; \qquad [h_j, e_{j}^\pm]=\pm
A_{ij}e_{j}^\pm \qquad \text{and} \qquad [h_i, h_j]=0 \qquad \text{for any $i, j$}.
\end{equation}
Let Lie (super)algebras with Cartan matrix $ \fg(A,
I)$\index{$\fg(A, I)$} be the quotient of $\tilde \fg(A, I)$ modulo
the ideal we explicitly described in \cite{BGL1, BGL2, BGLL}.

By abuse of notation we retain the notations $e_j^\pm$ and $h_j$~--
the elements of $\tilde{\fg}(A,I)$~-- for their images in $\fg(A,I)$ and $\fg^{(i)}(A,I)$.

The additional to \eqref{gArel_0} relations that turn $\tilde \fg(A,
I)$ into $\fg(A, I)$ are of the form $R_i=0$ whose left sides are
implicitly described, for the general Cartan matrix with entries in
$\Kee$, as
\begin{equation}\label{myst}
\text{\begin{minipage}[c]{14cm} the $R_i$ that generate the ideal
$\fr$ maximal among the ideals of $\tilde\fg (A, I)$\\ whose
intersection with the span of the above $h_i$ and the $d_j$
described\\ in equation~\eqref{central3} is zero.\end{minipage}}
\end{equation}

Set
\begin{equation}\label{central}
c_i=\sum_{j=1}^n T_{ij}h_j, \qquad \text{where~~} i=1,\dots,l.
\end{equation}
Then, from the properties of the matrix $T$, we deduce that
\begin{equation}\label{central1}
\begin{tabular}{l}
a) the elements $c_i$ are linearly independent;\\
b) the elements $c_i$ are central, because\\
\qquad $[c_i,e_j^\pm]=\pm\left(\sum\limits_{k=1}^n T_{ik}A_{kj}\right)
e_j^\pm=\pm (TA)_{ij} e_j^\pm $.
\end{tabular}
\end{equation}
The existence of central elements means that the linear span of all
the roots is of dimension $n-l$ only. (This can be explained even
without central elements: The weights can be considered as
column-vectors whose $i$-th coordinates are the corresponding
eigenvalues of $\ad_{h_i}$. The weight of $e_i$ is, therefore, the
$i$-th column of $A$. Since $\rk A=n-l$, the linear span of all
columns of $A$ is $(n-l)$-dimensional just by def\/inition of the
rank. Since any root is an (integer) linear combination of the
weights of the $e_i$, the linear span of all roots is
$(n-l)$-dimensional.)

This means that some elements which we would like to see having
dif\/ferent (even opposite if $p=2$) weights have, actually, identical
weights. To f\/ix this, we do the following: Let $B$ be an arbitrary
$l\times n$-matrix such that
\begin{equation}\label{matrixB}
\text{the~~}(n+l)\times
n\text{-matrix~~}\begin{pmatrix}A\\ B\end{pmatrix}\text{~~has
rank~}n.
\end{equation}
Let us add to the algebra $\fg=\tilde \fg(A, I)$ (and hence $\fg(A,
I)$) the grading elements $d_i$, where $i=1,\dots,l$, subject to the
following relations:
\begin{equation}\label{central3}
{}[d_i,e_j^\pm]=\pm B_{ij}e_j;\quad [d_i,d_j]=0;\quad [d_i,h_j]=0
\end{equation}
(the last two relations mean that the $d_i$ lie in the Cartan
subalgebra, and even in the maximal torus which will be denoted by
$\fh$).

Note that these $d_i$ are {\it outer} derivations of $\fg(A,
I)^{(1)}$, i.e., they can not be obtained as linear combinations of
brackets of the elements of $\fg(A, I)$ (i.e., the $d_i$ do not lie
in $\fg(A, I)^{(1)}$).

\subsection{Roots and weights}\label{roots}
In this subsection, $\fg$
denotes one of the algebras $\fg(A,I)$ or $\tilde{\fg}(A,I)$.

Let $\fh$ be the span of the $h_i$ and the $d_j$. The elements of
$\fh^*$ are called {\it weights}.\index{weight} For a given weight
$\alpha$, the {\it weight subspace} of a given $\fg$-module $V$ is
def\/ined as
\[
V_\alpha=\{x\in V\mid \text{an integer $N>0$ exists such that
$(\alpha(h)-\ad_h)^N x=0$ for any $h\in\fh$}\}.
\]

Any non-zero element $x\in V$ is said to be {\it of weight
$\alpha$}. For the roots, which are particular cases of weights if
$p=0$, the above def\/inition is inconvenient because it does not lead
to the modular analog of the following useful statement.

\sssbegin{Statement}[\cite{K}]\label{rootdec} Over $\Cee$, the space
of any Lie algebra $\fg$ can be represented as a~direct sum of
subspaces
\begin{equation}\label{rootdeceq}
\fg=\mathop{\bigoplus}\limits_{\alpha\in \fh^*} \fg_\alpha.
\end{equation}
\end{Statement}
Note that if $p=2$, it might happen that $\fh\subsetneq\fg_0$. (For
example, all weights of the form $2\alpha$ over $\Cee$ become 0 over
$\Kee$.)

To salvage the formulation of Statement in the modular case with
minimal changes, \textit{at least for the Lie $($super$)$algebras $\fg$ with
Cartan matrix}~-- and only this case we will have in mind speaking
of roots, we decree that the elements $e_i^\pm$ with the same
superscript (either $+$ or $-$) correspond to linearly independent
{\it roots} $\alpha_i$, and any root $\alpha$ such that
$\fg_\alpha\neq 0$ lies in the $\Zee$-span of
$\{\alpha_1,\dots,\alpha_n\}$, i.e.,
\begin{equation}\label{rootdeceqnew}
\fg=\mathop{\bigoplus}\limits_{\alpha\in
\Zee\{\alpha_1,\dots,\alpha_n\}} \fg_\alpha.
\end{equation}

Thus, $\fg$ has a $\Ree^n$-grading such that $e_i^\pm$ has grade
$(0,\dots,0,\pm 1,0,\dots,0)$, where $\pm 1$ stands in the $i$-th
slot (this can also be considered as $\Zee^n$-grading, but we use
$\Ree^n$ for simplicity of formulations). If $p=0$, this grading is
equivalent to the weight grading of $\fg$. If $p>0$, these gradings
may be inequivalent; in particular, if $p=2$, then the elements
$e_i^+$ and $e_i^-$ have the same weight. (That is why in what
follows we consider roots as elements of $\Ree^n$, not as weights.)

Any non-zero element $\alpha\in\Ree^n$ is called {\it a
root}\index{root} if the corresponding eigenspace of grade $\alpha$
(which we denote $\fg_\alpha$ by abuse of notation) is non-zero. The
set $R$ of all roots is called {\it the root system}\index{root
system} of $\fg$.

Clearly, the subspaces $\fg_\alpha$ are purely even or purely odd,
and the corresponding roots are said to be \textit{even} or
\textit{odd}.

\subsection{Systems of simple and positive roots} In this subsection,
$\fg=\fg(A,I)$, and $R$ is the root system of $\fg$.

For any subset $B=\{\sigma_{1}, \dots, \sigma_{m}\} \subset R$, we
set (we denote by $\Zee_{+}$ the set of non-negative integers):
\[
R_{B}^{\pm} =\Big\{ \alpha \in R \mid \alpha = \pm \sum n_{i}
\sigma_{i},\;\text{where }\;n_{i} \in \Zee_{+} \Big\}.
\]
The set $B$ is called a {\it system of simple roots}
of $R$ (or $\fg$) if $ \sigma_{1}, \dots ,
\sigma_{m}$ are linearly independent and $R=R_B^+\cup R_B^-$. Note
that $R$ contains basis coordinate vectors, and therefore spans
$\Ree^n$; thus, any system of simple roots contains exactly $n$
elements.

A subset $R^+\subset R$ is called a {\it system of positive
roots} 
 of $R$ (or $\fg$) if there
exists $x\in\Ree^n$ such that
\begin{equation}\label{x}
(\alpha,x)\in\Ree\backslash \{0\}\text{ for all $\alpha\in R$},\qquad
 R^+=\{\alpha\in R\mid (\alpha,x)>0\}.
\end{equation} (Here $(\cdot,\cdot)$ is the standard Euclidean inner
product in $\Ree^n$). Since $R$ is a f\/inite (or, at least, countable
if $\dim \fg(A)=\infty$) set, so the set
\[\{y\in\Ree^n\mid\text{there exists $\alpha\in R$ such that }
(\alpha,y)=0\} \] is a f\/inite/countable union of $(n-1)$-dimensional
subspaces in $\Ree^n$, so it has zero measure. So for almost every
$x$, condition (\ref{x}) holds.

By construction, any system $B$ of simple roots is contained in
exactly one system of positive roots, which is precisely $R_B^+$.

\sssbegin{Statement} Any finite system $R^+$ of positive roots of
$\fg$ contains exactly one system of simple roots. This system
consists of all the positive roots $($i.e., elements of $R^+)$ that
can not be represented as a sum of two positive
roots.\end{Statement}

We can not give an {\it a priori} proof of the fact that each set of
all positive roots each of which is not a sum of two other positive
roots consists of linearly independent elements. This is, however,
true for f\/inite dimensional Lie algebras and Lie superalgebras of
the form $\fg(A)$ if $p\neq 2$.

\subsection{Normalization convention}\label{normA}

Clearly,
\begin{equation}
\label{rescale} \text{the rescaling
$e_i^\pm\mapsto\sqrt{\lambda_i}e_i^\pm$, sends $A$ to $A':=
\diag(\lambda_1, \dots , \lambda_n)\cdot A$.} 
\end{equation}

Two pairs $(A, I)$ and $(A', I')$ are said to be {\it equivalent} if
$(A', I')$ is obtained from $(A, I)$ by a composition of a
permutation of parities and a rescaling $A' = \diag (\lambda_{1},
\dots, \lambda_{n})\cdot A$, where $\lambda_{1}\cdots \lambda_{n}\neq
0$. Clearly, equivalent pairs determine isomorphic Lie
superalgebras.

The rescaling af\/fects only the matrix $A_B$, not the set of parities
$I_B$. The Cartan matrix $A$ is said to be {\it
normalized}\index{Cartan matrix, normalized} if
\begin{equation}
\label{norm} A_{jj}=0\quad \text{or 1, or 2.}
\end{equation}
We let $A_{jj}=2$ only if $i_j=\ev$; in order to distinguish between
the cases where $i_j=\ev$ and $i_j=\od$, we write $A_{jj}=\ev$ or
$\od$, instead of 0 or 1, if $i_j=\ev$. \textit{We will only consider
normalized Cartan matrices; for them, we do not have to describe
$I$.}

The row with a 0 or $\ev$ on the main diagonal can be multiplied by
any nonzero factor; usually (not only in this paper) we multiply the
rows so as to make $A_{B}$ symmetric, if possible.

\subsection{Equivalent systems of simple roots} \label{EqSSR} Let
$B=\{\alpha_1,\dots,\alpha_n\}$ be a system of simple roots. Choose
non-zero elements $e_i^\pm$ in the 1-dimensional (by def\/inition)
superspaces $\fg_{\pm\alpha_i}$; set $h_{i}=[e_{i}^{+}, e_{i}^-]$,
let $A_{B} =(A_{ij})$, where the entries $A_{ij}$ are recovered from
relations \eqref{gArel_0}, and let $I_{B}=\{p(e_{1}), \cdots,
p(e_{n})\}$. Lemma \ref{serg} claims that all the pairs $(A_B,I_B)$
are equivalent to each other.

Two systems of simple roots $B_{1}$ and $B_{2}$ are said to be {\it
equivalent} if the pairs $(A_{B_{1}}, I_{B_{1}})$ and $(A_{B_{2}},
I_{B_{2}})$ are equivalent.

It would be nice to f\/ind a convenient way to f\/ix some distinguished
pair $(A_B,I_B)$ in the equivalence class. For the role of the
``best'' (f\/irst among equals) order of indices we propose the
one that minimizes the value
\begin{equation}\label{minCM}
\max\limits_{i,j\in\{1,\dots,n\}\text{~such that~}(A_B)_{ij}\neq
0}|i-j|
\end{equation}
(i.e., gather the non-zero entries of $A$ as close to the main
diagonal as possible). Observe that this numbering dif\/fers from the
one that Bourbaki use for the $\fe$ type Lie algebras.

\subsubsection{Chevalley generators and Chevalley bases}\label{SsChev} We
often denote the set of generators corresponding to a normalized
matrix by $X_{1}^{\pm},\dots , X_{n}^{\pm}$ instead of
$e_{1}^{\pm},\dots , e_{n}^{\pm}$; and call them, together with the
elements $H_i:=[X_{i}^{+}, X_{i}^{-}]$, and the derivatives $d_j$
added for convenience for all $i$ and $j$, the {\it Chevalley
generators}.\index{Chevalley generator}

For $p=0$ and normalized Cartan matrices of simple f\/inite
dimensional Lie algebras, there exists only one (up to signs) basis
containing $X_i^\pm$ and $H_i$ in which $A_{ii}=2$ for all $i$ and
all structure constants are integer, cf.~\cite{St}. Such a basis is
called the {\it Chevalley}\index{Basis! Chevalley}  basis.

Observe that, having normalized the Cartan matrix of $\fo(2n+1)$ so
that $A_{ii}=2$ for all $i\neq n$ but $A_{nn}=1$, we get {\bf
another} basis with integer structure constants. We think that this
basis also qualif\/ies to be called {\it Chevalley basis}; for Lie
superalgebras, and if $p=2$, such normalization is a must.

\begin{Conjecture} If $p>2$, then for finite dimensional Lie
$($super$)$algebras with indecomposable Cartan matrices normalized as in
$(\ref{norm})$, there also exists only one $($up to signs$)$ analog of
the Chevalley basis. \end{Conjecture}

We had no idea how to describe analogs of Chevalley bases for $p=2$
until recently; clearly, the methods of the recent paper \cite{CR}
should solve the problem.

\section{Restricted Lie superalgebras}\label{Srestr}
Let $\fg$ be a~Lie algebra of characteristic $p>0$. Then, for every
$x\in \fg$, the operator $\ad_x^{p}$ is a~derivation of $\fg$. If it
is an {\it inner} derivation for every $x\in\fg$, i.e., if
$\ad_x^{p}= \ad_{x^{[p]}}$ for some element denoted $x^{[p]}$, then
the corresponding map
\begin{equation}
\label{p-struct} [p]\colon x\mapsto x^{[p]}
\end{equation}
is called a $p$-{\it structure}\index{$p$-structure} on $\fg$, and
the Lie algebra $\fg$ endowed with a~$p$-structure is called a~{\it
restricted}\index{Lie (super)algebra, restricted} Lie algebra. If
$\fg$ has no center, then $\fg$ can have not more than one
$p$-structure. The Lie algebra $\fgl(n)$ possesses a $p$-structure,
unique up to the contribution of the center; this $p$-structure is
used in the next def\/inition.

The notion of a $p$-{\it representation} is naturally def\/ined as a
linear map $\rho:\fg\tto\fgl(V)$ such that
$\rho(x^{[p]})=(\rho(x))^{[p]})$; in this case $V$ is said to be a
$p$-{\it module}.

Passing to superalgebras, we see that, for any odd $D\in \fder A$,
we have
\begin{equation}
\label{6eq38}
D^{2n}([a, b]) = \sum \binom{n}{l}[D^{2l}(a), D^{2n-2l}(b)]\quad \text{for any $a, b\in A$.}
\end{equation}
So, if $\Charr\; \Kee=p$, then $D^{2p}$ is always an even derivation
for any odd $D\in \fder A$. Now, let $\fg$ be a Lie superalgebra of
characteristic $p>0$. Then
\begin{center}
\begin{minipage}[l]{12cm}
for every $x \in\fg_\ev$, the operator $\ad_x^p$ is a derivation of
$\fg$, i.e., $\fg_\ev$-action on $\fg_\od$ is a $p$-representation;

for every $x\in\fg_\od$, the operator $\ad_x^{2p}=\ad_{x^2}^p$ is a
derivation of $\fg$.
\end{minipage}
\end{center}

So, if for every $x\in\fg_\ev$, there is $x^{[p]}\in\fg_\ev$ such
that $\ad_x^p=\ad_{x^{[p]}}$ for any $x\in\fg_\ev$, then we can
def\/ine $x^{[2p]}:=(x^2)^{[p]}$ for any $x\in\fg_\od$. We demand that
for any $x\in\fg_\ev$, we have
\[\ad_x^p=\ad_{x^{[p]}}\text{~~as operators on the whole
$\fg$, i.e., $\fg_\od$ is a restricted $\fg_\ev$-module.}\] Then the
pair of maps \[ [p]\colon\fg_\ev\tto\fg_\ev\quad(x\mapsto
x^{[p]})\qquad\text{and}\qquad [2p]\colon
\fg_\od\tto\fg_\ev\quad(x\mapsto x^{[2p]})\] is called a $p|2p$-{\it
structure}\index{$p\vert 2p$-structure}~-- or just $p$-structure
-- on $\fg$, and the Lie
superalgebra $\fg$ endowed with a~$p$-structure is called a~{\it
restricted} Lie superalgebra.\index{Lie superalgebra! restricted}

\subsection[The case where $\fg_\ev$ has center]{The case where $\boldsymbol{\fg_\ev}$ has center}

The $p$-structure on
$\fg_\ev$ does not have to determine a $p|2p$-structure on $\fg$:
Even if the actions of $\ad_x^p$ and $\ad_{x^{[p]}}$ coincide on
$\fg_\ev$, they do not have to coincide on the whole of $\fg$. This
remark af\/fects even simple Lie superalgebras if $\fg_\ev$ has
center. We can not say if a $p$-structure on $\fg_{\ev}$ def\/ines a
$p|2p$-structure on $\fg$ in the case of centerless $\fg_{\ev}$: To
def\/ine it we need to have, separately, a $p$-module structure on
$\fg_{\od}$ over $\fg_{\ev}$.

For the case where the Lie superalgebra $\fg$ or even $\fg_\ev$ has
center, the following def\/inition is more appropriate: $\fg$ is said
to be {\it restricted} if
\begin{equation}
\label{p2p-struct}
\renewcommand{\arraystretch}{1.4}
\begin{array}{ll}
\text{there is given the map~~}[p]\colon x\mapsto x^{[p]}&\text{for
any $x\in\fg_\ev$}
\end{array}
\end{equation}
such that, for any $a\in \Kee$, we have
\begin{equation}
\label{p|2pprop1}
\renewcommand{\arraystretch}{1.4}
\begin{array}{ll}
1)& (a x)^{[p]}=a^p \cdot x^{[p]}\text{~~for any $x\in\fg_\ev$, $a\in\Kee$},\\
2)&
(x+y)^{[p]}=x^{[p]}+y^{[p]}+\mathop{\sum}\limits^{p-1}_{i=1}s_i(x,
y) \text{~~for any $x, y\in\fg_\ev$},\\
&\text{where $is_i(x, y)$ is the coef\/f\/icient of $a^{i-1}$ in
the expression of $(\ad_{a x+y})^{p-1}(x)$}\\
&\text{for an indeterminate $a$}, \\
3)& [x^{[p]}, y]=(\ad_x)^{p}(y)\text{~~for any $x\in\fg_\ev$, $y\in
\fg$}.
\end{array}
\end{equation}
We set \[\begin{array}{ll} {} [2p]\colon x\mapsto
x^{[2p]}:=(x^2)^{[p]}&\text{for any $x\in\fg_\od$.}
\end{array}
\]

\sssbegin{Remark} If $\fg$ is centerless, we do not need conditions
1) and 2) of (\ref{p|2pprop1}) since they follow from~3).
\end{Remark}

\ssbegin{Proposition}\qquad{}
\begin{enumerate}\itemsep=0pt
\item[$1)$] If $p>2$ $($or $p=2$ but $A_{ii}\neq \od$ for
all $i)$ and $\fg(A)$ is finite-dimensional, then $\fg(A)$ has a~$p|2p$-structure such that
\begin{equation}\label{p-str}
(x_\alpha)^{[p]}=0\text{~for any even $\alpha\in R$ and $x_\alpha\in\fg_\alpha$},\qquad \fh^{[p]}\subset \fh.
\end{equation}

\item[$2)$] If all the entries of $A$ are elements of $\Zee/p\Zee$, then we
can set $h_i^{[p]}=h_i$ for all $i=1,\dots,n$.

\item[$3)$] The quotient modulo center of $\fg(A)$ or $\fg^{(1)}(A)$ always
inherits the $p$-structure of $\fg(A)$ or $\fg^{(1)}(A)$ $($if any$)$
whereas $\fg^{(1)}(A)$ does not necessarily inherit the
$p$-structure of $\fg(A)$.
\end{enumerate}
\end{Proposition}

\begin{proof} For the simple Lie algebras, the $p$-structure is unique if any exists,
see \cite{J}. The same proof applies to simple Lie superalgebras and
$p|2p$-structures. The explicit construction completes the proof of
headings 1) and 2). To prove 3) a counterexample suf\/f\/ices; we leave
it as an exercise to the reader to produce one.
\end{proof}

\sssbegin{Remarks} \qquad
\begin{enumerate}\itemsep=0pt
\item[1)] It is not enough to def\/ine $p|2p$-structure on
generators, one has to def\/ine it on a basis.

\item[2)] If $p(X_i^\pm)=\od$, then $(X_i^\pm)^{[p]}$ is not def\/ined unless
$p=2$: Only $(X_i^\pm)^{[2p]}$ is def\/ined.

\item[3)] For examples of simple Lie superalgebras without Cartan matrix
but with a $p|2p$-structure, see~\cite{LCh}. In addition to the
expected examples of the modular versions of Lie superalgebras of
vector f\/ields, and the queer analog of the $\fgl$ series, there are
-- {\bf for $p=2$} --  numerous (and hitherto unexpected) queerif\/ications, see
\cite{LCh}.
\end{enumerate}
\end{Remarks}

\subsubsection{$\boldsymbol{(2,4)}$-structure on Lie algebras} If $p=2$, we encounter a
new phenomenon f\/irst mentioned in \cite{KL}. Namely, let
$\fg=\fg_+\oplus \fg_-$ be a $\Zee/2$-grading of a Lie algebra. We
say that $\fg$ has a {\it $(2, -)$-structure},\index{$(2,
-)$-structure} if there is a $2$-structure on~$\fg_+$ but not on~$\fg$. It sometimes happens that this $(2, -)$-structure can be
extended to a\footnote{Observe a slightly dif\/ferent notation:
$(2,4)$, not $2|4$.}~{\it $(2,4)$-structure}, which means
that\index{$(2,4)$-structure}
\begin{equation}\label{2,4str}\text{~~ for any $x\in
\fg_-$ there exists an $x^{[4]}\in \fg_+$ such that $\ad_x^4=\ad_{
x^{[4]}}$.}
\end{equation}
For example, if indecomposable symmetrizable matrix $A$ is such that
($\fg(A)=\fo^{(1)}(2n+1)$)
\[
A_{11}=\od;\qquad A_{ii}=\ev\text{~for~} i>1,
\]
and the Lie algebra $\fg(A)$ (i.e., $\fg(A,(\ev,\dots,\ev))$ is
f\/inite-dimensional, then $\fg(A)$ has no $2$-structure but has a
$(2,4)$-structure inherited from the Lie {\bf super}algebra
$\fg(A,(\od,\ev,\dots,\ev))$.

\paragraph{$\boldsymbol{(2,4|2)}$-structure on Lie superalgebras.} A generalization of
the $(2,4)$-structure from Lie algebras to Lie superalgebras (such
as $\fo\fo_{\rm I\Pi}^{(1)}(2n+1|2m)$) is natural: Def\/ine the
$\Zee/2$-grading $\fg=\fg_+\oplus \fg_-$ of a Lie superalgebra
having nothing to do with the parity similarly to that of
$\fg(A)=\fo^{(1)}(2n+1)$, and def\/ine the squaring on the {\it plus}
part and a $(2,4)$-structure on the {\it minus} part such that the
conditions
\begin{equation}\label{2,4|2str}
\renewcommand{\arraystretch}{1.4}
\begin{array}{l}\ad_{x^{[2]}}(y)=\ad_x^2(y)\text{~ for all $x\in
(\fo\fo_{\rm I\Pi}^{(1)}(2k_\ev+1|2k_\od)_\ev)_+$,}\\
\ad_{x^{[4]}}(y)=\ad_x^4(y)\text{~ for all $x\in
(\fo\fo_{\rm I\Pi}^{(1)}(2k_\ev+1|2k_\od)_\ev)_-$} \end{array}
\end{equation}
are satisf\/ied for any $y\in\fo\fo_{\rm I\Pi}^{(1)}(2k_\ev+1|2k_\od)$,
not only for $y\in\fo\fo_{\rm I\Pi}^{(1)}(2k_\ev+1|2k_\od)_\ev$.

\paragraph{$\boldsymbol{(2|2)}$-structure on Lie superalgebras.} Lebedev observed
that if $p=2$ and a Lie superalgebra $\fg$ possesses a
$2|4$-structure, then the Lie algebra $F(\fg)$ one gets from $\fg$
by forgetting the superstructure (this is possible since
$[x,x]=2x^2=0$ for any odd $x$) possesses a $2$-structure given by
\begin{enumerate}\itemsep=0pt
\item[]
the ``2'' part of $2|4$-structure on the former $\fg_\ev$;

\item[] the squaring on $\fg_\od$;

\item[] the rule $(x+y)^{[2]}=x^{[2]}+y^{[2]}+[x,y]$ on the formerly
inhomogeneous (with respect to parity) elements.
\end{enumerate}

So one can say that if $p=2$, then any restricted Lie superalgebra
$\fg$ (i.e., the one with a~$2|4$-structure) induces a {\it
$2|2$-structure}\index{$2\vert 2$-structure} on the Lie algebra
$F(\fg)$ which is def\/ined even on inhomogeneous elements (unlike
$p|2p$-structures def\/ined on homogeneous elements only).

\parbegin{Remark} The restricted Lie superalgebra structures resemble (somehow) a
hidden supersymmetry of the following well-known fact:
\begin{equation}\label{nosusy}
\begin{minipage}[l]{12cm}
{\sl The product of two vector fields is not necessarily a vector
field, whereas their commutator always is a vector field. }
\end{minipage}
\end{equation}
This fact was not considered to be supersymmetric until recently:
Dzhumadildaev investigated a similar phenomenon: For the general and
divergence-free Lie algebras of polynomial vector f\/ields in $n$
indeterminates over $\Cee$, he investigated for which $N=N(n)$ the
anti-symmetrization of the map $D\longmapsto D^N$ (i.e., the
expression $\mathop{\sum}\limits_{\sigma\in S_N} \sign(\sigma) \,
X_{\sigma (1)}\dots X_{\sigma (N)}$) yields a vector f\/ield. For the
answer for $n=2, 3$ and a conjecture, see \cite{Dz}. But the most
remarkable is Dzhumadildaev's discovery of a hidden supersymmetry of
the usual commutator described by a universal odd vector f\/ield.
Dzhumadildaev deduced the above fact \eqref{nosusy} from the
following property of odd vector f\/ields:
\begin{equation}\label{susy}
\begin{minipage}[l]{12cm}
{\sl The product of two vector fields is not necessarily a vector
field, whereas the square of any odd field always is a vector field.
}
\end{minipage}
\end{equation}
\end{Remark}

\section{Ortho-orthogonal and periplectic Lie superalgebras}\label{Soo}

In this section, $p=2$ and $\Kee$ is perfect. We also assume that
$n_\ev,n_\od>0$. Set $n:=n_\ev+n_\od$.

\subsection{Non-degenerate even supersymmetric bilinear forms\\
 and
ortho-orthogonal Lie superalgebras}

For $p=2$, there are, in
general, four equivalence classes of inequivalent non-degenerate
even supersymmetric bilinear forms on a given superspace. Any such
form $B$ on a superspace $V$ of superdimension $n_\ev|n_\od$ can be
decomposed as follows:
\[
B=B_\ev\oplus B_\od,
\]
where $B_\ev$, $B_\od$ are symmetric non-degenerate forms on $V_\ev$
and $V_\od$, respectively. For $i=\ev,\od$, the form $B_i$ is
equivalent to $1_{n_i}$ if $n_i$ is odd, and equivalent to $1_{n_i}$
or $\Pi_{n_i}$ if $n_i$ is even. So every non-degenerate even
symmetric bilinear form is equivalent to one of the following forms
(some of them are def\/ined not for all dimensions):
\begin{alignat*}{3}
& B_{\rm II}=1_{n_\ev}\oplus 1_{n_\od}; \qquad & &B_{\rm I\Pi}=1_{n_\ev}\oplus
\Pi_{n_\od}\text{~if~}n_\od \text{~is even;}& \\
& B_{\rm \Pi I}=\Pi_{n_\ev}\oplus 1_{n_\od}\text{~if~}n_\ev \text{~is
even}; \qquad & &B_{\Pi\Pi}=\Pi_{n_\ev}\oplus \Pi_{n_\od}\text{~if~}n_\ev,
n_\od \text{~are even.}&
\end{alignat*}
We denote the Lie superalgebras that preserve the respective forms
by $\fo\fo_{\rm II}(n_\ev|n_\od)$, $\fo\fo_{\rm I\Pi}(n_\ev|n_\od)$,
$\fo\fo_{\rm \Pi I}(n_\ev|n_\od)$, $\fo\fo_{\Pi\Pi}(n_\ev|n_\od)$,
respectively. Now let us describe these algebras.

\subsubsection[$\fo\fo_{\rm II}(n_\ev|n_\od)$]{$\boldsymbol{\fo\fo_{\rm II}(n_\ev|n_\od)}$}

If $n:=n_\ev+n_\od\geq 3$, then the Lie
superalgebra $\fo\fo_{\rm II}^{(1)}(n_\ev|n_\od)$ is simple. This Lie
superalgebra has a $2|4$-structure; it {\bf has no Cartan matrix}.

\parbegin{Remark} To prove that a given Lie (super)algebra $\fg$ has no Cartan
matrix, we have to consider its maximal tori, and for each of them,
take the corresponding root grading. Then, if the simple roots are
impossible to def\/ine, or the elements of weight 0 do not commute,
etc.~-- if any of the requirements needed to def\/ine the Lie
(super)algebra with Cartan matrix is violated~-- we are done. We
skip such proofs in what follows.
\end{Remark}

\subsubsection[$\fo\fo_{\rm I\Pi}(n_\ev|n_\od)$ ($n_\od=2k_\od$)]{$\boldsymbol{\fo\fo_{\rm I\Pi}(n_\ev|n_\od)}$ ($\boldsymbol{n_\od=2k_\od}$)}

The Lie superalgebra $\fo\fo_{\rm I\Pi}(2k_\ev+1|2k_\od)$ possesses a
$2|4$-structure.

The Lie superalgebra $\fo\fo_{\rm I\Pi}^{(1)}(n_\ev|n_\od)$ is simple;
\[\fo\fo_{\rm I\Pi}^{(1)}(2k_\ev+1|2k_\od)\text{~~ possesses
$\begin{cases}\text{$2|4$-structure}&\text{if
$n_\ev=1$ ($k_\ev=0$)},\\
\text{$(2,4|2)$-structure}&\text{if $n_\ev> 1$ ($k_\ev>0$)};
\end{cases}$}
\]
$\fo\fo_{\rm I\Pi}^{(1)}(n_\ev|n_\od)$ has a Cartan matrix if and only
if $n_\ev$ is odd; this matrix has the following form (up to a
format; all possible formats~-- corresponding to $\ast=0$ or
$\ast=\ev$~-- are described in Table Section~\ref{tbl} below):
\begin{equation}\label{oowith1}
\begin{pmatrix} \ddots&\ddots&\ddots&\vdots\\
\ddots&\ast&1&0\\
\ddots&1&\ast&1\\
\cdots&0&1&1\end{pmatrix}.
\end{equation}

\subsubsection[$\fo\fo_{\Pi\Pi}(n_\ev|n_\od)$ ($n_\ev=2k_\ev$, $n_\od=2k_\od$)]{$\boldsymbol{\fo\fo_{\Pi\Pi}(n_\ev|n_\od)}$ ($\boldsymbol{n_\ev=2k_\ev}$, $\boldsymbol{n_\od=2k_\od}$)}\label{oo-oo_PP}

If $n=n_\ev+n_\od\geq 6$, then
\begin{equation}
\label{oopipi}
\begin{split}
&\text{if $k_\ev+k_\od$ is odd, then the Lie superalgebra
$\fo\fo_{\Pi\Pi}^{(2)}(n_\ev|n_\od)$ is simple;} \cr &\text{if
$k_\ev+k_\od$ is even, then the Lie superalgebra
$\fo\fo_{\Pi\Pi}^{(2)}(n_\ev|n_\od)/\Kee 1_{n_\ev|n_\od}$ is
simple.}
\end{split}
\end{equation}

Each of these simple Lie superalgebras has a $2|4$-structure; it is
also close to a Lie superalgebra with Cartan matrix. To describe
this CM Lie superalgebra in most simple terms, we will choose a
slightly dif\/ferent realization of $\fo\fo_{\Pi\Pi}(2k_\ev|2k_\od)$:
Let us consider it as the algebra of linear transformations that
preserve the bilinear form $\Pi(2k_\ev+2k_\od)$ in the format
$k_\ev|k_\od|k_\ev|k_\od$. Then the algebra
$\fo\fo_{\Pi\Pi}^{(i)}(2k_\ev|2k_\od)$ is spanned by supermatrices
of format $k_\ev|k_\od|k_\ev|k_\od$ and the form
\begin{equation}\label{matform}
\begin{pmatrix}A&C\\D&A^T\end{pmatrix},\quad \text{~where~}
\begin{array}{l}
A\in\begin{cases}\fgl(k_\ev|k_\od)&\text{if~}i\leq
1,\\\fsl(k_\ev|k_\od)&\text{if~}i\geq 2,\end{cases} \\ 
C,D\text{~are~}\begin{cases}\text{symmetric matrices}&\text{if~}
i=0,\\\text{symmetric zero-diagonal matrices}&\text{if~} i\geq
1.\end{cases}\end{array}
\end{equation}
If $i\geq 1$, these derived algebras have a non-trivial central
extension given by the following cocycle:
\begin{equation}\label{cocycle}
F\left(\begin{pmatrix}A&C\\D&A^T\end{pmatrix},
\begin{pmatrix}A'&C'\\D'&A'^T\end{pmatrix}\right)=\sum\limits_{1\leq
i<j\leq k_\ev+k_\od} (C_{ij}D'_{ij}+C'_{ij}D_{ij})
\end{equation}
(note that this expression resembles $\frac 12\tr(CD'+C'D)$). We
will denote this central extension of
$\fo\fo_{\Pi\Pi}^{(i)}(2k_\ev|2k_\od)$ by
$\fo\fo\fc(i,2k_\ev|2k_\od)$.

Let\index{$I_0:=\diag(1_{k_\ev\vert k_\od},0_{k_\ev\vert k_\od})$}
\begin{equation}
\label{I_0osp}I_0:=\diag(1_{k_\ev|k_\od},0_{k_\ev|k_\od}).
\end{equation}
Then the corresponding CM Lie superalgebra is
\begin{equation}
\label{ooc}
\begin{split}
&\fo\fo\fc(2,2k_\ev|2k_\od)\subplus\Kee I_0\text{ if $k_\ev+k_\od$
is odd;} \cr &\fo\fo\fc(1,2k_\ev|2k_\od)\subplus\Kee I_0\text{ if
$k_\ev+k_\od$ is even.}
\end{split}
\end{equation}

The corresponding Cartan matrix has the form (up to format; all
possible formats~-- corresponding to $\ast=0$ or $\ast=\ev$~-- are
described in Table Section~\ref{tbl} below):
\begin{equation}\label{ooPPCM}
\begin{pmatrix}
\ddots&\ddots&\ddots&\vdots&\vdots\\
\ddots&\ast&1&0&0\\
\ddots&1&\ast&1&1\\
\cdots&0&1&\ev&0\\
\cdots&0&1&0&\ev\end{pmatrix}.
\end{equation}

\subsection{The non-degenerate odd supersymmetric bilinear forms.\\
Periplectic Lie superalgebras} \label{peLS}

In this subsection,
$m\geq 3$.
\begin{equation}
\label{pe0}
\begin{split}
&\text{If $m$ is odd, then the Lie superalgebra $\fpe_B^{(2)}(m)$ is
simple;} \cr &\text{If $m$ is even, then the Lie superalgebra
$\fpe_B^{(2)}(m)/\Kee 1_{m|m}$ is simple.}
\end{split}
\end{equation}

If we choose the form $B$ to be $\Pi_{m|m}$, then the algebras
$\fpe_B^{(i)}(m)$ consist of matrices of the form (\ref{matform});
the only dif\/ference from $\fo\fo_{\Pi\Pi}^{(i)}$ is the format which
in this case is $m|m$.

Each of these simple Lie superalgebras has a $2$-structure. Note
that if $p\neq 2$, then the Lie superalgebra $\fpe_B(m)$ and its
derived algebras are not close to CM Lie superalgebras (because, for
example, their root system is not symmetric). If $p=2$ and $m\geq
3$, then they {\bf are} close to CM Lie superalgebras; here we
describe them.

The algebras $\fpe_B^{(i)}(m)$, where $i>0$, have non-trivial
central extensions with cocycles (\ref{cocycle}); we denote these
central extensions by $\fpe\fc(i,m)$. Let us introduce one more
matrix \index{$I_0:=\diag(1_m,0_m)$}
\begin{equation}
\label{I_0pe}I_0:=\diag(1_m,0_m).
\end{equation} Then the CM Lie superalgebras are
\begin{equation}
\label{pec}
\begin{split}
&\fpe\fc(2,m)\subplus\Kee I_0\text{ if $m$ is odd;} \cr
&\fpe\fc(1,m)\subplus\Kee I_0\text{ if $m$ is even.}
\end{split}
\end{equation}

The corresponding Cartan matrix has the form (\ref{ooPPCM}); the
only condition on its format is that the last two simple roots must
have distinct parities. The corresponding Dynkin diagram is shown in
Table Section~\ref{tbl}; all its nodes, except for the ``horns'', may
be both $\otimes$ or~$\odot$, see (\ref{cm1}).

\subsection{Superdimensions}

The following expressions (with a $+$ sign)
are the superdimensions of the
 relatives of the ortho-orthogonal and
periplectic Lie superalgebras that possess Cartan matrices. To get
the superdimensions of the simple relatives, one should replace $+2$
and $+1$ by $-2$ and $-1$, respectively, in the two f\/irst lines and
the four last ones:
\begin{alignat}{3}
& \dim \fo\fc (1;2k)\subplus\Kee I_0 =2k^2-k\pm 2\qquad &&\text{if $k$ is even;}& \nonumber\\
& \dim \fo\fc (2;2k)\subplus\Kee I_0=2k^2-k\pm 1\qquad & &\text{if $k$ is odd;}& \nonumber\\
& \dim\fo^{(1)}(2k+1)=2k^2+k; &&& \nonumber\\
& \sdim\fo\fo^{(1)}(2k_\ev+1|2k_\od)=2k_\ev^2+k_\ev+ 2k_\od^2+k_\od\mid
2k_\od(2k_\ev+1); &&& \label{dimtable}\\
& \sdim \fo\fo\fc (1;2k_\ev|2k_\od)\subplus\Kee I_0=2k_\ev^2-k_\ev+
2k_\od^2-k_\od\pm 2\mid 4k_\ev k_\od\qquad &&\text{if $k_\ev+k_\od$ is even;}&\nonumber\\
& \sdim \fo\fo\fc (2;2k_\ev|2k_\od)\subplus\Kee I_0=2k_\ev^2-k_\ev+
2k_\od^2-k_\od\pm 1\mid 4k_\ev k_\od\qquad &&\text{if $k_\ev+k_\od$ is odd;}& \nonumber\\
& \sdim \fpe\fc (1;m)\subplus\Kee I_0=m^2\pm 2\mid m^2-m\qquad &&\text{if $m$ is even;}& \nonumber\\
& \sdim \fpe\fc (2;m)\subplus\Kee I_0=m^2\pm 1\mid m^2-m\qquad &&\text{if $m$
is odd}.&\nonumber
\end{alignat}

\subsubsection{Summary: The types of Lie superalgebras preserving\\
non-degenerate symmetric forms}

Let
\begin{equation}\label{hat}\widehat{\fg}:=\fg\subplus\Kee I_0.\end{equation}
We have the following types of non-isomorphic Lie (super) algebras
(except for an occasional isomorphism intermixing the types, e.g.,
$\foo_{\Pi\Pi}^{(1)}(6|2)\simeq\fpe^{(1)}(4)$):
\begin{equation}\label{oandoo}\renewcommand{\arraystretch}{1.4}
\begin{tabular}{|l|l|}
\hline no relative has Cartan matrix&with Cartan matrix\\
\hline
$\foo_{\rm II}(2n+1|2m+1),\;\;\foo_{\rm II}(2n+1|2m)$&$\widehat{\fo\fc(i;2n)},\;\;\fo^{(1)}(2n+1);\;
\;\widehat{\fpe\fc(i;k)}$\\
$\foo_{\rm II}(2n|2m),
\;\;\foo_{\rm I\Pi}(2n|2m);\;\;\fo_{\rm I}(2n);$&$\widehat{\foo\fc(i;2n|2m)},
\;\;\foo_{\rm I\Pi}^{(1)}(2n+1|2m)$ \\
\hline\end{tabular}
\end{equation}
The superdimensions are as follows (in the second and third column
stand the additions to the superdimensions in the f\/irst column):
\begin{equation}\label{sdimtable}
\begin{tabular}{|l|l|l|}
\hline $\sdim\fo\fo_{\rm II}(a|b)$&$\sdim\fo\fo^{(1)}_{\rm II}(a|b)$&
$\sdim\fo\fo^{(2)}_{\rm II}(a|b)$\\
\hline$\frac12a(a+1)+ \frac12b(b+1)\mid
ab$&$-1|0$&\\
\hline$\sdim\fo\fo_{\rm I\Pi}(a|b)$&$\sdim\fo\fo^{(1)}_{\rm I\Pi}(a|b)$&
$\sdim\fo\fo^{(2)}_{\rm I\Pi}(a|b)$\\
\hline$\frac12a(a+1)+ \frac12b(b+1)\mid
ab$&$-a|0$&\\
\hline$\sdim\fo\fo_{\Pi\Pi}(a|b)$&$\sdim\fo\fo^{(1)}_{\Pi\Pi}(a|b)$&
$\sdim\fo\fo^{(2)}_{\Pi\Pi}(a|b)$\\
$\frac12a(a+1)+ \frac12b(b+1)\mid
ab$&$-a-b|0$&$-1|0$\\
\hline\end{tabular}
\end{equation}

\section{Dynkin diagrams}

A usual way to represent simple Lie algebras over $\Cee$ with
integer Cartan matrices is via graphs called, in the f\/inite
dimensional case, {\it Dynkin diagrams} (DD). The Cartan matrices of
certain interesting inf\/inite dimensional simple Lie {\it
super}algebras $\fg$ (even over $\Cee$) can be non-symmetrizable or
(for any $p$ in the super case and for $p>0$ in the non-super case)
have entries belonging to the ground f\/ield $\Kee$. Still, it is
always possible to assign an analog of the Dynkin diagram to each
(modular) Lie (super)algebra (with Cartan matrix, of course)
provided the edges and nodes of the graph (DD) are rigged with an
extra information. Although these analogs of the Dynkin graphs are
not uniquely recovered from the Cartan matrix (and the other way
round), they give a graphic presentation of the Cartan matrices and
help to observe some hidden symmetries.

Namely, the {\it Dynkin diagram}\index{Dynkin diagram} of a
normalized $n\times n$ Cartan matrix $A$ is a set of $n$ nodes
connected by multiple edges, perhaps endowed with an arrow,
according to the usual rules~\cite{K} or their modif\/ication, most
naturally formulated by Serganova: compare \cite{Se, FLS} with
\cite{FSS}. In what follows, we recall these rules, and further
improve them to f\/it the modular case.

\subsection{Nodes}

To every simple root there corresponds
\begin{equation}\label{cm1}
\begin{cases}
\text{a node}\; \mcirc\; &\text{if $p(\alpha_{i})= \ev$ and $
A_{ii}=2$},\\
\text{a node}\; \ast \;&\text{if $p(\alpha_{i}) =\ev$ and
$A_{ii}=\od$};\\
\text{a node}\; \mbullet \;&\text{if $p(\alpha_{i}) =\od$ and
$A_{ii}=1$};\\
\text{a node}\; \motimes \;& \text{if $p(\alpha_{i})
=\od$ and $
A_{ii}=0$},\\
\text{a node}\; \odot\; &\text{if $p(\alpha_{i})= \ev$ and $
A_{ii}=\ev$}.
\end{cases}
\end{equation}

The Lie algebras $\fsl(2)$ and $\fo(3)^{(1)}$ with Cartan matrices
$(2)$ and $(\od)$, respectively, and the Lie superal\-geb\-ra
$\fosp(1|2)$ (which is $\foo_{\rm I\Pi}^{(1)}(1|2)$ if $p=2$) with
Cartan matrix $(1)$ are simple.

The Lie algebra with Cartan matrix $(\ev)$ and the Lie superalgebra
with Cartan matrix $(0)$ are solvable of $\dim 4$ and $\sdim 2|2$,
respectively. Their derived algebras are {\it Heisenberg al\-geb\-ra}
$\fhei(2)\simeq\fhei(2|0)$ and {\it Heisenberg superalgebra}
$\fhei(0|2)\simeq\fsl(1|1)$, respectively; their (super)dimensions
are 3 and $1|2$, respectively.

\subsubsection{Digression}\label{fock}

Let $\xi=(\xi_1, \dots, \xi_n)$ and
$\eta=(\eta_1, \dots, \eta_n)$ be odd elements, $p=(p_1, \dots,
p_m)$, $q=(q_1, \dots, q_m)$ and $z$ even elements.
$\fhei(2m|2n)=\Span(p, q, \xi, \eta, z)$, where the brackets are
\begin{equation}\label{heicr}{}[p_i, q_j]=\delta_{ij}z, \qquad
[\xi_i, \eta_j]=\delta_{ij}z,\qquad [z, \fhei(2m|2n)]=0.
\end{equation} In what follows we will need the Lie superalgebra
$\fhei(2m|2n)$ (for the cases where $mn=0$) and its only (up to the
change of parity) non-trivial irreducible representation, called the
{\it Fock space}, which in characteristic $p$ is $\Kee[q,
\xi]/(q_1^p, \dots, q_n^p)$ on which the elements $q_i$ and $\xi_j$
act as operators of left multiplication by $q_i$ and $\xi_j$,
respectively, whereas $p_i$ and $\eta_j$ act as $h\partial_{q_i}$
and $h\partial_{\xi_j}$, where $h\in\Kee\setminus\{0\}$ can be f\/ixed
to be equal to 1 by a change of the basis.

\sssbegin{Remark} {\it A posteriori} (from the classif\/ication of
simple Lie superalgebras with Cartan matrix and of polynomial growth
for $p=0$) we f\/ind out that the roots~$\odot$ can only occur if
$\fg(A, I)$ grows faster than polynomially. Thanks to classif\/ication
again, if $\dim \fg<\infty$, the roots of type $\odot$ can not occur
if $p>3$; whereas for $p=3$, the Brown Lie algebras are examples of
$\fg(A)$ with a simple root of type $\odot$; for $p=2$, such roots
are routine.
\end{Remark}

\subsection{Edges} If $p=2$ and $\dim \fg(A)<\infty$, the Cartan matrices
considered are symmetric. If $A_{ij}=a$, where $a\neq 0$ or 1, then
we rig the edge connecting the $i$th and $j$th nodes by a label $a$.

If $p>2$ and $\dim \fg(A)<\infty$, then $A$ is symmetrizable, so let
us symmetrize it, i.e., consider $DA$ for an invertible diagonal
matrix $D$. Then, if $(DA)_{ij}=a$, where $a\neq 0$ or $-1$, we rig
the edge connecting the $i$th and $j$th nodes by a label $a$.

If all of\/f-diagonal entries of $A$ belong to $\Zee/p$ and their
representatives are selected to be non-positive integers, we can
draw the DD as for $p=0$, i.e., connect the $i$th node with the
$j$th one by $\max(|A_{ij}|, |A_{ji}|)$ edges rigged with an arrow
$>$ pointing from the $i$th node to the $j$th if $|A_{ij}|>|A_{ji}|$
or in the opposite direction if $|A_{ij}|<|A_{ji}|$.

\subsection{Ref\/lections} Let $R^+$ be a system of positive roots of Lie
superalgebra $\fg$, and let $B=\{\sigma_1,\dots,\sigma_n\}$ be the
corresponding system of simple roots with some corresponding pair
$(A=A_B,I=I_B)$. Then for any $k\in \{1, \dots, n\}$, the set
$(R^+\backslash\{\sigma_k\})\coprod\{-\sigma_k\}$ is a system of
positive roots. This operation is called {\it the reflection in
$\sigma_k$}; it changes the system of simple roots by the formulas
\begin{equation}
\label{oddrefl}
r_{\sigma_k}(\sigma_{j})= \begin{cases}{-\sigma_j}&\text{if~}k=j,\\
\sigma_j+B_{kj}\sigma_k&\text{if~}k\neq j,\end{cases}\end{equation}
where
\begin{equation}
\label{Boddrefl}B_{kj}=\begin{cases}
-\displaystyle\frac{2A_{kj}}{A_{kk}}& \text{~if~}i_k=\ev, A_{kk}\neq
0,\text{~and~}
-\displaystyle\frac{2A_{kj}}{A_{kk}}\in \Zee/p\Zee,\\
p-1&\text{~if~}i_k=\ev, A_{kk}\neq 0\text{~and~}
 -\displaystyle\frac{2A_{kj}}{A_{kk}}\not\in \Zee/p\Zee,\\
-\displaystyle\frac{A_{kj}}{A_{kk}}&
\text{~if~}i_k=\od, A_{kk}\neq 0,\text{~and~}
-\displaystyle\frac{A_{kj}}{A_{kk}}\in \Zee/p\Zee,\\
p-1&\text{~if~}i_k=\od, A_{kk}\neq 0,
\text{~and~} -\displaystyle\frac{A_{kj}}{A_{kk}}\not\in \Zee/p\Zee,\\
1&\text{~if~}i_k=\od, A_{kk}=0,A_{kj}\neq 0,\\
0&\text{~if~}i_k=\od, A_{kk}=A_{kj}=0,\\
p-1&\text{~if~}i_k=\ev, A_{kk}=\ev,A_{kj}\neq 0,\\
0&\text{~if~}i_k=\ev, A_{kk}=\ev,A_{kj}=0,\end{cases}
\end{equation}
where we consider $\Zee/p\Zee$ as a subf\/ield of $\Kee$.

\sssbegin{Remark} In the second, fourth and penultimate cases, the
matrix entries in (\ref{Boddrefl}) can, in principle, be equal to
$kp-1$ for any $k\in\Nee$, and in the last case any element of
$\Kee$ may occur. We may only hope at this stage that, at least for
$\dim\fg<\infty$, this does not happen.\end{Remark}

The values $-\displaystyle\frac{2A_{kj}}{A_{kk}}$ and
$-\displaystyle\frac{A_{kj}}{A_{kk}}$ are elements of $\Kee$, while
the roots are elements of a vector space over $\Ree$. Therefore
\textit{These expressions in the first and third cases in}
(\ref{Boddrefl}) \textit{should be understood as} ``\textit{the minimal
non-negative integer congruent to
$-\displaystyle\frac{2A_{kj}}{A_{kk}}$ or
$-\displaystyle\frac{A_{kj}}{A_{kk}}$, respectively''. $($If
$\dim\fg<\infty$, these expressions are always congruent to
integers.$)$}

\textit{There is known just one exception: If $p=2$ and $A_{kk}=A_{jk}$,
the expression} $-\displaystyle\frac{2A_{jk}}{A_{kk}}$ \textit{should be
understood as $2$, not $0$.}

The name ``ref\/lection'' is used because in the case of
(semi)simple f\/inite-dimensional Lie algebras this action extended on
the whole $R$ by linearity is a map from $R$ to $R$, and it does not
depend on $R^+$, only on $\sigma_k$. This map is usually denoted by
$r_{\sigma_k}$ or just $r_{k}$. The map $r_{\sigma_i}$ extended to
the $\Ree$-span of $R$ is ref\/lection in the hyperplane orthogonal to
$\sigma_i$ relative the bilinear form dual to the Killing form.

The ref\/lections in the even (odd) roots are referred to as {\it
even} ({\it odd}) {\it reflections}.\index{Reflection!
odd}\index{Reflection! even! non-isotropic} \index{Reflection! even!
isotropic} A simple root is called {\it isotropic}, if the
corresponding row of the Cartan matrix has zero on the diagonal, and
{\it non-isotropic} otherwise. The ref\/lections that correspond to
isotropic or non-isotropic roots will be referred to accordingly.

If there are isotropic simple roots, the ref\/lections $r_\alpha$ do
not, as a rule, generate a version of the {\it Weyl group} because
the product of two ref\/lections in nodes not connected by one
(perhaps, multiple) edge is not def\/ined. These ref\/lections just
connect pair of ``neighboring'' systems of simple roots and
there is no reason to expect that we can multiply two distinct such
ref\/lections. In the general case (of Lie superalgebras and $p>0$),
the action of a given isotropic ref\/lections~(\ref{oddrefl}) can not,
generally, be extended to a linear map $R\tto R$. For Lie
superalgebras over $\Cee$, one can extend the action of ref\/lections
by linearity to the root lattice but this extension preserves the
root system only for $\fsl(m|n)$ and $\fosp(2m+1|2n)$, cf.~\cite{Se1}.

If $\sigma_i$ is an odd isotropic root, then the corresponding
ref\/lection sends one set of Chevalley generators into a new one:
\begin{equation}
\label{oddrefx} \tilde X_{i}^{\pm}=X_{i}^{\mp};\;\;
\tilde X_{j}^{\pm}=\begin{cases}[X_{i}^{\pm},
X_{j}^{\pm}]&\text{if $A_{ij}\neq 0, \ev$},\\
X_{j}^{\pm}&\text{otherwise}.\end{cases}
\end{equation}

\subsubsection{Lebedev's lemma} Serganova \cite{Se}
proved (for $p=0$) that there is always a chain of
ref\/lections connecting $B_1$ with some system of simple roots $B'_2$
equivalent to $B_2$ in the sense of def\/inition~\ref{EqSSR}. Here is
the modular version of Serganova's Lemma. Observe that Serganova's
statement is not weaker: Serganova used only odd ref\/lections.

\begin{Lemma}[\cite{LCh}]\label{serg} For any two systems of simple roots
$B_1$ and $B_2$ of any finite dimensional Lie superalgebra with
indecomposable Cartan matrix, there is always a chain of reflections
connecting~$B_1$ with $B_2$.\end{Lemma}

\section{A careful study of an example}\label{Sex}
Now let $p=2$ and let us apply all the above to the Lie superalgebra
$\fpe(k)$ (the situation with $\fo_\Pi(2k)$ and
$\fo\fo_{\Pi\Pi}(2k_\ev|2k_\od)$ is the same). For the Cartan matrix
(all possible formats~-- corresponding to $\ast=0$ or $\ast=\ev$~-- are listed in Table Section~\ref{tbl}) we take
\begin{equation}
\label{peCM} A=\begin{pmatrix}
\ddots & \ddots & \ddots & \ddots \\
\cdots & \ast& 1 & 1\\
\cdots & 1& 0& 0\\
\cdots& 1 & 0 &\ev
\end{pmatrix}.
\end{equation}

The Lie superalgebra $\fpe^{(i)}(k)$ consists of supermatrices of
the form
\[
\begin{pmatrix} B&C\\D&B^T\end{pmatrix},
\]
where
\begin{equation}\label{pe1}
\begin{tabular}{ll}
for $i=0$, we have &$B\in\fgl(k)$, $C$, $D$ are symmetric;\\
for $i=1$, we have &$B\in\fgl(k)$, $C$, $D$ are symmetric zero-diagonal;\\
for $i=2$, we have &$B\in\fsl(k)$, $C$, $D$ are symmetric
zero-diagonal.
\end{tabular}
\end{equation}

We expect (by analogy with the orthogonal Lie algebras in
characteristic $\neq 2$) that
\begin{equation}\label{pe2}
\begin{tabular}{l}
$e_i^+=E_{i,i+1}+E_{k+i+1,k+i};\qquad
e_i^-=E_{i+1,i}+E_{k+i,k+i+1} $\quad for $i=1,\dots, k-1$;\\
$ e_k^+=E_{k-1,2k}+E_{k,2k-1};\qquad e_k^-=E_{2k-1,k}+E_{2k,k-1} $.
\end{tabular}
\end{equation}

Let us f\/irst consider the (simpler) \underline{case of $k$ odd}.
Then $\rk A=k-1$ since the sum of the last two rows is zero. Let us
start with the simple algebra $\fpe^{(2)}(k)$. The Cartan subalgebra
(i.e., the subalgebra of diagonal matrices) is $(k-1)$-dimensional
because the elements $[e_1^+,e_1^-],\dots, [e_{k-1}^+,e_{k-1}^-]$
are linearly independent, whereas
$[e_k^+,e_k^-]=[e_{k-1}^+,e_{k-1}^-]$. Thus, we should f\/irst f\/ind a
non-trivial central extension, spanned by $z$ satisfying the
condition
\begin{equation}\label{pe6}
z=[e_k^+,e_k^-]+[e_{k-1}^+,e_{k-1}^-].
\end{equation}
Elucidation: The values of $e_i^\pm$ in (\ref{pe2}) are what we {\it
expect} them to be from their $p=0$ analogs. But from the def\/inition
of CM Lie superalgebra we see that the algebra must have a center
$z$, see~(\ref{pe6}). Thus, the CM Lie superalgebra is not
$\fpe^{(2)}(k)$ but is spanned by the central extension of~$\fpe^{(2)}(k)$ plus the grading operator def\/ined from
(\ref{central1}). The extension $\fpe\fc(2,k)$ described in~(\ref{cocycle}) satisf\/ies this condition.

Now let us choose $B$ to be $(0,\dots,0,1)$. Then we need to add to
the algebra a grading operator $d$ such that
\begin{equation}\label{gr1}
\begin{array}{l}
{}[d,e_i^\pm]=0 \text{~for all~}i=1,\dots, k-1;\\
{}[d,e_k^\pm]=e_k^\pm;\\
d\text{~commutes with all diagonal matrices.} \end{array}
\end{equation}
The matrix $I_0=\diag(1_k,0_k)$ satisf\/ies all these conditions.
Thus, the corresponding CM Lie superalgebra is
\begin{equation}\label{pec2}
\fpe\fc(2,k)\subplus\Kee I_0.
\end{equation}

\begin{Remark} Rather often we need ideals of CM Lie (super)algebras
that do not contain the outer grading operator(s), cf.\ Section~\ref{warn}. These ideals, such as $\fpe\fc(2,k)$ or $\fsl(n|n)$, do
not have Cartan matrix.
\end{Remark}

Now let us consider the \underline{case of $k$ even}. Then the
simple algebra is $\fpe^{(2)}(k)/(\Kee 1_{2k})$. The Cartan matrix
is of rank $k-2$:
\begin{equation}\label{pe3}
\begin{tabular}{l}
(a) the sum of
the last two rows is zero; \\
(b) the sum of all the rows with odd numbers is zero.
\end{tabular}
\end{equation}
The condition (\ref{pe3}a) gives us the same central extension and
the same grading operator an in the previous case.

To satisfy condition (\ref{pe3}b), we should f\/ind a non-trivial
central extension such that
\[
z=\sum\limits_{i\text{~is odd}} [e_i^+,e_i^-].
\]
(This formula follows from (\ref{central}) and the 2nd equality in
(\ref{pe3}).) But we can see that, in $\fpe^{(2)}(k)$, we have
\[
\sum\limits_{i\text{~is odd}} [e_i^+,e_i^-]= \sum\limits_{i\text{~is
odd}}(E_{i,i}+E_{i+1,i+1}+E_{k+i,k+i}+E_{k+i+1,k+i+1})=1_{2k}.
\]
It means that the corresponding central extension of
$\fpe^{(2)}(k)/(\Kee 1_{2k})$ is just $\fpe^{(2)}(k)$.

Now, concerning the grading operator: Let the second row of $B$ be
$(1,0,\dots,0)$ (the f\/irst row is, as in the previous case,
$(0,\dots,0,1)$). Then we need a grading operator $d_2$ such that
\begin{equation}\label{pe4}
\begin{tabular}{l}
$[d_2,e_1^\pm]=e_1^\pm;$\\ $[d_2,e_i^\pm]=0\text{~for
all~}i>1;$\\
$d_2\text{~commutes with all diagonal matrices.}$
\end{tabular}
\end{equation}
The matrix $d_2:=E_{1,1}+E_{k+1,k+1}$ satisf\/ies these conditions.
But $\fpe^{(2)}(k)\subplus \Kee(E_{1,1}+E_{k+1,k+1})$ is just
$\fpe^{(1)}(k)$. So, the resulting CM Lie superalgebra is
\[
\fpe\fc(1,k)\subplus\Kee I_0.
\]

\section{Main steps of our classif\/ication}\label{Ssteps}

In this section we deal with Lie (super)algebras of the form
$\fg(A)$ or their simple subquotients $\fg(A)^{(i)}/\fc$, where
$i=1$ or 2.

\subsection{Step 1: An overview of known results}\label{SSsteps}

\underline{Lie algebras (nothing super)}. There are known the two
methods of classif\/ication:

1) \underline{Over $\Cee$}, Cartan \cite{C} did not use any roots,
instead he used what is nowadays called in his honor {\it Cartan
prolongations} and a generalization (which he never formulated
explicitly) of this procedure which we call {\it CTS-ing} ({\it
Cartan--Tanaka--Shchepochkina prolonging}).

2) Nowadays, to get the shortest classif\/ication of the simple f\/inite
dimensional Lie algebras, everybody (e.g.~\cite{Bou, OV}) uses root
technique and the non-degenerate invariant symmetric bilinear form
(the Killing form).

In the {\bf modular} case, as well as in the {\bf super} case, and
in the mixture of these cases we consider here, the Killing form
might be identically zero. However, if the Cartan matrix $A$ is
symmetrizable (and indecomposable), on the Lie (super)algebra
$\fg(A)$ if $\fg(A)$ is simple (or on $\fg(A)^{(i)}/\fc$ if $\fg(A)$
is not simple), there is a non-degenerate replacement of the Killing
form. (Astonishingly, this replacement might sometimes be not coming
from any representation, see \cite{Ser}. Much earlier Kaplansky
observed a similar phenomenon in the modular case and associated the
non-degenerate bilinear form with a {\it projective} representation.
Kaplansky pointed at this phenomenon in his wonderful preprints~\cite{Kapp} which he modestly did not publish.)

In the modular case, and in the super case for $p=0$, this approach~-- to use a non-degenerate even invariant symmetric form in
order to classify the simple algebras~-- was pursued by Kaplansky~\cite{Kapp}.

\underline{For $p>0$}, Weisfeiler and Kac \cite{WK} gave a
classif\/ication, but although the idea of their proof is OK, the
paper has several gaps and vague notions (the Brown algebra
$\fbr(3)$ was missed, whereas Brown \cite{Br3} who discovered it did
not write that it possesses Cartan matrix, actually two inequivalent
matrices f\/irst observed by Skryabin \cite{Sk1, KWK}; the notion of
the Lie algebra with Cartan matrix nicely formulated in \cite{K} was
not properly developed at the time \cite{WK} was written; the Dynkin
diagrams mentioned there were not def\/ined at all in the modular
case; the algebras $\fg(A)$ and $\fg(A)^{(i)}/\fc$ were sometimes
identif\/ied). The case $p>3$ being completely investigated by Block,
Wilson, Premet and Strade \cite{PS, S} (see also \cite{BGP}), we
double-checked the cases where $p<5$. The answer of
\cite{WK}$\cup$\cite{Sk1} is correct.

\medskip

\underline{Lie {\bf super}algebras}.

\underline{Over $\Cee$}, for {\it any} Lie algebra $\fg_\ev$, Kac
\cite{K2} listed all
\begin{equation}\label{kacg(a)}
\text{$\fg_\ev$-modules $\fg_\od$ such that the Lie superalgebra
$\fg=\fg_\ev\oplus \fg_\od$ is simple.} \end{equation}

Kaplansky \cite{Kapp, FK, Kap}, Djokovi\'c and Hochschild \cite{Dj},
and also Scheunert, Nahm and Rittenberg \cite{SNR} had their own
approaches to the problem (\ref{kacg(a)}) and solved it without gaps
for various particular cases, but they did not investigate which of
the simple f\/inite dimensional Lie superalgebras possess Cartan
matrix.

Kac observed that (a) some of the simple Lie superalgebras
(\ref{kacg(a)}) possess analogs of Cartan matrix, (b) one Lie
superalgebra may have several inequivalent Cartan matrices. His
f\/irst list of inequivalent Cartan matrices (in other words, distinct
$\Zee$-gradings) for f\/inite dimensional Lie superalgebras $\fg(A)$
in~\cite{K2} had gaps; Serganova~\cite{Se} and (by a dif\/ferent
method and only for symmetrizable matrices) van de Leur~\cite{vdL}
f\/ixed the gaps and even classif\/ied Lie superalgebras of polynomial
growth (for the proof in the non-symmetrizable case, announced 20
years earlier, see~\cite{HS}). Kac also suggested analogs of Dynkin
diagrams to graphically encode the Cartan matrices.

Kaplansky was the f\/irst (see his newsletters in \cite{Kapp}) to
discover the exceptional algebras $\fag(2)$ and $\fab(3)$ (he dubbed
them $\Gamma_2$ and $\Gamma_3$, respectively) and a parametric
family $\fosp(4|2; \alpha)$ (he dubbed it $\Gamma(A, B, C))$); our
notations ref\/lect the fact that $\fag(2)_\ev=\fsl(2)\oplus\fg(2)$
and $\fab(3)_\ev=\fsl(2)\oplus\fo(7)$ ($\fo(7)$ is $B_3$ in Cartan's
nomenclature). Kaplansky's description (irrelevant to us at the
moment except for the fact that $A$, $B$ and $C$ are on equal
footing) of what we now identify as $\fosp(4|2; \alpha)$, a
parametric family of deforms of $\fosp(4|2)$, made
 an $S_3$-symmetry of the parameter manifest (to A.\ A.~Kirillov,
 and he informed us, in 1976).
Indeed, since $A+B+C=0$, and $\alpha\in \Cee\cup\infty$ is the ratio
of the two remaining parameters, we get an $S_3$-action on the plane
$A+B+C=0$ which in terms of $\alpha$ is generated by the
transformations:
\begin{equation}\label{osp42symm}
\alpha\longmapsto -1-\alpha, \qquad \alpha\longmapsto
\frac{1}{\alpha}.
\end{equation}
This symmetry should have immediately sprang to mind since
$\fosp(4|2; \alpha)$ is strikingly similar to $\fwk(3; a)$ found 5
years earlier, cf.~(\ref{wkiso}), and since $S_3\simeq \SL(2;
\Zee/2)$.

The following f\/igure depicts the fundamental domains of the
$S_3$-action. The other transformations generated by
(\ref{osp42symm}) are \[\alpha\longmapsto
-\frac{1+\alpha}{\alpha},\qquad\alpha\longmapsto
-\frac{1}{\alpha+1},\qquad\alpha\longmapsto -\frac{\alpha}{\alpha+1}.
\]
\begin{figure}[ht]\centering
\includegraphics[scale=1.05]{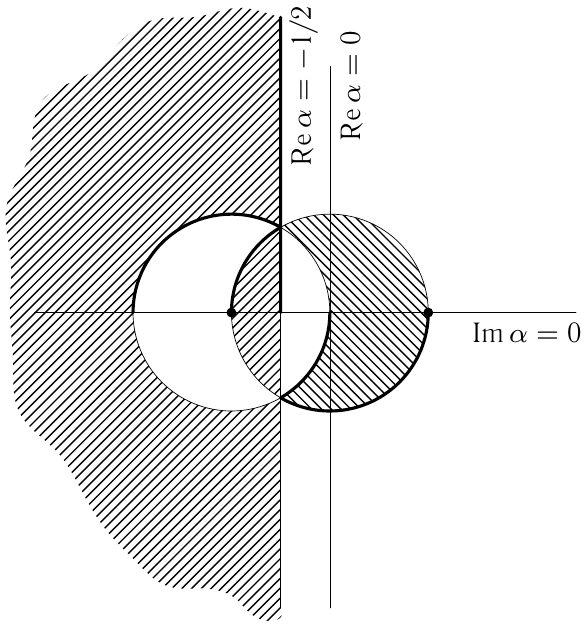}
\end{figure}

\subsubsection{Notation: On matrices with a ``--'' sign and other
notations\\ in the lists of inequivalent Cartan
matrices}\label{recmat}

The rectangular matrix at the beginning of
each list of inequivalent Cartan matrices for each Lie superalgebra
shows the result of odd ref\/lections (the number of the row is the
number of the Cartan matrix in the list below, the number of the
column is the number of the root (given by small boxed number) in
which the ref\/lection is made; the cells contain the results of
ref\/lections (the number of the Cartan matrix obtained) or a ``--'' if the ref\/lection is not appropriate because $A_{ii}\neq 0$.
Some of the Cartan matrices thus obtained are equivalent, as
indicated.

The number of the matrix $A$ such that $\fg(A)$ has only one odd
simple root is \boxed{boxed}, that with all simple roots odd is
\underline{underlined}. The nodes are numbered by small boxed
numbers; the curly lines with arrows depict odd ref\/lections.

Recall that $\fag(2)$ of $\sdim = 17|14$ has the following Cartan
matrices
\begin{figure}[ht!]\centering
\begin{minipage}[m]{0.2\linewidth}\centering
$
\begin{lmatrix}
 2 & - & - \\
 1 & 3 & - \\
 - & 2 & 4 \\
 - & - & 3
\end{lmatrix}
$
\end{minipage}
\begin{minipage}[m]{0.48\linewidth}\centering
\includegraphics{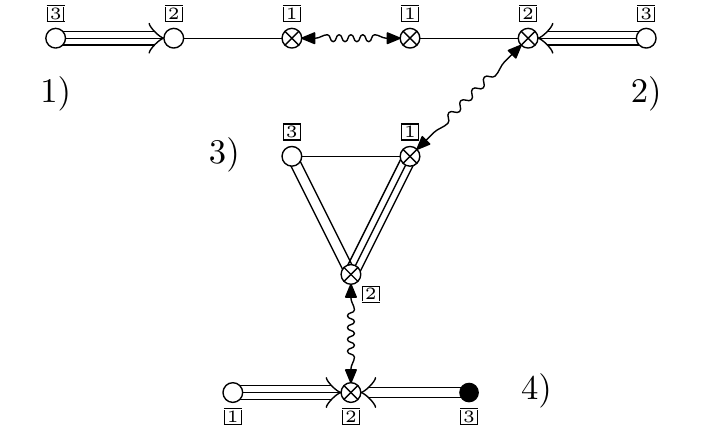}
\end{minipage}\hfill
\end{figure}
\begin{equation}\label{ag2cm}
\boxed{1)}\; \begin{pmatrix} 0 & -1 & 0 \\ -1 & 2 & -3 \\ 0 & -1 & 2
\end{pmatrix},\quad 2)\; \begin{pmatrix}
0 & -1 & 0 \\ -1 & 0 & 3 \\ 0 & -1 & 2
\end{pmatrix},\quad 3)\; \begin{pmatrix}
0 & -3 & 1 \\ -3 & 0 & 2 \\ -1 & -2 & 2
\end{pmatrix},\quad 4)\;
\begin{pmatrix}
2 & -1 & 0 \\ -3 & 0 & 2 \\ 0 & -1 & 1
\end{pmatrix}.
\end{equation}

Recall that $\fab(3)$ of $\sdim = 24|16$ has the following Cartan
matrices
\begin{figure}[ht]\centering
\begin{minipage}[m]{0.2\linewidth}\centering
$
\begin{lmatrix}
 - & 2 & - & - \\
 3 & 1 & 4 & - \\
 2 & - & - & - \\
 - & - & 2 & 5 \\
 - & 6 & - & 4 \\
 - & 5 & - & -
\end{lmatrix}
$
\end{minipage}
\begin{minipage}[m]{0.58\linewidth}\centering
\includegraphics{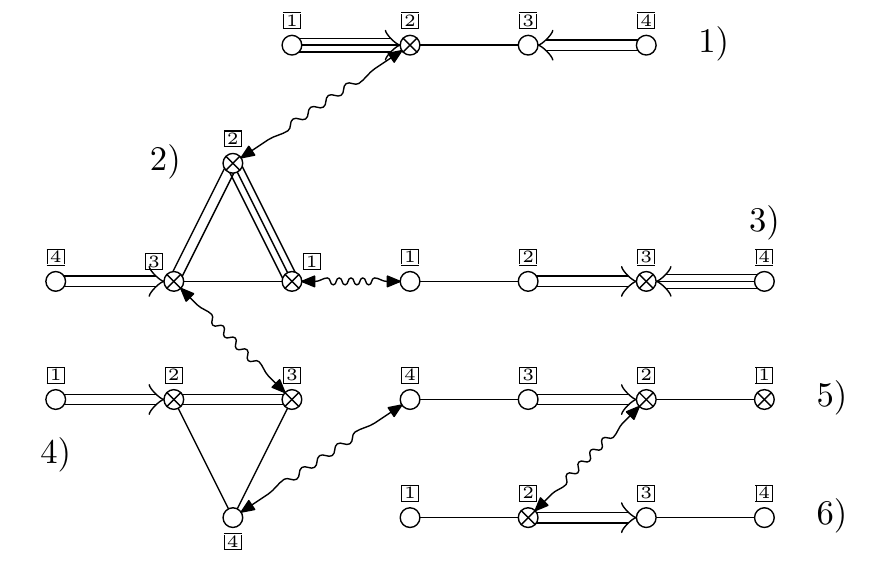}
\end{minipage}\hfill
\end{figure}
\begin{equation}\label{ab3cm}
\begin{matrix}
\boxed{1)}\; \begin{pmatrix}
2 & -1 & 0 & 0 \\ -3 & 0 & 1 & 0 \\ 0 & -1 & 2 & -2 \\
 0 & 0 & -1 & 2
\end{pmatrix},\quad 2)\;
\begin{pmatrix}
0 & -3 & 1 & 0 \\ -3 & 0 & 2 & 0 \\ 1 & 2 & 0 & -2 \\
0 & 0 & -1 & 2
\end{pmatrix},\quad \boxed{3)}\;
\begin{pmatrix}
2 & -1 & 0 & 0 \\ -1 & 2 & -1 & 0 \\ 0 & -2 & 0 & 3 \\
0 & 0 & -1 & 2
\end{pmatrix},\\
4)\; \begin{pmatrix}
2 & -1 & 0 & 0 \\ -2 & 0 & 2 & -1 \\ 0 & 2 & 0 & -1 \\
0 & -1 & -1 & 2
\end{pmatrix},\quad 5)\;
\begin{pmatrix}
0 & 1 & 0 & 0 \\ -1 & 0 & 2 & 0 \\ 0 & -1 & 2 & -1 \\
0 & 0 & -1 & 2
\end{pmatrix},\quad \boxed{6)}\;
\begin{pmatrix}
2 & -1 & 0 & 0 \\ -1 & 2 & -1 & 0 \\ 0 & -2 & 2 & -1 \\ 0 & 0 & -1 &
0
\end{pmatrix}.
\end{matrix}\!\!\!\!\!
\end{equation}

\underline{Modular Lie algebras and Lie superalgebras}.

\underline{$p=2$, Lie algebras}. Weisfeiler and Kac \cite{WK}
discovered two new parametric families that we denote $\fwk(3;a)$
and $\fwk(4;a)$ ({\it Weisfeiler and Kac
algebras}).\index{$\fwk(3;a)$, Weisfeiler and Kac algebra}
\index{$\fwk(3;a)$, Weisfeiler and Kac algebra}\index{Weisfeiler and
Kac algebra}

$\fwk(3;a)$, where $a\neq 0, -1$, of dim 18 is a non-super version
of $\fosp(4|2; a)$ (although no $\fosp$ exists for $p=2$); the
dimension of its simple subquotient $\fwk(3;a)^{(1)}/\fc$ is equal
to 16; the inequivalent Cartan matrices are:
\[
1)\; \begin{pmatrix} \ev &a &0\\
a&\overline{0}&1\\0&1&\overline{0} \end{pmatrix} ,\quad
2)\; \begin{pmatrix} \ev &1+a &a\\
1+a&\overline{0}&1\\
a&1&\overline{0} \end{pmatrix}.
\]

$\fwk(4;a)$, where $a\neq 0, -1$, of $\dim=34$; the inequivalent
Cartan matrices are:
\[
1)\; \begin{pmatrix} \ev &a &0&0\\
a &\overline{0}&1&0\\
0&1&\overline{0}&1\\
0&0&1&\overline{0} \end{pmatrix} ,\quad
2)\; \begin{pmatrix} \ev &1 &1+a&0\\
1 &\overline{0}& a & 0\\
a+1& a &\overline{0}&a\\
0&0&a&\overline{0} \end{pmatrix}
 ,\quad
3)\;\begin{pmatrix} \ev &a & 0 &0\\
a &\overline{0}& a+1 & 0\\
0& a+1 &\overline{0}&1\\
0&0&1&\overline{0} \end{pmatrix}.
\]

Weisfeiler and Kac investigated also which of these algebras are
isomorphic and the answer is as follows:
\begin{equation}\label{wkiso}
\renewcommand{\arraystretch}{1.4}\begin{array}{l}
\fwk(3;a)\simeq \fwk(3;a')\Longleftrightarrow
a'=\displaystyle\frac{\alpha a+\beta}{\gamma a+\delta}, \quad \text{where
$\begin{pmatrix}\alpha&\beta\\ \gamma &\delta\end{pmatrix}\in \SL(2; \Zee/2)$},\\
\fwk(4;a)\simeq \fwk(4;a')\Longleftrightarrow
a'=\displaystyle\frac{1}{a}.
\end{array}
\end{equation}

\subsubsection{2-structures on $\fwk$ algebras}\label{2srtwk}

1) Observe
that the center $\fc$ of $\fwk(3; a)$ is spanned by $a h_1 +h_3$.
The 2-structure on $\fwk(3;a)$ is given by the conditions
$(e_\alpha^{\pm})^{[2]}=0$ for all root vectors and the following
ones:

a) For the matrix $B=(0,0,1)$ in (\ref{matrixB}) for the grading
operator $d$, set:
\begin{equation}\label{2strwk3}
\begin{array}{l} (\ad_{h_1})^{[2]} =
(1+a t) h_1 + t h_3\equiv h_1\pmod{\fc},\\
(\ad_{h_2})^{[2]} = a t h_1 + h_2 + t h_3 + a(1 +a) d
\equiv h_2 + a(1 +a) d\pmod{\fc},\\
(\ad_{h_3})^{[2]} = (a t +a^2) h_1 + t h_3\equiv a^2h_1\pmod{\fc},\\
(\ad_d)^{[2]} = a t h_1 + t h_3 +d\equiv d\pmod{\fc},
\end{array}\end{equation} where $t$ is a parameter.

b) Taking $B=(1,0,0)$ in (\ref{matrixB}) we get a more symmetric
answer:
\begin{equation}\label{2strwk31}
\begin{array}{l} (\ad_{h_1})^{[2]} =
(1+a t) h_1 + t h_3\equiv h_1\pmod{\fc},\\
(\ad_{h_2})^{[2]} = a t h_1 + a h_2 + t h_3 + (1 +a) d
\equiv a h_2 + (1 +a) d\pmod{\fc},\\
(\ad_{h_3})^{[2]} = (a t +a^2) h_1 + t h_3\equiv a^2 h_1\pmod{\fc},\\
(\ad_d)^{[2]} = a t h_1 + t h_3 +d\equiv d\pmod{\fc},
\end{array}\end{equation}
(The expressions are somewhat dif\/ferent since we have chosen a
dif\/ferent basis but on this simple Lie algebra the 2-structure is
unique.)

2) The 2-structure on $\fwk(4;a)$ is given by the conditions
$(e_\alpha^{\pm})^{[2]}=0$ for all root vectors and
\begin{equation}
\label{2strwk4}
\begin{array}{l}
(\ad_{h_1})^{[2]} =a h_1 +(1+a) h_4,\\
(\ad_{h_2})^{[2]} = a h_2,\\
(\ad_{h_3})^{[2]} = h_3,\\
(\ad_{h_4})^{[2]} = h_4.\\
\end{array}
\end{equation}

\underline{$p=3$, Lie algebras}. {\it Brown\footnote{To interpret
the limit in \eqref{br2a}, set $\eps=1+\frac{1}{\alpha}$, and
$\fbr(2):=\fbr(2;\eps)$ for $\eps=1$.} algebras}:\index{$\fbr(2,
a)$, where $a\neq 0, -1$, Brown algebra}\index{$\fbr(2)$, Brown
algebra}\index{Brown algebra}\index{$\fbr(3)$, Brown algebra}
\begin{equation}\label{br2a}
\fbr(2, a)\text{~with CM~}
\begin{pmatrix}2&-1\\a&2\end{pmatrix}\text{~and $\fbr(2)=
``\mathop{\lim}\limits_{-\frac2a\tto 0}$''} \fbr(2,
a)\text{~with CM~}
\begin{pmatrix}2&-1\\-1&0\end{pmatrix}
\end{equation}
The ref\/lections change the value of the parameter, so
\begin{gather}\label{brso}
\fbr(2,a)\simeq \fbr(2,a')\Longleftrightarrow a'=-(1+a).
\\
\label{br3a}
1\fbr(3)\text{~with CM~}
\begin{pmatrix}2&-1&0\\-1&2&-1\\
0&-1&\ev\end{pmatrix}\text{~and ~} 2\fbr(3)\text{~with CM~}
\begin{pmatrix}2&-1&0\\-2&2&-1\\
0&-1&\ev\end{pmatrix}.
\end{gather}

\underline{$p=3$, Lie superalgebras}.

{\it Brown superalgebra} $\fbrj(2;3)$\index{$\fbrj(2;3)$, Brown
superalgebra}\index{Brown superalgebra} of $\sdim = 10|8$ (recently
discovered in \cite[Theorem 3.2(i)]{El1}; its Cartan matrices are
f\/irst listed in \cite{BGL3}) has the following Cartan matrices
\[
1)\ \begin{pmatrix}0&-1\\
-2&1\end{pmatrix}, \quad 2)\ \begin{pmatrix}0&-1\\
-1&\ev\end{pmatrix}, \quad 3)\ \begin{pmatrix}1&-1\\
-1&\ev\end{pmatrix}.
\]
The Lie superalgebra $\fbrj(2;3)$ is a super analog of the Brown
algebra $\fbr(2)=\fbrj(2;3)_\ev$, its even part;
$\fbrj(2;3)_\od=R(2\pi_1)$ is irreducible $\fbrj(2;3)_\ev$-module.

Elduque \cite{El1, El2, CE, CE2} considered a particular case of the
problem (\ref{kacg(a)}) and arranged the Lie (super)algebras he
discovered in a Supermagic Square all its entries being of the form
$\fg(A)$. These {\it Elduque and Cunha superalgebras}\index{Elduque
and Cunha superalgebra} are, indeed, exceptional ones. For the
complete list of their inequivalent Cartan matrices, reproduced
here, see \cite{BGL1}, where their presentation are also given; we
also reproduce the description of the even and odd parts of these
Lie superalgebras (all but one discovered by Elduque and whose
description in terms of symmetric composition algebras is due to
Elduque and Cunha), see Section~\ref{Secssr}.

\underline{$p=5$, Lie superalgebras}. {\it Brown
superalgebra}\index{$\fbrj(2;5)$, Brown superalgebra}\index{Brown
superalgebra} $\fbrj(2;5)$ of $\sdim = 10|12$, recently discovered
in~\cite{BGL3}, such that $\fbrj(2;5)_\ev=\fsp(4)$ and
$\fbrj(2;5)_\od=R(\pi_1+\pi_2)$ is an irreducible
$\fbrj(2;5)_\ev$-module\footnote{To the incredulous reader: The
Cartan subalgebra of $\fsp(4)$ is generated by $h_2$ and $2h_1+h_2$.
The highest weight vector is
$x_{10}=[[x_2,[x_2,[x_1,\;x_2]]],[[x_1,x_2],[x_1,x_2]]]$
and its weight is not a multiple of a fundamental weight, but
$(1,1)$. We encounter several more instances of non-fundamental
weights in descriptions of exceptions for $p=2$.}. The Lie
superalgebra $\fbrj(2;5)$ has the following Cartan matrices:
\[
\begin{lmatrix}
 2 & - \\
 1 & -
\end{lmatrix}
 \qquad 1)\ \begin{pmatrix}0&-1\\
-2&1\end{pmatrix}, \quad 2)\ \begin{pmatrix}0&-1\\
-3&2\end{pmatrix}.
\]

{\it Elduque superalgebra} $\fel(5; 5)$\index{$\fel(5; 5)$, Elduque
superalgebra}\index{Elduque superalgebra} of $\sdim=55|32$, where
$\fel(5; 5)_\ev=\fo(11)$ and $\fel(5; 5)_\od=\spin_{11}$. Its
inequivalent Cartan matrices, f\/irst described in \cite{BGL2}, are as
follows: \label{SSel5ssr}

Instead of joining nodes by four segments in the cases where
$A_{ij}=A_{ji}=1\equiv -4\mod 5$ we use one dotted segment.
\begin{figure}[ht]
\centerline{\parbox{.75\linewidth}{\includegraphics{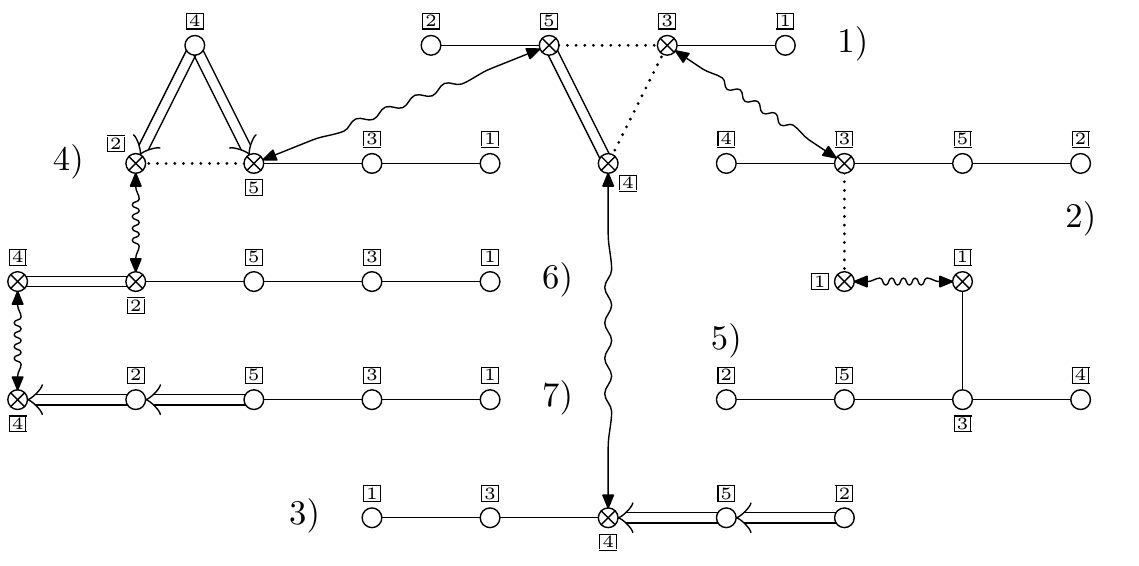}}}
\end{figure}
\begin{gather*}
1) \begin{pmatrix}
2&0&-1&0&0 \\
 0&2&0&0&-1 \\
 -1&0&0&-4&-4 \\
 0&0&-4&0&-2 \\
 0&-1&-4&-2&0
\end{pmatrix},\quad
2) \begin{pmatrix}
 0&0&-4&0&0 \\
 0&2&0&0&-1 \\
 -4&0&0&-1&-1 \\
 0&0&-1&2&0 \\
 0&-1&-1&0&2
 \end{pmatrix},
\\
 \boxed{3)} \begin{pmatrix}
 2&0&-1&0&0 \\
 0&2&0&0&-1 \\
 -1&0&2&-1&0 \\
 0&0&-1&0&2 \\
 0&-2&0&-1&2
 \end{pmatrix},
\quad
4) \begin{pmatrix}
 2&0&-1&0&0 \\
 0&0&0&2&-4 \\
 -1&0&2&0&-1 \\
 0&-1&0&2&-1 \\
 0&-4&-1&2&0
 \end{pmatrix},
\\
\boxed{5)} \begin{pmatrix}
 0&0&-1&0&0 \\
 0&2&0&0&-1 \\
 -1&0&2&-1&-1 \\
 0&0&-1&2&0 \\
 0&-1&-1&0&2
 \end{pmatrix},\quad
6) \begin{pmatrix}
 2&0&-1&0&0 \\
 0&0&0&-2&-1 \\
 -1&0&2&0&-1 \\
 0&-2&0&0&0 \\
 0&-1&-1&0&2
 \end{pmatrix}
\\
 \boxed{7)} \begin{pmatrix}
 2&0&-1&0&0 \\
 0&2&0&-1&-2 \\
 -1&0&2&0&-1 \\
 0&2&0&0&0 \\
 0&-1&-1&0&2
 \end{pmatrix},\quad
 8) \begin{pmatrix}
 -&-&2&3&4 \\
 5&-&1&-&- \\
 -&-&-&1&- \\
 -&6&-&-&1 \\
 2&-&-&-&- \\
 -&4&-&7&- \\
 -&-&-&6&-
\end{pmatrix}.
\end{gather*}

\subsection[Step 2: Studying $2\times 2$ and $3\times 3$ Cartan matrices]{Step 2: Studying $\boldsymbol{2\times 2}$ and $\boldsymbol{3\times 3}$ Cartan matrices}

1) We ask \texttt{Mathematica} to construct all possible matrices of
a specif\/ic size. The matrices are not normalized and they must not
be symmetrizable: we can not eliminate non-symmetrizable matrices at
this stage. Fortunately, all $2\times 2$ matrices are symmetrizable.

2) We ask \texttt{Mathematica} to eliminate the matrices with the
following properties:
\begin{equation}\label{a,b}
\begin{array}{ll}
\text{a)}&\text{Matrices $A$ for whose submatrix $B$ we know that
$\dim \fg(B)=\infty$;}\\
\text{b)}& \text{decomposable matrices.}
\end{array}
\end{equation}

3) Matrices with a row in each that dif\/fer from each other by a
nonzero factor are counted once, e.g.,
\[
\begin{pmatrix}
1&1\\
3&2
\end{pmatrix}  \cong   \begin{pmatrix}
2&2\\
3&2
\end{pmatrix}\cong   \begin{pmatrix}
6&6\\
6&4
\end{pmatrix}  .
\]

4) Equivalent matrices are counted once, where equivalence means
that one matrix can be obtained from the other one by simultaneous
transposition of rows and columns with the same numbers and {\it the
same parity}. For example,
\[
 \begin{pmatrix} 0 &\alpha &0&0\\
\alpha &\overline{0}&0&1\\
0&0&\overline{0}&1\\
0&1&1&\overline{0} \end{pmatrix}
\sim
 \begin{pmatrix} 0 &\alpha &1&0\\
\alpha &\overline{0}&0&0\\
1&0&\overline{0}&1\\
0&0&1&\overline{0} \end{pmatrix}
\sim
\begin{pmatrix} 0 &0 &0&1\\
0 &\overline{0}&\alpha&1\\
0&\alpha&\overline{0}&0\\
1&1&0&\overline{0} \end{pmatrix}.
\]
\textit{At the substeps $1.1)$--$1.4)$ we thus get a store of
Cartan matrices to be tested further.}

5) Now, we ask \texttt{SuperLie}\label{conj}, see \cite{Gr}, to construct the Lie
superalgebras $\mathfrak{g}(A)$ up to certain dimension (say, 256).
Having stored the Lie superalgebras $\mathfrak{g}(A)$ of dimension
$<256$ we increase the range again if there are any algebras left
(say, to 1024 or 2048). At this step, we {\bf conjecture} that the
dimension of any f\/inite dimensional simple Lie (super)algebra of the
form $\fg(A)$, where $A$ is of size $n\times n$, does not grow too
rapidly with $n$. Say, at least, not as fast as $n^{10}$.

If the dimension of $\mathfrak{g}(A)$ increases accordingly, then we
{\bf conjecture} that $\mathfrak{g}(A)$ is inf\/inite dimensional and
this Lie superalgebra is put away for a while (but not completely
eliminated as decomposable matrices that correspond to non-simple
algebras: The progress of science might require soon to investigate
{\bf how fast} the dimension grows with $n$: polynomially or
faster).

6) For the stored Cartan matrices $A$, we have $\dim \fg(A)<\infty$.
Once we get the full list all of such Cartan matrices of a given
size, we have to check if $\fg(A)$ is simple, one by one.

7) The vectors of parities of the generators $\Pty=(p_1, \dots ,
p_n)$ are only considered of the form
$(\od,\dots,\od,\ev,\dots,\ev)$.

\subsubsection[The case of $2\times 2$ Cartan matrices]{The case of $\boldsymbol{2\times 2}$ Cartan matrices}

On the diagonal we
may have 2, $\od$ or $\ev$, if the corresponding root is even; 0 or
1 if the root is odd. To be on the safe side, we redid the purely
even case. We have the following options to consider:
\begin{gather}\label{supercm21}
\begin{array}{ll}
\Pty=(\ev,
\ev)&\fa_1\quad\begin{pmatrix}2&2a\\2b&2\end{pmatrix}\simeq
\begin{pmatrix}2&2a\\b&\od\end{pmatrix}\simeq
\begin{pmatrix}\od&a\\2b&2\end{pmatrix}\simeq
\begin{pmatrix}\od&a\\b&\od\end{pmatrix}\simeq
\begin{pmatrix}b&ab\\ab&a\end{pmatrix},\\
&\fa_2\quad\begin{pmatrix}\od&2a\\-1&\ev\end{pmatrix},\qquad
\fa_3\quad\begin{pmatrix}\ev&-1\\-1&\ev\end{pmatrix};
\end{array}
\\
\label{supercm22}
\begin{array}{ll}
\Pty=(\od,
\ev)&\fa_4\quad\begin{pmatrix}0&-1\\2a&2\end{pmatrix}\simeq
\begin{pmatrix}0&-1\\ a&\overline{1}\end{pmatrix},\\
&\fa_5\quad\begin{pmatrix}\ev&-1\\-1&0\end{pmatrix},
\qquad\fa_6\quad\begin{pmatrix}1&a\\-1&\ev\end{pmatrix},\\
&\fa_7\quad\begin{pmatrix}1&a\\2b&2\end{pmatrix}\simeq
\begin{pmatrix}1&a\\b&\overline{1}\end{pmatrix};
\end{array}
\\
\label{supercm23}
\begin{array}{ll}
\Pty=(\od, \od)&\fa_8\quad\begin{pmatrix}1&a\\b&1\end{pmatrix}\simeq
\begin{pmatrix}b&ab\\ab&a\end{pmatrix},\\
&\fa_9\quad\begin{pmatrix}0&-1\\-1&0\end{pmatrix},
\qquad\fa_{10}\quad\begin{pmatrix}0&-1\\b&1\end{pmatrix}.
\end{array}
\end{gather}
Obviously, some of these CMs had appeared in the study of (twisted)
loops and the correspon\-ding Kac--Moody Lie (super)algebras. One could expect that
the reduction of the entries of $A$ modulo~$p$ might yield a f\/inite dimensional
algebra, but this does not happen.

\sssbegin{Lemma}\label{th_non-sy} If $A$ is non-symmetrizable, then
$\dim \fg(A)=\infty$.
\end{Lemma}

\begin{proof} We prove this by
inspection for $3\times 3$ matrices, but the general case does not
follow by reduction and induction: For example, for $p=2$ and the
normalized non-symmetrizable matrix (here the value of $\ast$ is
irrelevant) \[\begin{pmatrix} \ast& 1& 1& 0 \\ 1 &\ast& 0& 1 \\ 1
&0& \ast& 1 \\ 0 &1& a& \ast
\end{pmatrix},
\] where $a\neq 0, 1$, or analogous $n\times n$ matrix whose Dynkin
diagram is a loop, any $3\times 3$ submatrix is symmetrizable.

To eliminate non-symmetrizable Cartan matrices, and any loops of
length $>3$ in Dynkin diagrams, is, nevertheless, possible using
Lemmas~3.1, 3.3, and 3.10, 3.11 of \cite{WK}. (Van de Leur~\cite{vdL} used these Lemmas for $p=0$.)
\end{proof}

That was the idea of the proof. Now we pass to the case-by-case
study.

\subsection[Step 3: Studying $n\times n$ Cartan matrices for $n>3$]{Step 3: Studying $\boldsymbol{n\times n}$ Cartan matrices for $\boldsymbol{n>3}$}

By
Lemma \ref{th_non-sy} we will assume that $A$ is symmetrizable. The
idea is to use induction and the information found at each step.

\sssbegin{Hypothesis}\label{cj} {\em Each finite dimensional Lie
superalgebra of the form $\fg(A)$ possesses \linebreak a~``simplest''
Dynkin diagram~-- the one with only one odd node.}
\end{Hypothesis}

Therefore passing from $n\times n$ Cartan matrices to $(n+1)\times
(n+1)$ Cartan matrices it suf\/f\/ices to consider just two types of
$n\times n$ Cartan matrices: Purely even ones and the ``simplest'' ones~-- with only one odd node on their Dynkin diagrams. To the latter
ones only even node should be added.

\subsubsection{Further simplif\/ication of the algorithm}
Enlarging Cartan
matrices by adding new row and column, we let, for $n>4$, its only
non-zero elements occupy at most four slots (apart from the
diagonal). Justif\/ication: Lemmas from \S~3 in~\cite{WK} and Lemma~\ref{th_non-sy}.

Even this simplif\/ication still leaves lots of cases: To the 5 cases
to be enlarged for Cartan matrices of size $\leq 8$ that we
encounter for $p=0$, we have to add $16$ super cases, each producing
tens of possibilities in each of the major cases $p=2$, 3 and $5$.
To save several pages per each $n$ for each $p$, we have omitted the
results of enlargements of each Cartan matrix and give only the
f\/inal summary.

\subsection{On a quest for parametric families} Even for $2\times 2$
Cartan matrices we could have proceeded by ``enlarging'' but to
be on the safe side we performed the selection independently. We
considered only one or two parameters using the function called {\bf
ParamSolve} (of \texttt{SuperLie}, see~\cite{Gr}). It shows all cases where the
division by an expression possibly equal to zero occurred. Every
time \texttt{SuperLie} shows such a possibility we check it by hand;
these possibility are algebraic equations of the form
$\beta=f(\alpha)$, where~$\alpha$ and~$\beta$ are the parameters of
the CM. We saw that whenever $\alpha$ and $\beta$ are generic $\dim
\fg(A)$ grows too fast as compared with the {\it height} of the
element (i.e., the number of brackets in expressions like
$[a,[b,[c,d]]]$) that \texttt{SuperLie} should not exceed
constructing a Lie (super)algebra. We did not investigate if the
growth is polynomial or exponential, but def\/initely $\dim
\fg(A)=\infty$. For each pair of singular values of parameters
$\beta=f(\alpha)$, we repeat the computations again. In most cases,
the algebra is inf\/inite-dimensional, the exceptions being
$\beta=\alpha+1$ that nicely correspond to some of CMs we already
know, like $\fwk$ algebras.

For three parameters, we have equations of the form
$\gamma=f(\alpha,\beta)$. For generic $\alpha$ and $\beta$. the Lie
superalgebra $\fg(A)$ is inf\/inite-dimensional. For the singular
cases given by \texttt{SuperLie}, the constraints are of the form
$\beta=g(\alpha)$. Now we face two possibilities: If $\gamma$ is a
constant, then we just use the result of the previous step, when we
dealt with two parameters. In the rare cases where $\gamma$ is not a
constant and depends on the parameter $\alpha$, we have to recompute
again and again the $\dim\fg(A)$ is inf\/inite in these cases.

We f\/ind Cartan matrices of size $4\times 4$ and larger by ``enlarging''. For $p=2$, we see that $3\times 3$ CMs with parameters
can be extended to $4\times 4$ CMs. However, $4\times 4$ CMs cannot
be extended to $5\times 5$ CMs whose Lie (super)algebras are of
f\/inite dimension. For $p>2$, even $3\times 3$ CMs cannot be
extended.

\subsection[Super and modular cases: Summary of new features (as compared
with simple Lie algebras over $\Cee$)]{Super and modular cases: Summary of new features\\ (as compared
with simple Lie algebras over $\boldsymbol{\Cee}$)}

\underline{The super case, $p=0$}.
\begin{enumerate}\itemsep=0pt
\item[1)] There are three types of nodes ($\mbullet$, $\motimes$ and $\mcirc$),

\item[2)] there may occur a loop but only of length 3;

\item[3)] there is at most 1 parameter, but 1 parameter may occur;

\item[4)] to one algebra several inequivalent Cartan matrices can
correspond.
\end{enumerate}

\underline{The modular case}. \underline{For Lie algebras}, new
features are the same as in the $p=0$ super case; additionally there
appear new types of nodes ($\odot$ and $\ast$).

\section[The answer: The case where $p>5$]{The answer: The case where $\boldsymbol{p>5}$} \label{Sans>5}

This case is the simplest one since it does not dif\/fer much from the
$p=0$ case, where the answer is known.

\underline{Simple Lie algebras}:
\begin{enumerate}\itemsep=0pt
\item[1)]Lie algebras obtained from their $p=0$ analogs by reducing modulo
$p$. We thus get
\begin{enumerate}\itemsep=0pt
\item[]
the CM versions of $\fsl$, namely: either simple $\fsl(n)$ or
$\fgl(pn)$ whose ``simple core'' is $\fp\fsl(pn)$;
\item[]
the orthogonal algebras $\fo(2n+1)$ and $\fo(2n)$;
\item[]
the symplectic algebras $\fsp(2n)$;
\item[]
the exceptional algebras are $\fg(2)$, $\ff(4)$, $\fe(6)$, $\fe(7)$,
$\fe(8)$.
\end{enumerate}
\end{enumerate}

\underline{Simple Lie superalgebras}

Lie superalgebras obtained from their $p=0$ analogs by reducing
modulo $p$. We thus get
\begin{enumerate}\itemsep=0pt
\item[1)] the CM versions of $\fsl$, namely: either simple $\fsl(m|n)$ or
$\fgl(a|pk+a)$ whose ``simple core'' is $\fp\fsl(a|pk+a)$ and
$\fp\fsl^{(1)}(a|pk+a)$ if $a=kn$;
\item[2)] the ortho-symplectic algebras $\fosp(m|2n)$;
\item[3)] a parametric family $\fosp(4|2; a)$;
\item[4)] the exceptional algebras are $\fag(2)$ and $\fab(3)$.
\end{enumerate}

\section[The answer: The case where $p=5$]{The answer: The case where $\boldsymbol{p=5}$}\label{Sans5}

\underline{Simple Lie algebras}:
\begin{enumerate}\itemsep=0pt
\item[1)] same as in Section~\ref{Sans>5} for $p=5$.
\end{enumerate}

\underline{Simple Lie superalgebras}
\begin{enumerate}\itemsep=0pt
\item[1)] same as in Section~\ref{Sans>5} for $p=5$ and several new exceptions:
\item[2)] The Brown superalgebras \cite{BGL3}: $\fbrj(2;5)$ such that
$\fbrj(2;5)_\ev=\fsp(4)$ and the $\fbrj(2;5)_\ev$-module
$\fbrj(2;5)_\od=R(\pi_1+\pi_2)$ is irreducible with the highest
weight vector
\[
x_{10}=[[x_2, [x_2, [x_1, x_2]]],  [[x_1, x_2], [x_1, x_2]]]
\] (for the CM 2): with the two Cartan matrices
\[
\left(
\begin{array}{ll}
 2 & - \\
 1 & -
\end{array}
\right)
\qquad 1)\; \begin{pmatrix}0&-1\\
-2&1\end{pmatrix}, \quad 2)\; \begin{pmatrix}0&-1\\
-3&2\end{pmatrix}.
\]
\item[3)] The Elduque superalgebra $\fel(5; 5)$. Having found out one
Cartan matrix of $\fel(5; 5)$, we have listed them all, see~\ref{SSel5ssr}.
\end{enumerate}

\section[The answer: The case where $p=3$]{The answer: The case where $\boldsymbol{p=3}$}\label{Sans3}

\underline{Simple Lie algebras}:
\begin{enumerate}\vspace{-1mm}\itemsep=-0.5pt
\item[1)] same as in Section~\ref{Sans>5} for $p=3$, except $\fg(2)$ which is
not simple but contains a unique minimal ideal isomorphic to
$\fpsl(3)$, and the following additional exceptions:
\item[2)] the Brown algebras $\fbr(2; a)$~and $\fbr(2)$ as well as
$\fbr(3)$, see Section~\ref{SSsteps}.\vspace{-1mm}
\end{enumerate}

\underline{Simple Lie superalgebras}
\begin{enumerate}\vspace{-1mm}\itemsep=-0.5pt
\item[1)] same as in Section~\ref{Sans>5} for $p=3$ and $\fe(6)$ (with CM) which
is not simple but has a ``simple core'' $\fe(6)/\fc$;

\item[2)] the Brown superalgebras, see Section~\ref{SSsteps};

\item[3)] the Elduque and Cunha superalgebras, see \cite{CE2, BGL1}. They
are respective ``enlargements'' of the following Lie algebras
(but can be also obtained by enlarging certain Lie superalgebras):
\begin{enumerate}\vspace{-1mm}\itemsep=-0.5pt
\item[]$\fg(2,3)$ ($\fgl(3)$ yields $2\fg(1,6)$ and $1\fg(2,3)$) (with CM)
has a simple core $\fbj:=\fg(2,3)^{(1)}\!/\fc$;

\item[]$\fg(3,6)$ ($\fsl(4)$ yields $7\fg(3,6)$);

\item[]$\fg(3,3)$ ($\fsp(6)$ yields $1\fg(3,3)$ and $10\fg(3,3)$);

\item[]$\fg(4,3)$ ($\fo(7)$ yields $1\fg(4,3)$);

\item[]$\fg(8,3)$ ($\ff(4)$ yields $1\fg(8,3)$);

\item[]$\fg(2,6)$ ($\fsl(5)$ yields $3\fg(2,6)$) (with CM) has a simple
core $\fg(2,6)^{(1)}/\fc$;

\item[]$\fg(4,6)$ ($\fgl(6)$ yields $3\fg(4,6)$ and $\fo(10)$ yields
$7\fg(4,6)$);

\item[]$\fg(6,6)$ ($\fo(11)$ yields $21\fg(6,6)$);

\item[]$\fg(8,6)$ ($\fsl(7)$ yields $8\fg(8,6)$ and $\fe(6)$ yields
$3\fg(8,6)$);\vspace{-1mm}
\end{enumerate}

\item[4)] the Lie superalgebra $\fel(5; 3)$ we have discovered is a $p=3$
version of the Elduque superalgebra $\fel(5; 5)$: Their Cartan
matrices (whose elements are represented by non-positive integers)
7) for $\fel(5; 5)$ and 1) for $\fel(5; 3)$ are identical after a
permutation of indices (that is why we baptized $\fel(5; 3)$ so). It
can be obtained as an ``enlargement'' of any of the following
Lie (super)algebras: $\fsp(8)$, $\fsl(1|4)$, $\fsl(2|3)$,
$\fosp(4|4)$, $\fosp(6|2)$, $\fg(3, 3)$.\vspace{-1mm}
\end{enumerate}

\vspace{-1mm}

\subsection{Elduque and Cunha superalgebras: Systems of simple
roots}\label{Secssr}

For details of description of Elduque and Cunha superalgebras in
terms of symmetric composition algebras, see \cite{El1,CE,CE2}. Here
we consider the simple Elduque and Cunha superalgebras with Cartan
matrix for $p=3$. In what follows, we list them using somewhat
shorter notations as compared with the original ones: Hereafter
$\fg(A,B)$ denotes the superalgebra occupying $(A, B)$th slot in the
Elduque Supermagic Square; the f\/irst Cartan matrix is usually the
one given in~\cite{CE}, where only one Cartan matrix is given; the
other matrices are obtained from the f\/irst one by means of
ref\/lections. Accordingly, $i\fg(A, B)$ is the shorthand for the
realization of $\fg(A, B)$ by means of the $i$th Cartan matrix.
There are no instances of isotropic even ref\/lections. On notation in
the following tables, see Section~\ref{recmat}.

\vspace{-1mm}

\subsubsection[$\fg{}(1,6)$ of $\sdim =21|14$]{$\boldsymbol{\fg{}(1,6)}$ of $\boldsymbol{\sdim =21|14}$}

We have $\fg(1,
6)_\ev=\fsp(6)$ and $\fg(1,6)_\od=R(\pi_3)$.
\begin{figure}[h!]\centering
\parbox{.5\linewidth}{\mbox{}\hfill$
\begin{lmatrix}
 -&-&2 \\
 -&-&1
 \end{lmatrix}
$\quad\mbox{}}\hfill
\parbox{.5\linewidth}{\includegraphics[scale=0.9]{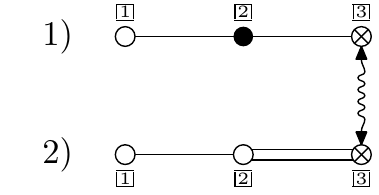}}
\end{figure}
 \[ 1) \begin{pmatrix}
 2&-1&0 \\
 -1&1&-1 \\
 0&-1&0
 \end{pmatrix}, \quad
\boxed{2)} \begin{pmatrix}
 2&-1&0 \\
 -1&2&-2 \\
 0&-2&0
 \end{pmatrix}.
\]

\subsubsection[$\fg(2,3)$ of $\sdim =12/10|14$]{$\boldsymbol{\fg(2,3)}$ of $\boldsymbol{\sdim =12/10|14}$}

We have
$\fg(2,3)_\ev=\fgl(3)\oplus \fsl(2)$ and
$\fg(2,3)_\od=\fpsl(3)\otimes \id$.
\begin{figure}[h!]\centering
\parbox{.35\linewidth}{$
\begin{lmatrix}
 -&-&2 \\
 3&4&1 \\
 2&5&- \\
 5&2&- \\
 4&3&-
 \end{lmatrix}
$} \quad
\parbox{.6\linewidth}{\includegraphics{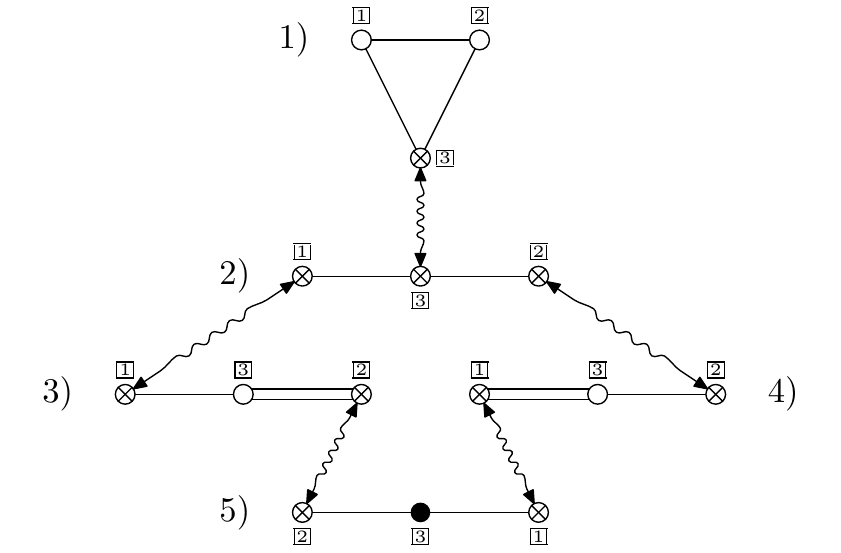}}
\end{figure}
\begin{gather*}
\boxed{1)} \begin{pmatrix}
 2&-1&-1 \\
 -1&2&-1 \\
 -1&-1&0
 \end{pmatrix},\quad
\underline{2)} \begin{pmatrix}
 0&0&-1 \\
 0&0&-1 \\
 -1&-1&0
 \end{pmatrix},\quad
3) \begin{pmatrix}
 0&0&-1 \\
 0&0&-2 \\
 -1&-2&2
 \end{pmatrix},\\
4) \begin{pmatrix}
 0&0&-2 \\
 0&0&-1 \\
 -2&-1&2
 \end{pmatrix},\quad
\underline{5)} \begin{pmatrix}
 0&0&-1 \\
 0&0&-1 \\
 -1&-1&1
 \end{pmatrix}.
\end{gather*}

\subsubsection[$\fg{}(3,6)$ of $\sdim =36|40$]{$\boldsymbol{\fg{}(3,6)}$ of $\boldsymbol{\sdim =36|40}$}

We have $\fg(3,
6)_\ev=\fsp(8)$ and $\fg(3,6)_\od=R(\pi_3)$.
\begin{gather*}
\begin{lmatrix}
 2&-&-&3 \\
 1&4&-&5 \\
 5&-&-&1 \\
 -&2&-&6 \\
 3&6&-&2 \\
 -&5&7&4 \\
 -&-&6&-
 \end{lmatrix}, \qquad
1) \begin{pmatrix}
 0&-1&0&0 \\
 -1&2&-1&0 \\
 0&-1&1&-1 \\
 0&0&-1&0
 \end{pmatrix},\quad
\underline{2)} \begin{pmatrix}
 0&-1&0&0 \\
 -1&0&-1&0 \\
 0&-1&1&-1 \\
 0&0&-1&0
 \end{pmatrix},
\\
3) \begin{pmatrix}
 0&-1&0&0 \\
 -1&2&-1&0 \\
 0&-1&2&-2 \\
 0&0&-1&0
 \end{pmatrix}, \quad
4)
\begin{pmatrix}
 2&-1&0&0 \\
 -1&0&-2&0 \\
 0&-2&2&-1 \\
 0&0&-1&0
 \end{pmatrix},\quad
5) \begin{pmatrix}
 0&-1&0&0 \\
 -2&0&-1&0 \\
 0&-1&2&-2 \\
 0&0&-1&0
 \end{pmatrix}, \\
6) \begin{pmatrix}
 2&-1&0&0 \\
 -1&0&-2&0 \\
 0&-2&0&-2 \\
 0&0&-1&0
 \end{pmatrix},\quad
\boxed{7)} \begin{pmatrix}
 2&-1&0&0 \\
 -1&2&-1&-1 \\
 0&-1&0&-1 \\
 0&-1&-1&2
 \end{pmatrix}.
\end{gather*}
\begin{figure}[h!]\centering
\includegraphics{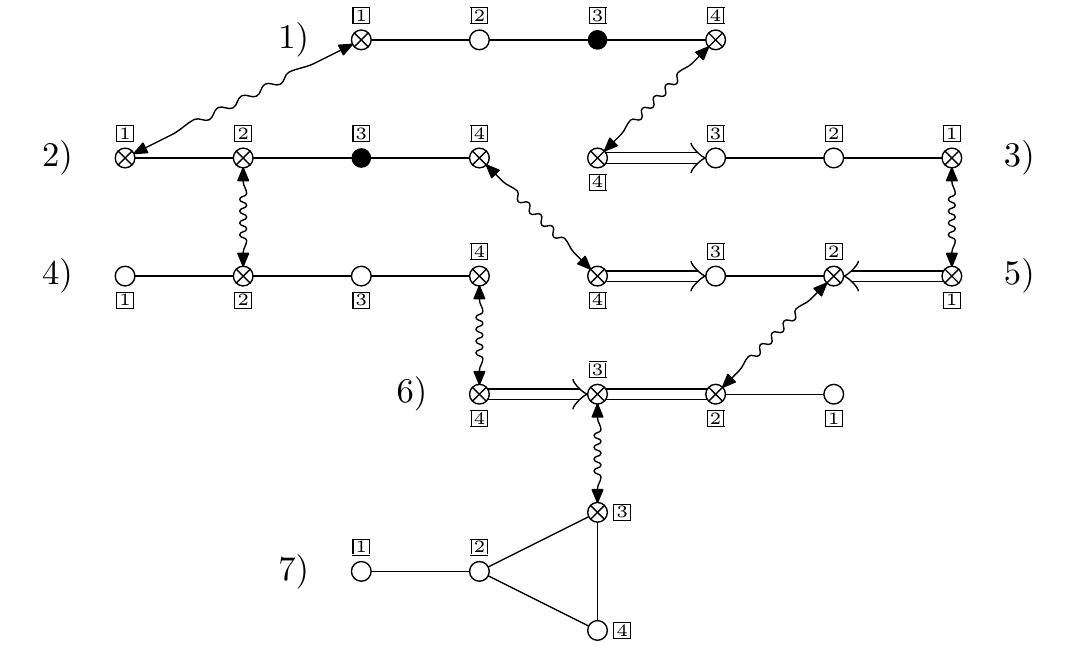}
\end{figure}

\subsubsection[$\fg{}(3,3)$ of $\sdim =23/21|16$]{$\boldsymbol{\fg{}(3,3)}$ of $\boldsymbol{\sdim =23/21|16}$}

 We have
$\fg(3,3)_\ev=(\fo(7)\oplus \Kee z)\oplus \Kee d$ and
$\fg(3,3)_\od=(\spin_7)_+\oplus (\spin_7)_-$; the action of $d$
separates the summands~-- identical $\fo(7)$-modules $\spin_7$~--
acting on one as the scalar multiplication by~1, on the other one by~$-1$.
\begin{figure}[h!]\centering
\includegraphics{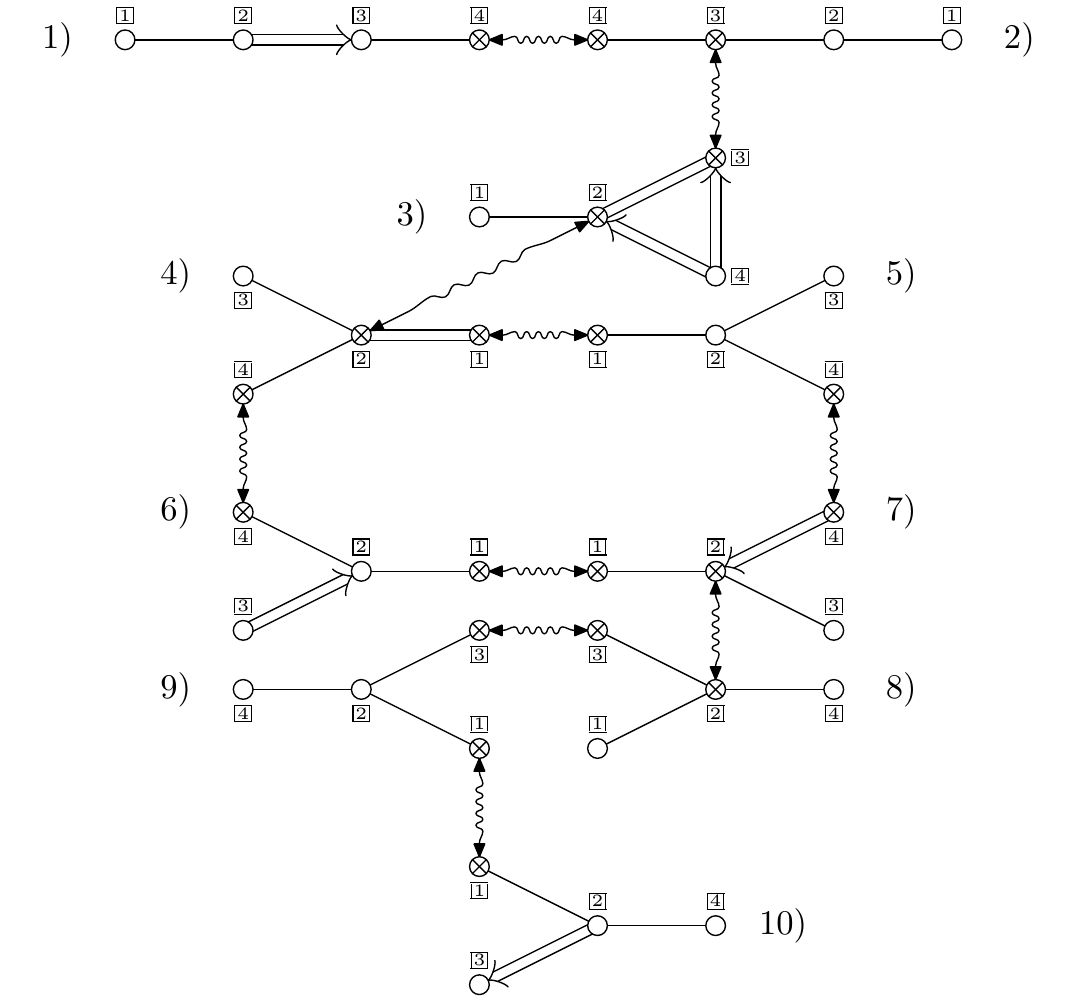}
\end{figure}\newpage
\vspace*{-12mm}

\begin{gather*}
\begin{lmatrix}
 -&-&-&2 \\
 -&-&3&1 \\
 -&4&2&- \\
 5&3&-&6 \\
 4&-&-&7 \\
 7&-&-&4 \\
 6&8&-&5 \\
 -&7&9&- \\
 10&-&8&- \\
 9&-&-&-
 \end{lmatrix}, \quad \begin{matrix}\boxed{1)} \begin{pmatrix}
 2&-1&0&0 \\
 -1&2&-1&0 \\
 0&-2&2&-1 \\
 0&0&-1&0
 \end{pmatrix},\quad
2) \begin{pmatrix}
 2&-1&0&0 \\
 -1&2&-1&0 \\
 0&-1&0&-1 \\
 0&0&-1&0
 \end{pmatrix},\\
 3) \begin{pmatrix}
 2&-1&0&0 \\
 -1&0&-2&-2 \\
 0&-2&0&-2 \\
 0&-1&-1&2
 \end{pmatrix}, \quad
4) \begin{pmatrix}
 0&-1&0&0 \\
 -2&0&-1&-1 \\
 0&-1&2&0 \\
 0&-1&0&0
 \end{pmatrix},
 \end{matrix}
\\
\quad 5) \begin{pmatrix}
 0&-1&0&0 \\
 -1&2&-1&-1 \\
 0&-1&2&0 \\
 0&-1&0&0
 \end{pmatrix},\quad
6) \begin{pmatrix}
 0&-1&0&0 \\
 -1&2&-2&-1 \\
 0&-1&2&0 \\
 0&-1&0&0
 \end{pmatrix}, \quad
7)\begin{pmatrix}
 0&-1&0&0 \\
 -1&0&-1&-2 \\
 0&-1&2&0 \\
 0&-1&0&0
 \end{pmatrix},\\
8) \begin{pmatrix}
 2&-1&-1&0 \\
 -2&0&-2&-1 \\
 -1&-1&0&0 \\
 0&-1&0&2
 \end{pmatrix},\quad
9)
\begin{pmatrix}
 0&0&-1&0 \\
 0&2&-1&-1 \\
 -1&-1&0&0 \\
 0&-1&0&2
 \end{pmatrix},\quad
\boxed{10)} \begin{pmatrix}
 0&0&-1&0 \\
 0&2&-1&-1 \\
 -1&-2&2&0 \\
 0&-1&0&2
 \end{pmatrix}.
\end{gather*}

\subsubsection[$\fg{}(4,3)$ of $\sdim =24|26$]{$\boldsymbol{\fg{}(4,3)}$ of $\boldsymbol{\sdim =24|26}$}

We have
$\fg(4,3)_\ev=\fsp(6)\oplus \fsl(2)$ and
$\fg(4,3)_\od=R(\pi_2)\otimes \id$.
\begin{figure}[h!]\centering
\includegraphics{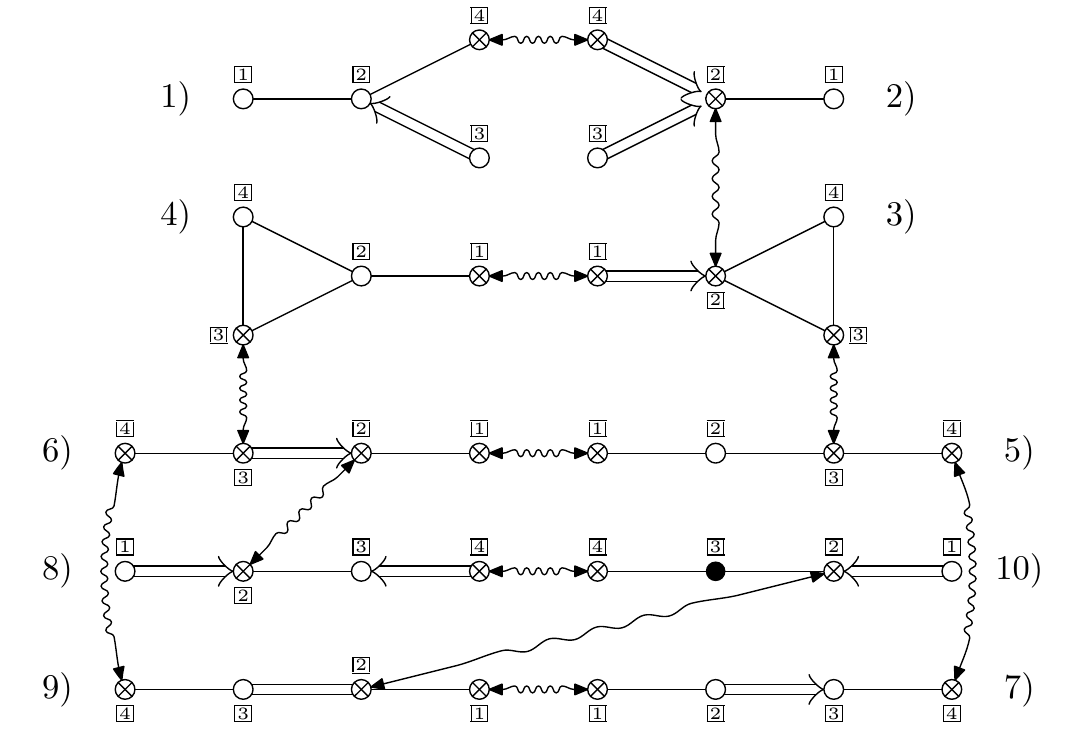}
\end{figure}
\begin{gather*}
\begin{lmatrix}
 -&-&-&2 \\
 -&3&-&1 \\
 4&2&5&- \\
 3&-&6&- \\
 6&-&3&7 \\
 5&8&4&9 \\
 9&-&-&5 \\
 -&6&-&10 \\
 7&10&-&6 \\
 -&9&-&8
 \end{lmatrix},\quad \begin{matrix}\boxed{1)} \begin{pmatrix}
 2&-1&0&0 \\
 -1&2&-2&-1 \\
 0&-1&2&0 \\
 0&-1&0&0
 \end{pmatrix},\quad
2) \begin{pmatrix}
 2&-1&0&0 \\
 -1&0&-2&-2 \\
 0&-1&2&0 \\
 0&-1&0&0
 \end{pmatrix},\\
3) \begin{pmatrix}
 0&-1&0&0 \\
 -2&0&-1&-1 \\
 0&-1&0&-1 \\
 0&-1&-1&2
 \end{pmatrix}, \quad
4) \begin{pmatrix}
 0&-1&0&0 \\
 -1&2&-1&-1 \\
 0&-1&0&-1 \\
 0&-1&-1&2
 \end{pmatrix},\end{matrix}
\\
5) \begin{pmatrix}
 0&-1&0&0 \\
 -1&2&-1&0 \\
 0&-1&0&-1 \\
 0&0&-1&0
 \end{pmatrix},\quad
\underline{6)} \begin{pmatrix}
 0&-1&0&0 \\
 -1&0&-2&0 \\
 0&-1&0&-1 \\
 0&0&-1&0
 \end{pmatrix}, \quad
7) \begin{pmatrix}
 0&-1&0&0 \\
 -1&2&-1&0 \\
 0&-2&2&-1 \\
 0&0&-1&0
 \end{pmatrix},\\
8) \begin{pmatrix}
 2&-1&0&0 \\
 -2&0&-1&0 \\
 0&-1&2&-2 \\
 0&0&-1&0
 \end{pmatrix},\quad
9) \begin{pmatrix}
 0&-1&0&0 \\
 -1&0&-2&0 \\
 0&-2&2&-1 \\
 0&0&-1&0
 \end{pmatrix},\quad
10) \begin{pmatrix}
 2&-1&0&0 \\
 -2&0&-1&0 \\
 0&-1&1&-1 \\
 0&0&-1&0
 \end{pmatrix}.
\end{gather*}

\subsubsection[$\fg{}(2,6)$ of $\sdim =36/34|20$]{$\boldsymbol{\fg{}(2,6)}$ of $\boldsymbol{\sdim =36/34|20}$}

We have
$\fg(2,6)_\ev=\fgl(6)$ and $\fg(2,6)_\od=R(\pi_3)$.
\begin{figure}[h!]\centering
\includegraphics{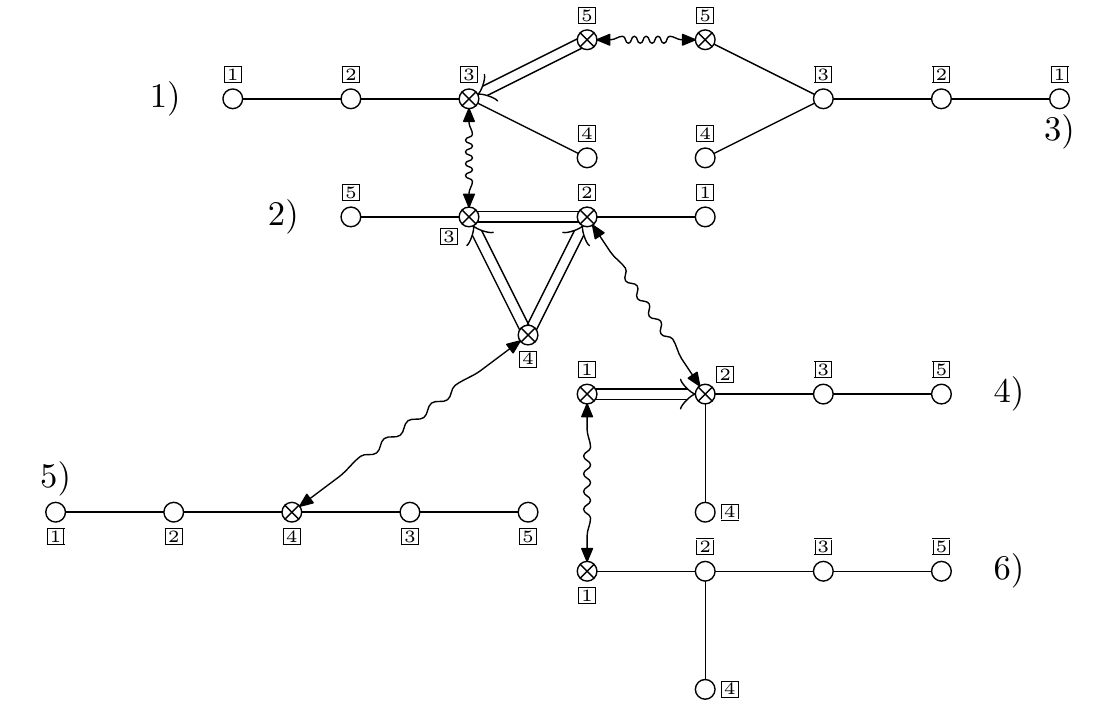}
\end{figure}
\begin{gather*}
\begin{pmatrix}
 -&-&2&-&3 \\
 -&4&1&5&- \\
 -&-&-&-&1 \\
 6&2&-&-&- \\
 -&-&-&2&- \\
 4&-&-&-&-
 \end{pmatrix}, \quad
 1) \begin{pmatrix}
 2&-1&0&0&0 \\
 -1&2&-1&0&0 \\
 0&-1&0&-1&-2 \\
 0&0&-1&2&0 \\
 0&0&-1&0&0
 \end{pmatrix}, \quad
2) \begin{pmatrix}
 2&-1&0&0&0 \\
 -1&0&-2&-2&0 \\
 0&-2&0&-2&-1 \\
 0&-1&-1&0&0 \\
 0&0&-1&0&2
 \end{pmatrix},
\\
 \boxed{3)} \begin{pmatrix}
 2&-1&0&0&0 \\
 -1&2&-1&0&0 \\
 0&-1&2&-1&-1 \\
 0&0&-1&2&0 \\
 0&0&-1&0&0
 \end{pmatrix} \sim
\boxed{6)} \begin{pmatrix}
 0&-1&0&0&0 \\
 -1&2&-1&-1&0 \\
 0&-1&2&0&-1 \\
 0&-1&0&2&0 \\
 0&0&-1&0&2
 \end{pmatrix}, \\
4)
 \begin{pmatrix}
 0&-1&0&0&0 \\
 -2&0&-1&-1&0 \\
 0&-1&2&0&-1 \\
 0&-1&0&2&0 \\
 0&0&-1&0&2
 \end{pmatrix}, \quad\boxed{5)}\begin{pmatrix}
 2&-1&0&0&0 \\
 -1&2&0&-1&0 \\
 0&0&2&-1&-1 \\
 0&-1&-1&0&0 \\
 0&0&-1&0&2
 \end{pmatrix}.
\end{gather*}

\newpage

\subsubsection[$\fg{}(8,3)$ of $\sdim =55|50$]{$\boldsymbol{\fg{}(8,3)}$ of $\boldsymbol{\sdim =55|50}$}

We have
$\fg(8,3)_\ev=\ff(4)\oplus \fsl(2)$ and
$\fg(8,3)_\od=R(\pi_4)\otimes \id$.
\begin{center}
\includegraphics{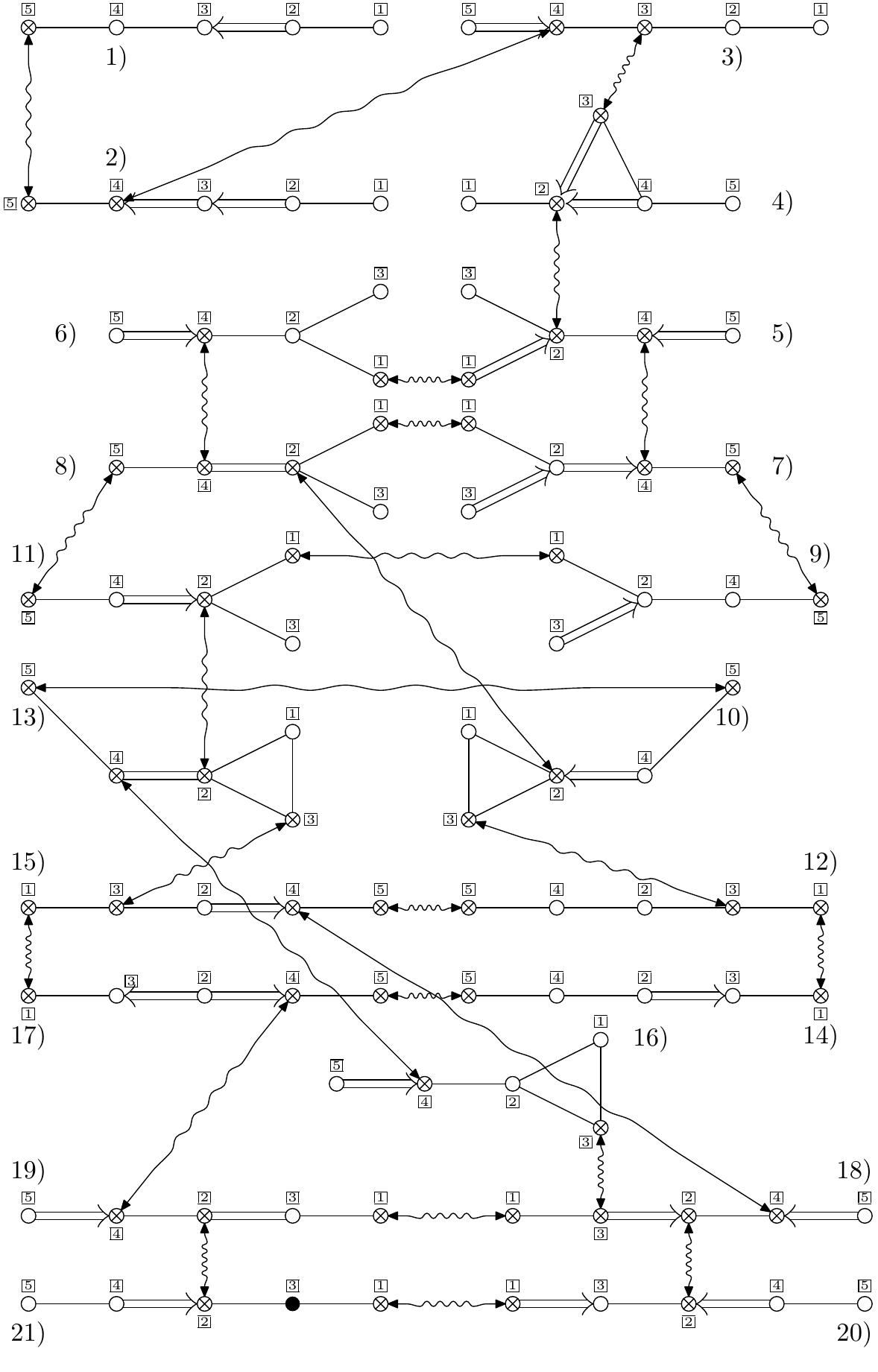}
\end{center}
\newpage
\vspace*{-12mm}

\begin{gather*}
{\arraycolsep=2.0pt \begin{pmatrix}
 -&-&-&-&2 \\
 -&-&-&3&1 \\
 -&-&4&2&- \\
 -&5&3&-&- \\
 6&4&-&7&- \\
 5&-&-&8&- \\
 8&-&-&5&9 \\
 7&10&-&6&11 \\
 11&-&-&-&7 \\
 -&8&12&-&13 \\
 9&13&-&-&8 \\
 14&-&10&-&15 \\
 -&11&15&16&10 \\
 12&-&-&-&17 \\
 17&-&13&18&12 \\
 -&-&18&13&- \\
 15&-&-&19&14 \\
 19&20&16&15&- \\
 18&21&-&17&- \\
 21&18&-&-&- \\
 20&19&-&-&-
 \end{pmatrix}},\quad
{\arraycolsep=2pt \begin{matrix}
\boxed{1)}
\begin{pmatrix}
 2&-1&0&0&0 \\
 -1&2&-1&0&0 \\
 0&-2&2&-1&0 \\
 0&0&-1&2&-1 \\
 0&0&0&1&0
 \end{pmatrix}, \\
2) \begin{pmatrix}
 2&-1&0&0&0 \\
 -1&2&-1&0&0 \\
 0&-2&2&-1&0 \\
 0&0&-2&0&-1 \\
 0&0&0&-1&0
 \end{pmatrix},\\
3) \begin{pmatrix}
 2&-1&0&0&0 \\
 -1&2&-1&0&0 \\
 0&-1&0&-1&0 \\
 0&0&-1&0&-2 \\
 0&0&0&-1&2
 \end{pmatrix}, \\
4) \begin{pmatrix}
 2&-1&0&0&0 \\
 -1&0&-2&-2&0 \\
 0&-1&0&-1&0 \\
 0&-1&-1&2&-1 \\
 0&0&0&-1&2
 \end{pmatrix},
 \end{matrix}} \quad
{\arraycolsep=2pt\begin{matrix}
5)  \begin{pmatrix}
 0&-1&0&0&0 \\
 -2&0&-1&-1&0 \\
 0&-1&2&0&0 \\
 0&-1&0&0&-2 \\
 0&0&0&-1&2
 \end{pmatrix},\\
6) \begin{pmatrix}
 0&-1&0&0&0 \\
 -1&2&-1&-1&0 \\
 0&-1&2&0&0 \\
 0&-1&0&0&-2 \\
 0&0&0&-1&2
 \end{pmatrix}, \\
7) \begin{pmatrix}
 0&-1&0&0&0 \\
 -1&2&-2&-1&0 \\
 0&-1&2&0&0 \\
 0&-2&0&0&-1 \\
 0&0&0&-1&0
 \end{pmatrix},\\
8) \begin{pmatrix}
 0&-1&0&0&0 \\
 -1&0&-1&-2&0 \\
 0&-1&2&0&0 \\
 0&-2&0&0&-1 \\
 0&0&0&-1&0
 \end{pmatrix},
 \end{matrix}}
\\
{\arraycolsep=2pt 9) \begin{pmatrix}
 0&-1&0&0&0 \\
 -1&2&-2&-1&0 \\
 0&-1&2&0&0 \\
 0&-1&0&2&-1 \\
 0&0&0&-1&0
 \end{pmatrix},  \quad
10)
  \begin{pmatrix}
 2&-1&-1&0&0 \\
 1&0&1&2&0 \\
 1&1&0&0&0 \\
 0&-1&0&2&-1 \\
 0&0&0&-1&0
 \end{pmatrix} , \quad
11) \begin{pmatrix}
 0&-1&0&0&0 \\
 -1&0&-1&-2&0 \\
 0&-1&2&0&0 \\
 0&-1&0&2&-1 \\
 0&0&0&-1&0
 \end{pmatrix}},\\
{\arraycolsep=2pt 12) \begin{pmatrix}
 0&0&-1&0&0 \\
 0&2&-1&-1&0 \\
 -1&-1&0&0&0 \\
 0&-1&0&2&-1 \\
 0&0&0&-1&0
 \end{pmatrix},
\quad
13)
 \begin{pmatrix}
 2&-1&-1&0&0 \\
 -1&0&-1&-2&0 \\
 -1&-1&0&0&0 \\
 0&-2&0&0&-1 \\
 0&0&0&-1&0
 \end{pmatrix},\quad
14) \begin{pmatrix}
 0&0&-1&0&0 \\
 0&2&-1&-1&0 \\
 -1&-2&2&0&0 \\
 0&-1&0&2&-1 \\
 0&0&0&-1&0
 \end{pmatrix}} ,\\
{\arraycolsep=2pt 15) \begin{pmatrix}
 0&0&-1&0&0 \\
 0&2&-1&-1&0 \\
 -1&-1&0&0&0 \\
 0&-2&0&0&-1 \\
 0&0&0&-1&0
 \end{pmatrix}, \quad
16)
 \begin{pmatrix}
 2&-1&-1&0&0 \\
 -1&2&-1&-1&0 \\
 -1&-1&0&0&0 \\
 0&-1&0&0&-2 \\
 0&0&0&-1&2
 \end{pmatrix},\quad
17) \begin{pmatrix}
 0&0&-1&0&0 \\
 0&2&-1&-1&0 \\
 -1&-2&2&0&0 \\
 0&-2&0&0&-1 \\
 0&0&0&-1&0
 \end{pmatrix}},\\
{\arraycolsep=2pt 18) \begin{pmatrix}
 0&0&-1&0&0 \\
 0&0&-2&-1&0 \\
 -1&-1&0&0&0 \\
 0&-1&0&0&-2 \\
 0&0&0&-1&2
 \end{pmatrix}, \quad
19)
 \begin{pmatrix}
 0&0&-1&0&0 \\
 0&0&-2&-1&0 \\
 -1&-2&2&0&0 \\
 0&-1&0&0&-2 \\
 0&0&0&-1&2
 \end{pmatrix},\quad
20) \begin{pmatrix}
 0&0&-1&0&0 \\
 0&0&-1&-2&0 \\
 -2&-1&2&0&0 \\
 0&-1&0&2&-1 \\
 0&0&0&-1&2
 \end{pmatrix}} ,\\
{\arraycolsep=2pt 21) \begin{pmatrix}
 0&0&-1&0&0 \\
 0&0&-1&-2&0 \\
 -1&-1&1&0&0 \\
 0&-1&0&2&-1 \\
 0&0&0&-1&2
 \end{pmatrix}}.
\end{gather*}
\newpage

\subsubsection[$\fg{}(4,6)$ of $\sdim =66|32$]{$\boldsymbol{\fg{}(4,6)}$ of $\boldsymbol{\sdim =66|32}$}

We have
$\fg(4,6)_\ev=\fo(12)$ and $\fg(4,6)_\od=R(\pi_{5})$.
\begin{figure}[h!]\centering
\includegraphics{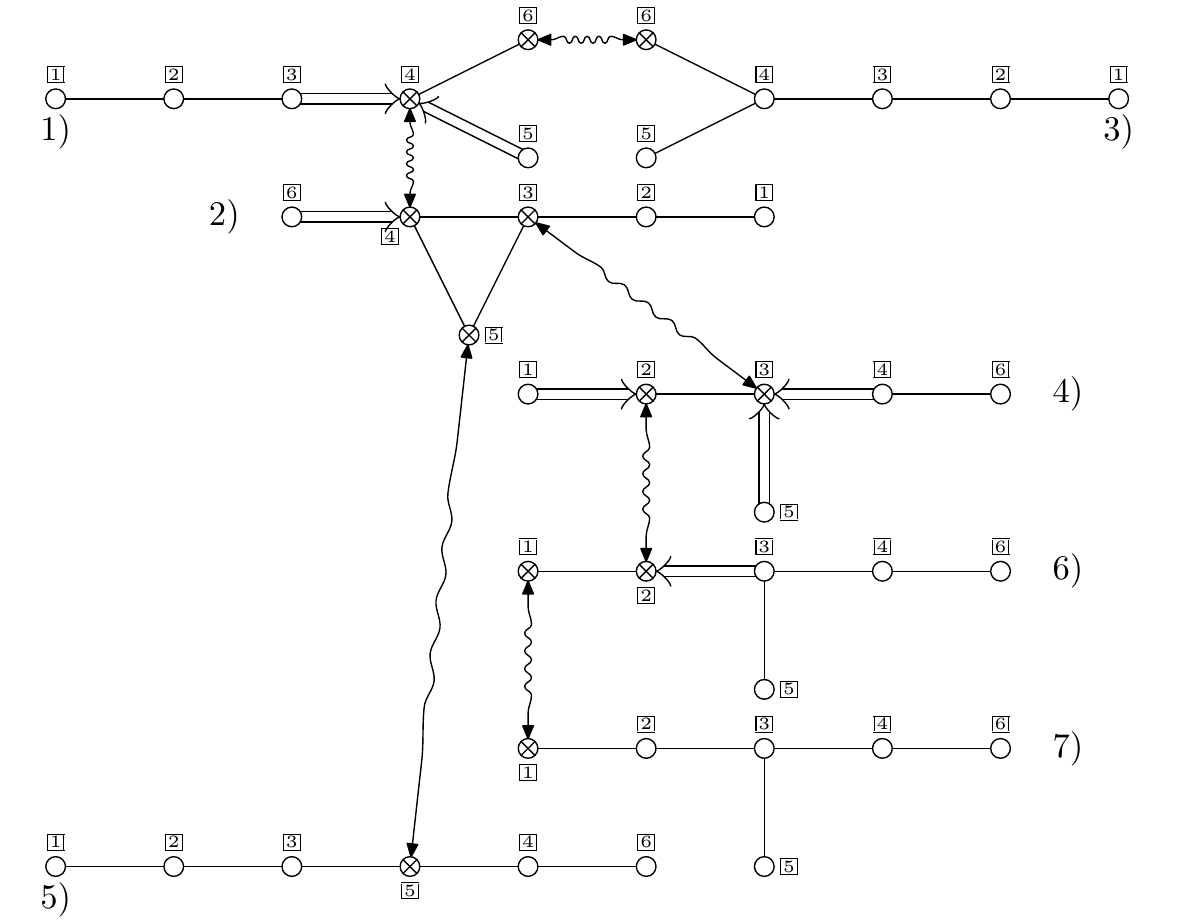}
\end{figure}
\begin{gather*} \begin{pmatrix}
 -&-&-&2&-&3 \\
 -&-&4&1&5&- \\
 -&-&-&-&-&1 \\
 -&6&2&-&-&- \\
 -&-&-&-&2&- \\
 7&4&-&-&-&- \\
 6&-&-&-&-&-
 \end{pmatrix},\quad 1)
\begin{pmatrix}
 2&-1&0&0&0&0 \\
 -1&2&-1&0&0&0 \\
 0&-1&2&-1&0&0 \\
 0&0&-2&0&-2&-1 \\
 0&0&0&-1&2&0 \\
 0&0&0&-1&0&0
 \end{pmatrix},
\\
2) \begin{pmatrix}
 2&-1&0&0&0&0 \\
 -1&2&-1&0&0&0 \\
 0&-2&0&-1&-1&0 \\
 0&0&-1&0&-1&-2 \\
 0&0&-1&-1&0&0 \\
 0&0&0&-1&0&2
 \end{pmatrix},
\quad \boxed{3)} \begin{pmatrix}
 2&-1&0&0&0&0 \\
 -1&2&-1&0&0&0 \\
 0&-1&2&-1&0&0 \\
 0&0&-1&2&-1&-1 \\
 0&0&0&-1&2&0 \\
 0&0&0&-1&0&0
 \end{pmatrix}, \\
4)  \begin{pmatrix}
 2&-1&0&0&0&0 \\
 -2&0&-1&0&0&0 \\
 0&-1&0&-2&-2&0 \\
 0&0&-1&2&0&-1 \\
 0&0&-1&0&2&0 \\
 0&0&0&-1&0&2
 \end{pmatrix},\quad
\boxed{5)} \begin{pmatrix}
 2&-1&0&0&0&0 \\
 -1&2&-1&0&0&0 \\
 0&-1&2&0&-1&0 \\
 0&0&0&2&-1&-1 \\
 0&0&-1&-1&0&0 \\
 0&0&0&-1&0&2
 \end{pmatrix}, \\
6)
 \begin{pmatrix}
 0&-1&0&0&0&0 \\
 -1&0&-2&0&0&0 \\
 0&-1&2&-1&-1&0 \\
 0&0&-1&2&0&-1 \\
 0&0&-1&0&2&0 \\
 0&0&0&-1&0&2
 \end{pmatrix},\quad
\boxed{7)} \begin{pmatrix}
 0&-1&0&0&0&0 \\
 -1&2&-1&0&0&0 \\
 0&-1&2&-1&-1&0 \\
 0&0&-1&2&0&-1 \\
 0&0&-1&0&2&0 \\
 0&0&0&-1&0&2
 \end{pmatrix}.
\end{gather*}

\newpage

\subsubsection[$\fg{}(6,6)$ of $\sdim =78|64$]{$\boldsymbol{\fg{}(6,6)}$ of $\boldsymbol{\sdim =78|64}$}

We have
$\fg(6,6)_\ev=\fo(13)$ and $\fg(6,6)_\od=\spin_{13}$.
\begin{gather*}
\begin{pmatrix}
 2&3&-&4&-&5 \\
 1&-&-&6&-&7 \\
 -&1&8&9&-&10 \\
 6&9&11&1&12&- \\
 7&10&-&-&-&1 \\
 4&-&13&2&14&- \\
 5&-&-&-&-&2 \\
 -&-&3&-&-&15 \\
 -&4&-&3&16&- \\
 -&5&15&-&-&3 \\
 13&-&4&-&-&- \\
 14&16&-&-&4&- \\
 11&17&6&-&-&- \\
 12&-&-&-&6&- \\
 -&-&10&18&-&8 \\
 -&12&19&-&9&- \\
 -&13&-&-&-&- \\
 -&-&-&15&20&- \\
 -&-&16&-&-&- \\
 -&-&-&-&18&21 \\
 -&-&-&-&-&20
 \end{pmatrix}, \quad \begin{matrix}
1) 
\begin{pmatrix}
 0&-1&0&0&0&0 \\
 -1&0&-2&0&0&0 \\
 0&-1&2&-1&0&0 \\
 0&0&-2&0&-2&-1 \\
 0&0&0&-1&2&0 \\
 0&0&0&-1&0&0
 \end{pmatrix},\\
2) \begin{pmatrix}
 0&-2&0&0&0&0 \\
 -1&2&-1&0&0&0 \\
 0&-1&2&-1&0&0 \\
 0&0&-2&0&-2&-1 \\
 0&0&0&-1&2&0 \\
 0&0&0&-1&0&0
 \end{pmatrix},\\
3) \begin{pmatrix}
 2&-1&0&0&0&0 \\
 -2&0&-1&0&0&0 \\
 0&-1&0&-2&0&0 \\
 0&0&-2&0&-2&-1 \\
 0&0&0&-1&2&0 \\
 0&0&0&-1&0&0
 \end{pmatrix},
 \end{matrix}
\\
4) \arraycolsep=1pt
\begin{pmatrix}
 0&-1&0&0&0&0 \\
 -1&0&-2&0&0&0 \\
 0&-2&0&-1&-1&0 \\
 0&0&-1&0&-1&-2 \\
 0&0&-1&-1&0&0 \\
 0&0&0&-1&0&2
 \end{pmatrix},\quad
5) \begin{pmatrix}
 0&-1&0&0&0&0 \\
 -1&0&-2&0&0&0 \\
 0&-1&2&-1&0&0 \\
 0&0&-1&2&-1&-1 \\
 0&0&0&-1&2&0 \\
 0&0&0&-1&0&0
 \end{pmatrix},\quad
6) \begin{pmatrix}
 0&-1&0&0&0&0 \\
 -1&2&-1&0&0&0 \\
 0&-2&0&-1&-1&0 \\
 0&0&-1&0&-1&-2 \\
 0&0&-1&-1&0&0 \\
 0&0&0&-1&0&2
 \end{pmatrix}, \\
7) \arraycolsep=1pt
 \begin{pmatrix}
 0&-1&0&0&0&0 \\
 -1&2&-1&0&0&0 \\
 0&-1&2&-1&0&0 \\
 0&0&-1&2&-1&-1 \\
 0&0&0&-1&2&0 \\
 0&0&0&-1&0&0
 \end{pmatrix},\quad
8) \begin{pmatrix}
 2&-1&0&0&0&0 \\
 -1&2&-1&0&0&0 \\
 0&-2&0&-1&0&0 \\
 0&0&-1&2&-2&-1 \\
 0&0&0&-1&2&0 \\
 0&0&0&-1&0&0
 \end{pmatrix}, \quad
9)
\begin{pmatrix}
 2&-1&0&0&0&0 \\
 -2&0&-1&0&0&0 \\
 0&-1&2&-1&-1&0 \\
 0&0&-1&0&-1&-2 \\
 0&0&-1&-1&0&0 \\
 0&0&0&-1&0&2
 \end{pmatrix},\\
10) \arraycolsep=1pt
 \begin{pmatrix}
 2&-1&0&0&0&0 \\
 -2&0&-1&0&0&0 \\
 0&-1&0&-2&0&0 \\
 0&0&-1&2&-1&-1 \\
 0&0&0&-1&2&0 \\
 0&0&0&-1&0&0
 \end{pmatrix}\!,\!\!\!\quad
11) \begin{pmatrix}
 0&-1&0&0&0&0 \\
 -1&2&-1&0&0&0 \\
 0&-1&0&-2&-2&0 \\
 0&0&-1&2&0&-1 \\
 0&0&-1&0&2&0 \\
 0&0&0&-1&0&2
 \end{pmatrix}\!,\!\!\! \quad
12) \arraycolsep=1pt
 \begin{pmatrix}
 0&-1&0&0&0&0 \\
 -1&0&-2&0&0&0 \\
 0&-1&2&0&-1&0 \\
 0&0&0&2&-1&-1 \\
 0&0&-1&-1&0&0 \\
 0&0&0&-1&0&2
 \end{pmatrix}\!,\\
13) \arraycolsep=1pt
 \begin{pmatrix}
 0&-1&0&0&0&0 \\
 -2&0&-1&0&0&0 \\
 0&-1&0&-2&-2&0 \\
 0&0&-1&2&0&-1 \\
 0&0&-1&0&2&0 \\
 0&0&0&-1&0&2
 \end{pmatrix}\!,\!\!\!\quad
14) \begin{pmatrix}
 0&-1&0&0&0&0 \\
 -1&2&-1&0&0&0 \\
 0&-1&2&0&-1&0 \\
 0&0&0&2&-1&-1 \\
 0&0&-1&-1&0&0 \\
 0&0&0&-1&0&2
 \end{pmatrix}\!,\!\!\! \quad
15) \begin{pmatrix}
 2&-1&0&0&0&0 \\
 -1&2&-1&0&0&0 \\
 0&-2&0&-1&0&0 \\
 0&0&-1&0&-2&-2 \\
 0&0&0&-1&2&0 \\
 0&0&0&-1&0&0
 \end{pmatrix}\!,\\
16) \arraycolsep=1pt
 \begin{pmatrix}
 2&-1&0&0&0&0 \\
 -2&0&-1&0&0&0 \\
 0&-1&0&0&-2&0 \\
 0&0&0&2&-1&-1 \\
 0&0&-1&-1&0&0 \\
 0&0&0&-1&0&2
 \end{pmatrix}\!,\!\!\!\quad
\boxed{17)} \begin{pmatrix}
 2&-1&0&0&0&0 \\
 -1&0&-2&0&0&0 \\
 0&-1&2&-1&-1&0 \\
 0&0&-1&2&0&-1 \\
 0&0&-1&0&2&0 \\
 0&0&0&-1&0&2
 \end{pmatrix}\!, \!\!\!\quad
18) \begin{pmatrix}
 2&-1&0&0&0&0 \\
 -1&2&-1&0&0&0 \\
 0&-1&2&-1&0&0 \\
 0&0&-2&0&-1&-1 \\
 0&0&0&-1&0&-1 \\
 0&0&0&-1&-1&2
 \end{pmatrix}\!,\\
\boxed{19)} \arraycolsep=1pt
 \begin{pmatrix}
 2&-1&0&0&0&0 \\
 -1&2&-1&0&0&0 \\
 0&-2&0&0&-1&0 \\
 0&0&0&2&-1&-1 \\
 0&0&-1&-2&2&0 \\
 0&0&0&-1&0&2
 \end{pmatrix}\!, \!\!\!\quad
20) \begin{pmatrix}
 2&-1&0&0&0&0 \\
 -1&2&-1&0&0&0 \\
 0&-1&2&-1&0&0 \\
 0&0&-1&2&-1&0 \\
 0&0&0&-1&0&-1 \\
 0&0&0&0&-1&0
 \end{pmatrix}\!,\!\!\!\quad
\boxed{21)} \begin{pmatrix}
 2&-1&0&0&0&0 \\
 -1&2&-1&0&0&0 \\
 0&-1&2&-1&0&0 \\
 0&0&-1&2&-1&0 \\
 0&0&0&-2&2&-1 \\
 0&0&0&0&-1&0
 \end{pmatrix}\!.
\end{gather*}
\begin{figure}[h!]\centering
\includegraphics[scale=1]{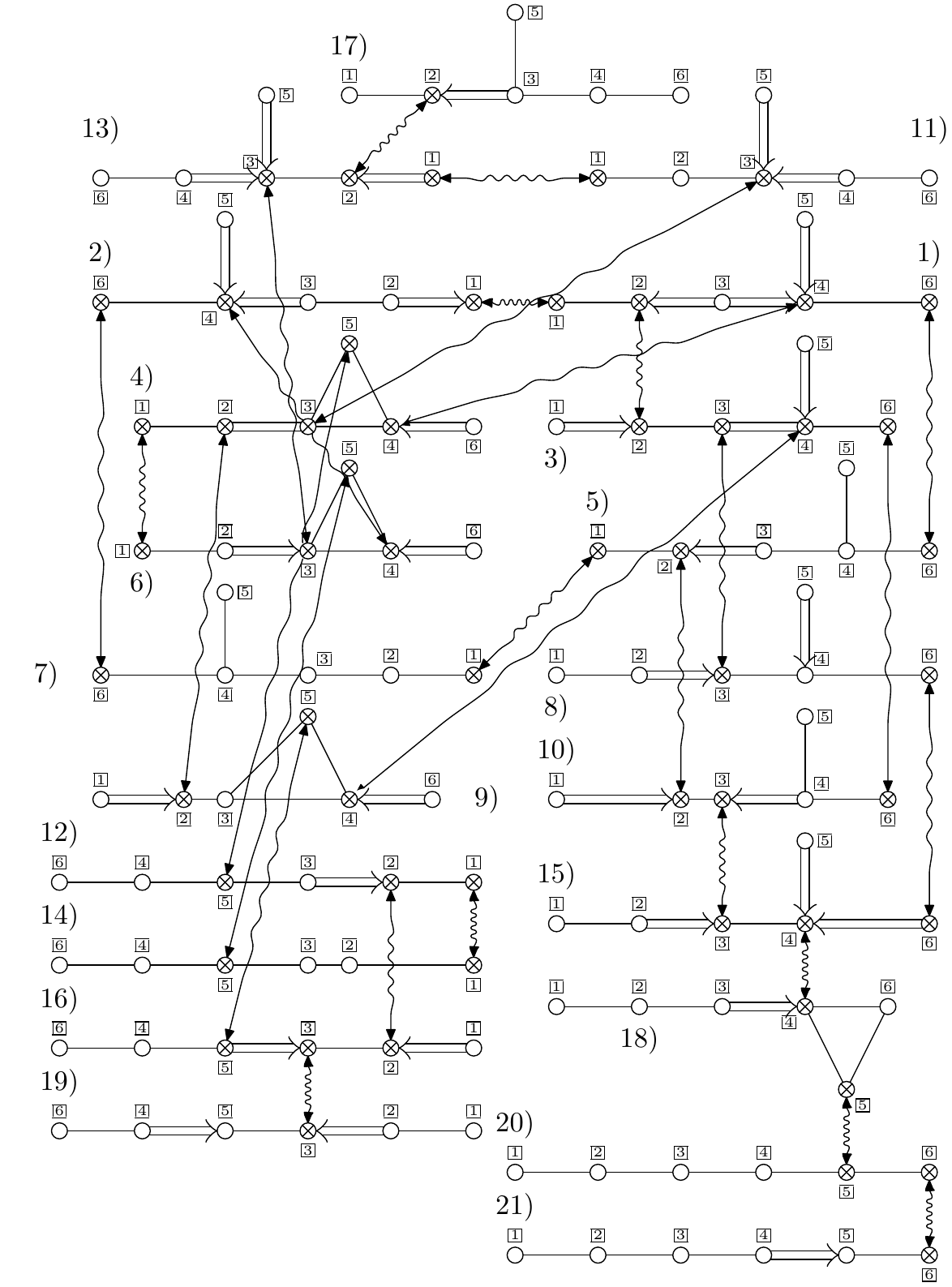}
\end{figure}

\newpage

\subsubsection[$\fg{}(8,6)$ of $\sdim =133|56$]{$\boldsymbol{\fg{}(8,6)}$ of $\boldsymbol{\sdim =133|56}$}

We have
$\fg(8,6)_\ev=\fe(7)$ and $\fg(8,6)_\od=R(\pi_1)$.
\begin{figure}[h!]\centering
\includegraphics[angle=90,scale=1.05]{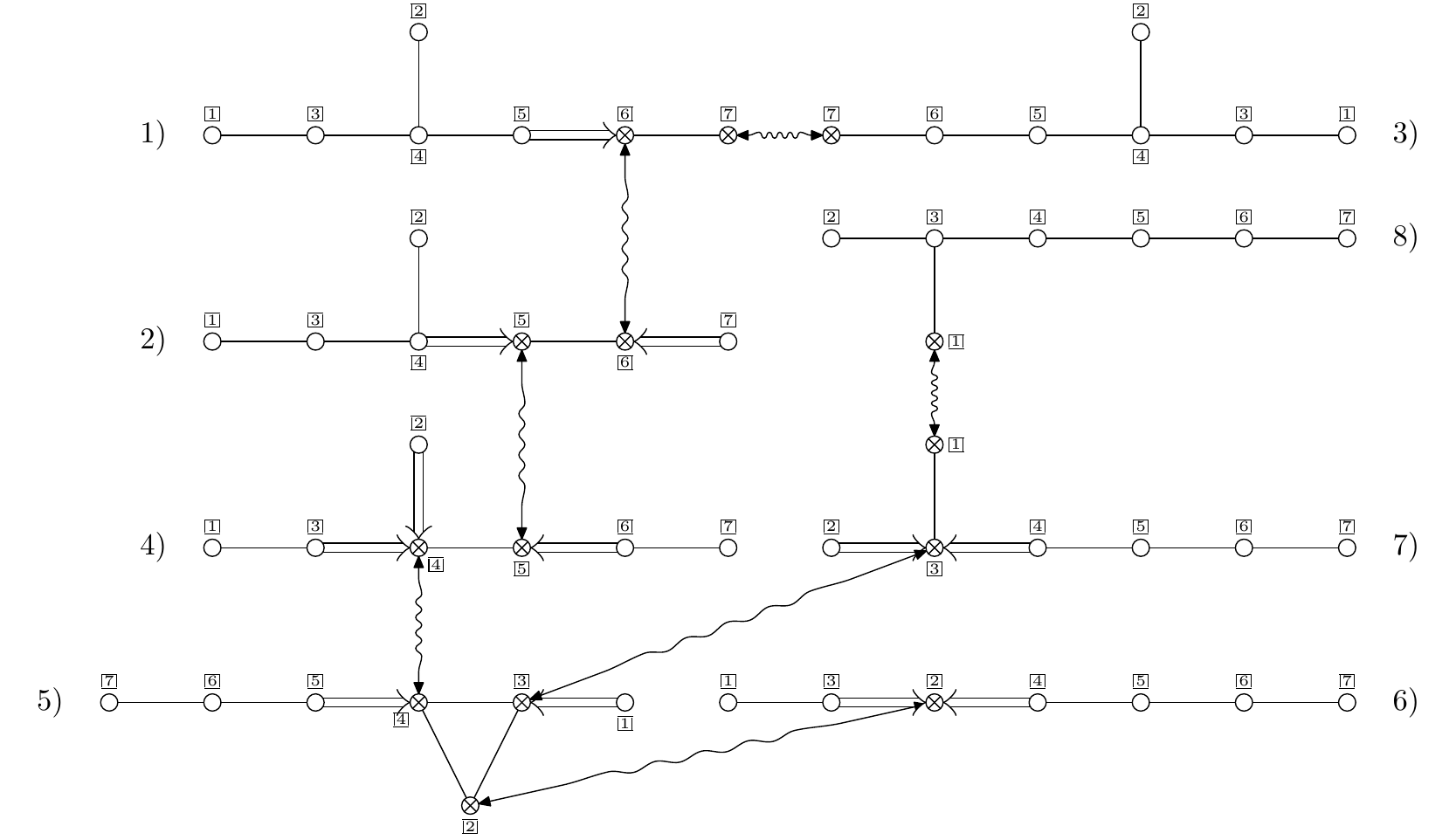}
\end{figure}
\begin{gather*}
\begin{pmatrix}
 -&-&-&-&-&2&3 \\
 -&-&-&-&4&1&- \\
 -&-&-&-&-&-&1 \\
 -&-&-&5&2&-&- \\
 -&6&7&4&-&-&- \\
 -&5&-&-&-&-&- \\
 8&-&5&-&-&-&- \\
 7&-&-&-&-&-&-
 \end{pmatrix},
\\ 1)
\begin{pmatrix}
 2&0&-1&0&0&0&0 \\
 0&2&0&-1&0&0&0 \\
 -1&0&2&-1&0&0&0 \\
 0&-1&-1&2&-1&0&0 \\
 0&0&0&-1&2&-1&0 \\
 0&0&0&0&-2&0&-1 \\
 0&0&0&0&0&-1&0
 \end{pmatrix},\quad
2) \begin{pmatrix}
 2&0&-1&0&0&0&0 \\
 0&2&0&-1&0&0&0 \\
 -1&0&2&-1&0&0&0 \\
 0&-1&-1&2&-1&0&0 \\
 0&0&0&-2&0&-1&0 \\
 0&0&0&0&-1&0&-2 \\
 0&0&0&0&0&-1&2
 \end{pmatrix}, \\
\boxed{3)} \arraycolsep=5pt
 \begin{pmatrix}
 2&0&-1&0&0&0&0 \\
 0&2&0&-1&0&0&0 \\
 -1&0&2&-1&0&0&0 \\
 0&-1&-1&2&-1&0&0 \\
 0&0&0&-1&2&-1&0 \\
 0&0&0&0&-1&2&-1 \\
 0&0&0&0&0&-1&0
 \end{pmatrix},\quad
4) \begin{pmatrix}
 2&0&-1&0&0&0&0 \\
 0&2&0&-1&0&0&0 \\
 -1&0&2&-1&0&0&0 \\
 0&-2&-2&0&-1&0&0 \\
 0&0&0&-1&0&-2&0 \\
 0&0&0&0&-1&2&-1 \\
 0&0&0&0&0&-1&2
 \end{pmatrix},\\
5) \arraycolsep=5pt
 \begin{pmatrix}
 2&0&-1&0&0&0&0 \\
 0&0&-1&-1&0&0&0 \\
 -2&-1&0&-1&0&0&0 \\
 0&-1&-1&0&-2&0&0 \\
 0&0&0&-1&2&-1&0 \\
 0&0&0&0&-1&2&-1 \\
 0&0&0&0&0&-1&2
 \end{pmatrix},\quad
\boxed{6)} \begin{pmatrix}
 2&0&-1&0&0&0&0 \\
 0&0&-2&-2&0&0&0 \\
 -1&-1&2&0&0&0&0 \\
 0&-1&0&2&-1&0&0 \\
 0&0&0&-1&2&-1&0 \\
 0&0&0&0&-1&2&-1 \\
 0&0&0&0&0&-1&2
 \end{pmatrix},
\\
7) \arraycolsep=5pt
 \begin{pmatrix}
 0&0&-1&0&0&0&0 \\
 0&2&-1&0&0&0&0 \\
 -1&-2&0&-2&0&0&0 \\
 0&0&-1&2&-1&0&0 \\
 0&0&0&-1&2&-1&0 \\
 0&0&0&0&-1&2&-1 \\
 0&0&0&0&0&-1&2
 \end{pmatrix},\quad
\boxed{8)}\begin{pmatrix}
 0&0&-1&0&0&0&0 \\
 0&2&-1&0&0&0&0 \\
 -1&-1&2&-1&0&0&0 \\
 0&0&-1&2&-1&0&0 \\
 0&0&0&-1&2&-1&0 \\
 0&0&0&0&-1&2&-1 \\
 0&0&0&0&0&-1&2
 \end{pmatrix}.
\end{gather*}

\subsection[The Elduque superalgebra $\fel(5; 3)$: Systems of simple
roots]{The Elduque superalgebra $\boldsymbol{\fel(5; 3)}$: Systems of simple
roots}\label{SSel3ssr}

Its superdimension is $39|32$; the even part
is $\fel(5; 3)_\ev=\fo(9)\oplus\fsl(2)$ and its odd part is
irreducible: $\fel(5; 3)_\od=R(\pi_4)\otimes \id$.

The following are all its Cartan matrices:
\begin{figure}[h!]\centering
\parbox{.75\linewidth}{\includegraphics{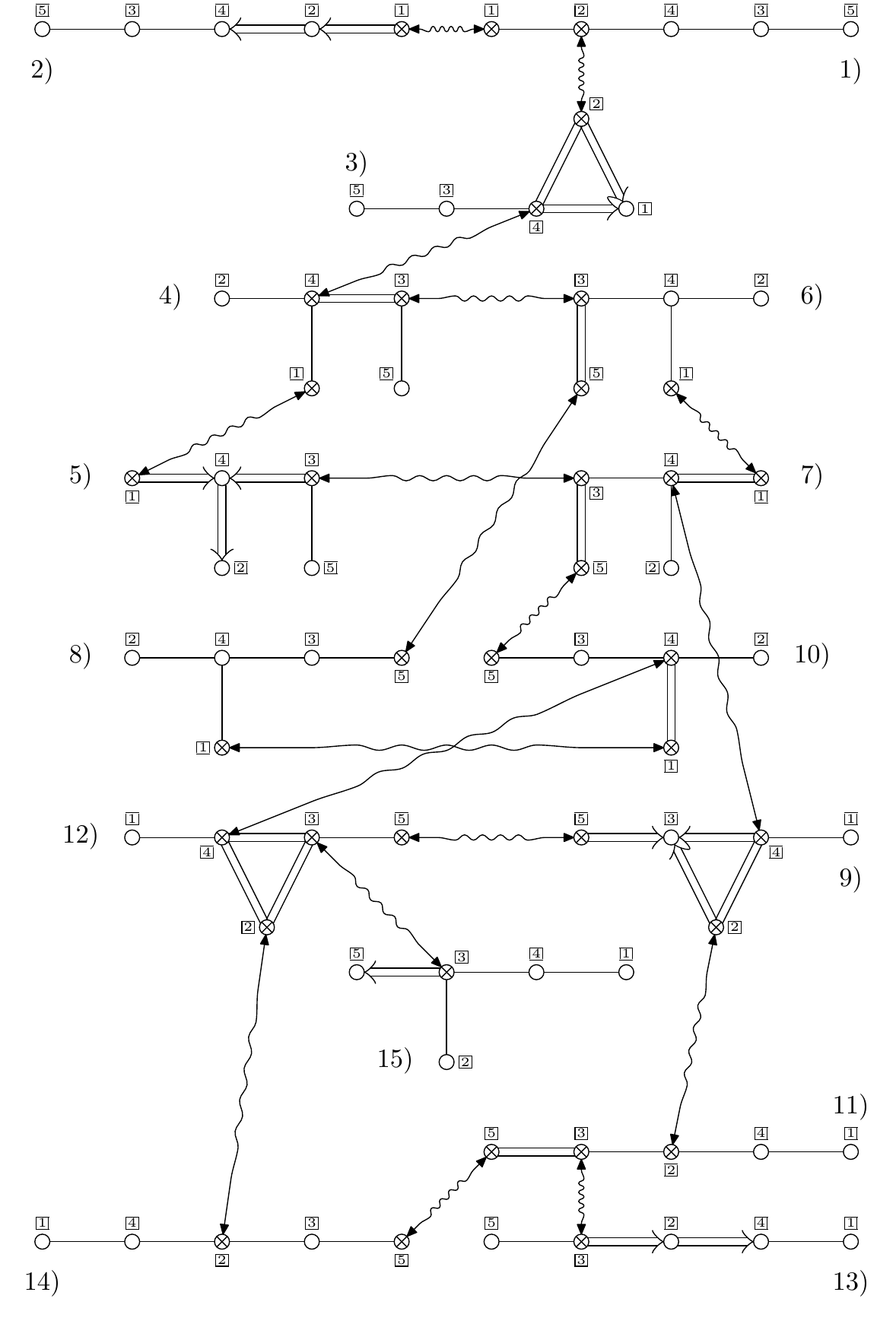}}\vspace{-3mm}
\end{figure}
\begin{gather*}
{\arraycolsep=3pt \begin{pmatrix}
2 & 3 & - & -& -\\[-1.5pt]
1 & - & - & -& -\\[-1.5pt]
- & 1 & - & 4& -\\[-1.5pt]
5 & - & 6 & 3& -\\[-1.5pt]
4 & - & 7 & -& -\\[-1.5pt]
7 & - & 4 & -& 8\\[-1.5pt]
6 & - & 5 & 9 & 10\\[-1.5pt]
10 & - & - & -& 6\\[-1.5pt]
- & 11& - &7& 12\\[-1.5pt]
8 & - & - & 12& 7\\[-1.5pt]
- & 9 & 13 & -& 14\\[-1.5pt]
- &14 & 15 & 10& 9\\[-1.5pt]
- & - & 11 & -& -\\[-1.5pt]
- & 12 & - & -& 11\\[-1.5pt]
- & - & 12 & -& -
\end{pmatrix}},\quad {\arraycolsep=3pt \begin{matrix}
1) \;
\begin{pmatrix}
 0 & -1 & 0 & 0 & 0 \\
 -1 & 0 & 0 & -1 & 0 \\
 0 & 0 & 2 & -1 & -1 \\
 0 & -1 & -1 & 2 & 0 \\
 0 & 0 & -1 & 0 & 2
\end{pmatrix},
\quad \boxed{2)} \;
\begin{pmatrix}
 0 &-2 & 0 & 0 & 0 \\
 -1 & 2 & 0 & -2 & 0 \\
 0 & 0 & 2 & -1 & -1 \\
 0 & -1 & -1 & 2 & 0 \\
 0 & 0 & -1 & 0 & 2
\end{pmatrix},
\\
 3) \;
\begin{pmatrix}
 2 & -1 & 0 & -1 & 0 \\
-2 & 0 & 0 &-2 & 0 \\
 0 & 0 & 2 & -1 & -1 \\
-2 &-2 & -1 & 0 & 0 \\
 0 & 0 & -1 & 0 & 2
\end{pmatrix},
\quad 4) \;
\begin{pmatrix}
 0 & 0 & 0 & -1 & 0 \\
 0 & 2 & 0 & -1 & 0 \\
 0 & 0 & 0 &-2 & -1 \\
 -1 & -1 &-2 & 0 & 0 \\
 0 & 0 & -1 & 0 & 2
\end{pmatrix},
\end{matrix}}
\\
{\arraycolsep=3pt 5) \;
\begin{pmatrix}
 0 & 0 & 0 &-2 & 0 \\
 0 & 2 & 0 & -1 & 0 \\
 0 & 0 & 0 &-2 & -1 \\
 -1 & -2 & -1 & 2 & 0 \\
 0 & 0 & -1 & 0 & 2
\end{pmatrix},
\quad 6) \;
\begin{pmatrix}
 0 & 0 & 0 & -1 & 0 \\
 0 & 2 & 0 & -1 & 0 \\
 0 & 0 & 0 & -1 &-2 \\
 -1 & -1 & -1 & 2 & 0 \\
 0 & 0 &-2 & 0 & 0
\end{pmatrix},
\quad
7) \;
\begin{pmatrix}
 0 & 0 & 0 &-2 & 0 \\
 0 & 2 & 0 & -1 & 0 \\
 0 & 0 & 0 & -1 &-2 \\
-2 & -1 & -1 & 0 & 0 \\
 0 & 0 &-2 & 0 & 0
\end{pmatrix}},\\
{\arraycolsep=3pt 8)
\begin{pmatrix}
 0 & 0 & 0 & -1 & 0 \\
 0 & 2 & 0 & -1 & 0 \\
 0 & 0 & 2 & -1 & -1 \\
 -1 & -1 & -1 & 2 & 0 \\
 0 & 0 & -1 & 0 & 0
\end{pmatrix},\!\quad
9)
\begin{pmatrix}
 2 & 0 & 0 & -1 & 0 \\
 0 & 0 &-2 &-2 & 0 \\
 0 & -1 & 2 & -1 & -1 \\
 -1 &-2 &-2 & 0 & 0 \\
 0 & 0 &-2 & 0 & 0
\end{pmatrix},
\quad \!10)
\begin{pmatrix}
 0 & 0 & 0 &-2 & 0 \\
 0 & 2 & 0 & -1 & 0 \\
 0 & 0 & 2 & -1 & -1 \\
-2 & -1 & -1 & 0 & 0 \\
 0 & 0 & -1 & 0 & 0
\end{pmatrix}},\\
{\arraycolsep=2.7pt 11)
\begin{pmatrix}
 2 & 0 & 0 & -1 & 0 \\
 0 & 0 & -1 & -1 & 0 \\
 0 & -1 & 0 & 0 &-2 \\
 -1 & -1 & 0 & 2 & 0 \\
 0 & 0 &-2 & 0 & 0
\end{pmatrix}\!,
\!\! \quad 12)
\begin{pmatrix}
 2 & 0 & 0 & -1 & 0 \\
 0 & 0 &-2 &-2 & 0 \\
 0 &-2 & 0 &-2 & -1 \\
 -1 &-2 &-2 & 0 & 0 \\
 0 & 0 & -1 & 0 & 0
\end{pmatrix}\!,\!\!
\quad \boxed{13)}
\begin{pmatrix}
 2 & 0 & 0 & -1 & 0 \\
 0 & 2 & -1 & -2 & 0 \\
 0 &-2 & 0 & 0 & -1 \\
 -1 & -1 & 0 & 2 & 0 \\
 0 & 0 & -1 & 0 & 2
\end{pmatrix}},\\
{\arraycolsep=3pt 14)
\begin{pmatrix}
 2 & 0 & 0 & -1 & 0 \\
 0 & 0 & -1 & -1 & 0 \\
 0 & -1 & 2 & 0 & -1 \\
 -1 & -1 & 0 & 2 & 0 \\
 0 & 0 & -1 & 0 & 0
\end{pmatrix},\quad
\boxed{15)}
\begin{pmatrix}
 2 & 0 & 0 & -1 & 0 \\
 0 & 2 & -1 & 0 & 0 \\
 0 & -1 & 0 & -1 &-2 \\
 -1 & 0 & -1 & 2 & 0 \\
 0 & 0 & -1 & 0 & 2
\end{pmatrix}}.
\end{gather*}

\section[The answer: The case where $p=2$]{The answer: The case where $\boldsymbol{p=2}$} \label{Sans2}

\underline{Simple Lie algebras}:
\begin{enumerate}\itemsep=0pt
\item[1)] The Lie algebras obtained from their Cartan matrices by reducing
modulo 2 (for $\fo(2n+1)$ one has, f\/irst of all, to divide the last
row by 2 in order to adequately normalize CM). We thus get:
\begin{enumerate}\itemsep=0pt
\item[] the CM versions of $\fsl$, namely: $\fsl(2n+1)$, and $\fgl(2n)$
whose ``simple core'' is $\fp\fsl(2n)$; in the ``second''
integer basis of $\fg(2)$ given in \cite[p.~346]{FH}, all structure
constants are integer and $\fg(2)$ becomes, after reduction modulo~2, a simple Lie algebra $\fpsl(4)$ (without Cartan matrix, as we
know);

\item[]the ``simple cores'' of the orthogonal algebras, namely, of
$\fo^{(1)}(2n+1)$ and $\foc(2n)$;

\item[]$\fe(6)$, $\fe(7)^{(1)}/\fc$, $\fe(8)$;
\end{enumerate}
\item[2)] the Weisfeiler and Kac algebras $\fwk(3;a)^{(1)}/\fc$ and
$\fwk(4;a)$.
\end{enumerate}

\underline{Simple Lie superalgebras}

In the list below the term ``super version'' of a Lie algebra
$\fg(A)$ stands for a Lie superalgebra with the ``same'' root
system as that of $\fg(A)$ but with some of the simple roots
considered odd.
\begin{enumerate}\itemsep=0pt
\item[1)] The Lie superalgebras obtained from their $p=0$ analogs that have
no $-2$ in of\/f-diagonal slots of the Cartan matrix by reducing the
structure constants modulo 2 (for $\fosp(2n+1|2m)$ one has, f\/irst of
all, to divide the last row by 2 in order to normalize CM), we thus
get
\begin{enumerate}\itemsep=0pt
\item[]
the CM versions of $\fsl$, namely: either simple $\fsl(a|a+2k+1)$ or
$\fgl(a|2k+a)$ whose ``simple core'' is $\fp\fsl(a|a+2k)$ for
$a$ odd or $\fp\fsl(a|a+2k)^{(1)}$ for $a$ even;
\end{enumerate}
\item[2)] the ortho-orthogonal algebras, namely: $\foo^{(1)}$ and $\fooc$
whose ``simple cores'' are described in Section~\ref{Soo};

\item[3)] $\fbgl(3;a)^{(1)}/\fc$ which is an analog of $\fwk(3;
a)^{(1)}/\fc$ with the ``same'' Cartan matrices but dif\/ferent
root systems;

\item[4)] the CM versions of periplectic algebras, namely: $\fpec$; these
are at the same time super versions of $\foc$; their ``simple
cores'' are described in Section~\ref{Soo};

\item[5)] a super version of $\fwk(4;a)$, namely: $\fbgl(4;a)$;

\item[6)] the super versions of $\fe(6)$, namely: $\fe(6, 1)$, $\fe(6, 6)$;

\item[7)] the super versions of $\fe(7)$, namely: $\fe(7, 1)$, $\fe(7, 6)$,
$\fe(7, 7)$ whose ``simple cores'' are described in Section~\ref{etype};

\item[8)] the super versions of $\fe(8)$, namely: $\fe(8, 1)$, $\fe(8, 8)$.
\end{enumerate}

\subsection[On the structure of $\fbgl(3;\alpha)$, $\fbgl(4;\alpha)$, and
$\fe(a,b)$]{On the structure of $\boldsymbol{\fbgl(3;\alpha)}$, $\boldsymbol{\fbgl(4;\alpha)}$, and
$\boldsymbol{\fe(a,b)}$}\label{Sstruct2}

In this section we describe the even parts $\fg_\ev$ of the new Lie
superalgebras $\fg=\fg(A)$ and their odd parts $\fg_\od$ as
$\fg_\ev$-modules. \texttt{SuperLie} enumerates the elements of the
Chevalley basis the $x_i$ (positive), starting with the generators,
then their brackets, etc., and the $y_i$ are negative root vectors
opposite to the $x_i$. Since the irreducible representations of the
Lie algebras may have neither highest nor lowest weight, observe
that the $\fg_\ev$-modules $\fg_\od$ always have both highest and
lowest weights.

\subsubsection[Notation  $\fA\oplus_c \fB$ needed to describe
$\fbgl(4; \alpha)$, $\fe(6, 6)$, $\fe(7,6)$, and $\fe(8,1)$]{Notation  $\boldsymbol{\fA\oplus_c \fB}$ needed to describe $\boldsymbol{\fbgl(4; \alpha)}$, $\boldsymbol{\fe(6, 6)}$, $\boldsymbol{\fe(7,6)}$, and $\boldsymbol{\fe(8,1)}$}\label{4cases}

This notation describes the case where $\fA$ and
$\fB$ are nontrivial central extensions of the Lie algebras $\fa$
and $\fb$, respectively, and $\fA\oplus_c \fB$~-- a nontrivial central extension of $\fa\oplus \fb$ (or, perhaps,
a more complicated $\fa\subplus \fb$) with 1-dimensional center
spanned by $c$~-- is such that the restriction of the extension of
$\fa\oplus \fb$ to $\fa$ gives $\fA$ and that to $\fb$ gives $\fB$.
(In other words, the situation resembles the (nontrivial) central
extension of the Lie algebra of derivations of the loop algebra,
namely, $\fg\otimes\Cee[t^{-1}, t]\subplus \fder(\Cee[t^{-1}, t])$,
where one central element serves both central extensions: That of
$\fg\otimes\Cee[t^{-1}, t]$ and of $\fder(\Cee[t^{-1}, t])$.)

In these four cases, $\fg(A)_\ev$ is of the form
\[\fg(B)\oplus_c\fhei(2)\simeq \fg(B)\oplus\Span(X^+,X^-),\] where the
matrix $B$ is not invertible (so $\fg(B)$ has a grading element $d$
and a central element $c$), and where $X^+$, $X^-$ and $c$ span the
Heisenberg Lie algebra $\fhei(2)$. The brackets are:
\begin{gather}
{}[\fg^{(1)}(B), X^\pm]=0,\nonumber\\
{}[d,X^\pm]=X^\pm, \qquad([d,X^\pm]=\alpha X^\pm\text{~for $\fbgl(4; \alpha)$}),\label{strange}\\
{} [X^+,X^-]=c. \nonumber
\end{gather}

The odd part of $\fg(A)$ (at least in two of the four cases)
consists of two copies of the same $\fg(B)$-module $N$, the
operators $\ad_{X^\pm}$ permute these copies, and $\ad_{X^\pm}^2=0$,
so each of the operators maps one of the copies to the other, and
this other copy to zero.

\subsubsection[$\fbgl(3;\alpha)$, where $\alpha\neq 0, 1$; $\sdim=10/8|8$]{$\boldsymbol{\fbgl(3;\alpha)}$, where $\boldsymbol{\alpha\neq 0, 1}$; $\boldsymbol{\sdim=10/8|8}$}

We consider the following Cartan matrix and the corresponding
positive root vectors (odd $|$ even)
\[
\begin{pmatrix}
 0 & 1 & 0 \\
 1 & \ev & \alpha \\
 0 & \alpha & \ev
\end{pmatrix}
\quad
\begin{array}{l}
x_1\mid x_2,\quad x_3,\\
x_4=[x_1,x_2]\mid x_5=[x_2,x_3],\quad x_6=[x_3,[x_1,x_2]],\mid\\
x_7=[[x_1,x_2],[x_2,x_3]]\mid.
 \end{array}
\]

Then $\fg_\ev\simeq\fgl(3)\oplus \mathbb{K} Z$. The $\fg_\ev$-module
$\fg_\od$ is reducible, with the two highest weight vectors, $x_7$
and $y_1$. The Cartan subalgebra of $\fgl(3)\oplus \mathbb{K}Z$ is
spanned by $\alpha h_1+h_3, h_2, h_3$ and $Z$. In this basis, the
weight of $x_7$ is $(0,1+\alpha,0,1)$. The weight of $y_1$ is
$(0,1,0,1)$, if for the grading operator we take
$(1,0,0)\in\fgl(3)$.

The lowest weight vectors of these modules are $x_1$ and $y_7$ and
their weights are $(0,1,0,1)$ and $(0,1+\alpha,0,1)$.
The module generated by $x_7$ is
$\Span\left\{x_1,x_4,x_6,x_7\right\}$.
The module generated by $y_1$ is
$\Span\left\{y_1,y_4,y_6,y_7\right\}$.

All inequivalent Cartan matrices are{\samepage
\[
 \begin{pmatrix} d_1 &\alpha &1\\ \alpha &d_2&0\\1&0&d_3\end{pmatrix}
, \qquad
\begin{pmatrix} d_1 &\alpha &1+\alpha\\ \alpha &d_2&1\\1+\alpha &1&d_3
\end{pmatrix},
\]
where $(d_1, d_2, d_3)$ is any distribution of $0$'s and $\ev$'s,
except $(\ev, \ev, \ev)$.}

\subsubsection[$\mathfrak{bgl}(4;\alpha)$, where $\alpha\neq 0, 1$, of $\sdim =18|16$]{$\boldsymbol{\mathfrak{bgl}(4;\alpha)}$, where $\boldsymbol{\alpha\neq 0, 1}$, of $\boldsymbol{\sdim =18|16}$}

We consider the following Cartan matrix and the
corresponding positive root vectors (odd $|$ even)
\[
\begin{pmatrix}
 0 & \alpha & 1 & 0 \\
 \alpha & \ev & 0 & 0 \\
 1 & 0 & \ev & 1 \\
 0 & 0 & 1 & \ev
\end{pmatrix}
\quad
\begin{array}{l}
x_1\mid x_2,\quad x_3,\quad x_4,\\
x_5=[x_1,x_2],\quad x_6=[x_1,x_3]\mid x_7=[x_3,x_4],\\
x_8=[x_3,[x_1,x_2]],\quad x_9=[x_4,[x_1,x_3]]\mid \\
 x_{11}=[[x_1,x_2],[x_3,x_4]]\mid x_{10}=[[x_1,x_2],[x_1,x_3]]\\
\mid x_{12}=[[x_1,x_2],[x_4,[x_1,x_3]]],\\
\mid x_{13}=[[x_3,[x_1,x_2]],[x_4,[x_1,x_3]]],\\
x_{14}=[[x_4,[x_1,x_3]],[[x_1,x_2],[x_1,x_3]]]\mid\\
x_{15}=[[[x_1,x_2],[x_1,x_3]],[[x_1,x_2],[x_3,x_4]]]\mid.
 \end{array}
\]
In this case $\fg_\ev\simeq\fgl(4) \oplus_c \fhei(2)$ see
Section~\ref{4cases} with commutation relations (\ref{strange}). The
$\fg_\ev$-module $\fg_\od$ is irreducible: $\fg_\od\simeq N\otimes
\id$, where $\id$ is the standard 2-dimensional $\fhei(2)$-module
and $N$ is an 8-dimensional $\fgl(4)$-module.

The highest weight vector $x_{15}$ has weight $(\alpha,
0,0,0,\alpha)$ with respect to \[
c=h_2,\qquad d=h_1,\qquad
H_1=h_3,\qquad H_2=h_3,\qquad H_3=h_2+h_3,\] where the $h_i$'s are the
Chevalley generators of the Cartan subalgebra of $\fbgl(4;\alpha)$.
The lowest weight vector is $y_{15}$ of the same weight as $x_{15}$.

All inequivalent Cartan matrices of $\fbgl(4; \alpha)$ are
\[
\begin{pmatrix} d_1 &\alpha &0&0\\
\alpha &d_2&1&0\\
0&1&d_3&1\\
0&0&1&d_4 \end{pmatrix} ,\qquad
\begin{pmatrix} d_1 &1 &1+\alpha&0\\
1 &d_2& \alpha & 0\\
\alpha+1& \alpha &d_3&\alpha\\
0&0&\alpha&d_4\end{pmatrix}
 ,\qquad
\begin{pmatrix} d_1 &\alpha & 0 &0\\
\alpha &d_2& \alpha+1 & 0\\
0& \alpha+1 &d_3&1\\
0&0&1&d_4 \end{pmatrix},
\]
where $\{d_1,d_2,d_3,d_4\}$ is any distribution of $0$'s and
$\ev$'s, except $\{\ev,\ev,\ev,\ev\}$.

\sssbegin{Proposition} [cf.~(\ref{wkiso}) and (\ref{osp42symm})]\qquad

$1)$ We have
\begin{equation}\label{swkiso}\begin{array}{l}
\renewcommand{\arraystretch}{1.4}
\fbgl(3;a)\simeq \fbgl(3;a')\Longleftrightarrow
a'=\displaystyle\frac{\alpha a+\beta}{\gamma a+\delta},\qquad
\text{where $\begin{pmatrix}\alpha&\beta\\
\gamma &\delta\end{pmatrix}\in \SL(2; \Zee/2)$}\\
\fbgl(4;a)\simeq \fbgl(4;a')\Longleftrightarrow
a'=\displaystyle\frac{1}{a}.
\end{array}
\end{equation}

$2)$ The $2|4$-structures on $\fbgl(3;a)$ and $\fbgl(4;a)$ are given
by the same formulas $\eqref{2strwk3}$, $\eqref{2strwk31}$,
$\eqref{2strwk4}$ as for $\fwk(3;a)$ and $\fwk(4;a)$ with the
following amendment:
\begin{equation}\label{swkiso1}
\begin{array}{ll}(e_\alpha^{\pm})^{[2]}=0&\text{for any root
vector $e_\alpha$ even},\\
(e_\alpha^{\pm})^{[4]}=((e_\alpha^{\pm})^2)^{[2]}&\text{for any root
vector $e_\alpha$  odd}.\end{array}\end{equation}
\end{Proposition}

\subsubsection[The $\fe$-type superalgebras]{The $\fe$-type superalgebras}\label{etype}

 {\bf Notation.}
The $\fe$-type superalgebras will be denoted by (one of) their
simplest Dynkin diagrams, i.e., $\fe(n, i)$ denotes the Lie
superalgebra whose diagram is of the same shape as that of the Lie
algebra $\fe(n)$ but with the only~-- $i$th~-- node $\motimes$. This, and other ``simplest'',
Cartan matrices are boxed. We enumerate the nodes of the Dynkin
diagram of $\fe(n)$ as in \cite{Bou, OV}: We f\/irst enumerate the
nodes in the row corresponding to $\fsl(n)$ (from the end-point of
the ``longest'' twig towards the branch point and further on
along the second long twig), and the $n$th node is the end-point of
the shortest ``twig''.

\paragraph{$\boldsymbol{\fe(6, 1)\simeq \fe(6,5)}$ of $\boldsymbol{\sdim = 46|32}$.} We have
$\fg_\ev\simeq \fo\fc(2; 10)\oplus\Kee Z$ and $\fg_\od$ is a
reducible module of the form $R(\pi_{4})\oplus R(\pi_{5})$ with the
two highest weight vectors \[x_{36}=[[[x_4,x_5],[x_6,[x_2,x_3]]],
[[x_3,[x_1,x_2]],[x_6,[x_3,x_4]]]]\] and $y_5$. Denote the basis
elements of the Cartan subalgebra by $Z$, $h_1$, $h_2$, $h_3$, $h_4$, $h_6$.
The weights of $x_{36}$ and $y_5$ are respectively, $(0,0,0,0,0,1)$
and $(0,0,0,0,1,0)$. The module generated by $x_{36}$ gives all odd
positive roots and the module generated by $y_5$ gives all odd
negative roots.

\paragraph{$\boldsymbol{\fe(6, 6)}$ of $\boldsymbol{\sdim = 38|40}$.}\label{3cases}
In this case, $\fg(B)\simeq \fgl(6)$, see Section~\ref{4cases}. The module
$\fg_\od$ is irreducible with the highest weight vector
\[x_{35}=[[[x_3,x_6], [x_4,[x_2,x_3]]], [[x_4,x_5], [x_3,[x_1,x_2]]]]
\text{~of weight $(0,0,1,0,0,1)$.}
\]
We consider the highest weight with respect to the elements
\[
h_1:=E_{11}-E_{22},\ \dots,\ h_5:=E_{55}-E_{66},\ h_6:=E_{11}+E_{22}+E_{33}\]
in $\fgl(6)$. We can equally well set
\[h_6:=E_{11}+E_{22}+E_{33}+ac \text{~for any $a\in\Kee$, where $c$ is the
non-zero central element of $\fg(B)$}\] but in our choice of
$h_6=E_{11}+E_{22}+E_{33}$, we have $M=\bigwedge^3(\id)$, as a
$\fgl(6)$-module (note that to write $M=R(\pi_3)$ is not enough
since this only describes $M$ as an $\fsl(6)$-module).

\paragraph{$\boldsymbol{\fe(7,1)}$ of $\boldsymbol{\sdim=80/78|54}$.} Since the Cartan matrix of
this Lie superalgebra is of rank 6, a~grading operator $d_1$ should
be (and is) added. Now if we take $d_1=(1,0,0,0,0,0,0)$, then
$\fg_\ev\simeq (\fe(6)\oplus\Kee z) \oplus \Kee I_0$. The Cartan
subalgebra is spanned by $h_1+h_3+h_7$, $h_2$, $h_3$, $h_4$, $h_5$, $h_6$, $h_7$
and $d_1$. We see that $\fg_\od$ has the two highest weight vectors:
\[
x_{63}=[[[[x_2,x_3],[x_4,x_7]], [[x_3,x_4],[x_5,x_6]]], [[[x_4,x_7],
[x_5,x_6]], [[x_4,x_5],[x_3,[x_1,x_2]]]]]\] and $y_1$. Their
respective weights (if we take $d_1=(1,0,0,0,0,0,0)$) are
$(0,0,0,0,0,1,0,1)$ and $(0,1,0,0,0,0,0,0)$. The module generated by
$x_{63}$ gives all odd positive roots and the module generated by
$y_1$ gives all odd negative roots.

\paragraph{$\boldsymbol{\fe(7,6)$ of $\sdim=70/68|64}$.} We are in the
same situation as before (Section~\ref{4cases}). We have $\fg(B)\simeq
\foc(1;12)\subplus \Kee I_0$. Note that in this case
$\text{size}(B)-\rk(B)=2$, so the center of $\fg(B)$ is
2-dimensional, and $\dim\fg(B)-\dim\fg^{(1)}(B)=2$. So we should be
a bit more specif\/ic than in~(\ref{strange}); namely, we have
\begin{gather*}
[\foc(1;12),X^\pm]=0,\\
[I_0, X^\pm]=X^\pm,\\
[X^+,X^-]=h_1+h_3+h_5~~\text{(which corresponds to $1_{12}$ in
$\foc(1;12)$).}
\end{gather*}
The module $\fg_\od$ is irreducible
with the highest weight vector
\[
x_{62}=[[[x_7,[x_5,[x_3,x_4]]],[[x_1,x_2], [x_3,x_4]]],
[[[x_2,x_3],[x_4,x_5]], [[x_4,x_7],[x_5,x_6]]]].
\]
The Cartan subalgebra is spanned by $h_1+h_3+h_5$, $h_1$, $h_2$, $h_3$,
$h_4$, $h_7$ and also $h_6$ and $d_1$. The weight of $x_{62}$ is
$(1,0,0,0,0,0,1,0)$. The highest weight vector of $\fg_\od$ is the
highest weight vector of one of the copies of the $\fg(B)$-module
$N$, see Section~\ref{4cases}, so the highest weight of $N$ is the same as
the highest weight of $\fg_\od$. (Of course, this is true for the
other two similar cases as well; in the case of $\fe(6,6)$, we used
Lebedev's choice~-- another basis of $\fh$~-- and expressed the
weight with respect to it.)

\paragraph{$\boldsymbol{\fe(7,7)}$ of $\boldsymbol{\sdim=64/62|70}$.} Since the Cartan matrix of this Lie
superalgebra is of rank 6, a grading operator $d_1$ should be (and
is) added. Then $\fg_\ev\simeq \fgl(8)$. The module $\fg_\od$ has
the two highest weight vectors:
\[
x_{58}=[[[x_3,[x_1,x_2]], [x_6,[x_4,x_5]]], [[x_7,[x_3,x_4]],
[[x_2,x_3], [x_4,x_5]]]]\] and $y_7$. The Cartan subalgebra is
spanned by $h_1$, $h_2$, $h_3$, $h_4$, $h_5$, $h_6$ and also $h_1+h_3+h_7$ and
$d_1$. The weight of $x_{58}$ with respect to these elements of the
Cartan subalgebra is $(0,0,1,0,0,0,0,1)$ and the weight of $y_7$ is
$(0,0,0,1,0,0,0,1)$. The module generated by $x_{58}$ gives all odd
positive roots and the module generated by $y_7$ gives all odd
negative roots.

\paragraph{$\boldsymbol{\fe(8, 1)}$ of $\boldsymbol{\sdim=136|112}$.} We have (cf.\ Section~\ref{4cases})
$\fg(B)\simeq \fe(7)$. (Recall that, in our notation, $\fe(7)^{(1)}$
has a center but not the grading operator, see Section ``Warning''~\ref{warn}.) The Cartan subalgebra is spanned by
$h_2+h_4+h_8$ and $h_1$, $h_2$, $h_3$, $h_4$, $h_5$, $h_6$, $h_7$.
The $\fg_\ev$-module $\fg_\od$ is irreducible with the highest
weight vector:
\begin{gather*}
x_{119}= [[[[x_4, [x_2, x_3]], [[x_5, x_8], [x_6, x_7]]]
 , [[x_8, [x_4, x_5]],  [[x_3, x_4], [x_5, x_6]]]] ,\\
\phantom{x_{119}=}{} \  [[[x_7,[x_5,x_6]],[[x_1,x_2],[x_3,x_4]]]
 ,[[x_8,[x_5,x_6]],
 [[x_2,x_3],[x_4,x_5]]]]]
\end{gather*}
of weight $(1,1,0,0,0,0,0,1)$ and one lowest weight vector $y_{119}$
whose expression is as above the $x$'s changed by the $y$'s, of the
same weight as that of $x_{119}$. (Again, the highest weight of the
$\fg(B)$-module $N$, see Section~\ref{4cases}, is the same as the highest
weight of $\fg_\od$.)

\paragraph{$\boldsymbol{\fe(8, 8)}$ of $\boldsymbol{\sdim=120|128}$.} In the $\Zee$-grading
with the 1st CM with $\deg e_8^\pm=\pm 1$ and $\deg e_i^\pm=0$ for
$i\neq 8$, we have $\fg_0=\fgl(8)=\fgl(V)$. There are dif\/ferent
isomorphisms between $\fg_0$ and $\fgl(8)$; using the one where
$h_i=E_{i,i}+E_{i+1,i+1}$ for all $i=1,\dots, 7$, and
$h_8=E_{6,6}+E_{7,7}+E_{8,8}$, we see that, as modules over
$\fgl(V)$,
\begin{alignat*}{4}
& \fg_1=\bigwedge^5V^*, \qquad &&\fg_2=\bigwedge^6 V, \qquad & &\fg_3=V, &\\
& \fg_{-1}=\bigwedge^5V, \qquad &&\fg_{-2}=\bigwedge^6 V^*, \qquad &&\fg_{-3}=V^*.&
\end{alignat*}
We can also set $h_8=E_{1,1}+E_{2,2}+E_{3,3}+E_{4,4}+E_{5,5}$. Then
we get
\begin{alignat*}{4}
& \fg_1=\bigwedge^3V, \qquad &&\fg_2=\bigwedge^6 V, \qquad &&\fg_3=\bigwedge^7V^*, &\\
& \fg_{-1}=\bigwedge^3V^*, \qquad &&\fg_{-2}=\bigwedge^6 V^*, \qquad & &\fg_{-3}=\bigwedge^7V.&
\end{alignat*}

The algebra $\fg_\ev$ is isomorphic to
$\fo_\Pi^{(2)}(16)\subplus\Kee d$, where
$d=E_{6,6}+\dots+E_{13,13}$, and $\fg_\od$ is an irreducible
$\fg_\ev$-module with the highest weight the highest weight element
$x_{120}$ of weight $(1,0,\dots,0)$ with respect to $h_1,\dots,h_8$;
$\fg_\od$ also possesses a lowest weight vector.

\subsection[Systems of simple roots of the $\fe$-type Lie superalgebras]{Systems of simple roots of the $\boldsymbol{\fe}$-type Lie superalgebras}\label{Sssr-e}

\sssbegin{Remark} Observe that if $p=2$ and the Cartan matrix has no
parameters, the ref\/lections do not change the shape of the Dynkin
diagram. Therefore, for the $\fe$-superalgebras, it suf\/f\/ices to list
distributions of parities of the nodes in order to describe the
Dynkin diagrams. Since there are tens and even hundreds of diagrams
in these cases, this possibility saves a lot of space, see the lists
of all inequivalent Cartan matrices of the $\fe$-type Lie
superalgebras.
\end{Remark}

\subsubsection[$\fe(6, 1)\simeq \fe(6,5)$ of $\sdim 46|32$]{$\boldsymbol{\fe(6, 1)\simeq \fe(6,5)}$ of $\boldsymbol{\sdim 46|32}$}

All inequivalent
Cartan matrices are as follows (none of the matrices corresponding
to the symmetric pairs of Dynkin diagrams is excluded but are placed
one under the other for clarity, followed by three symmetric
diagrams):\tiny
\begin{equation*}\label{e65}
\begin{array}{llllllllllll}
\boxed{1)}&000010&3)&010001&5)&100110&7)&000011&9)&000110&11)&000111\\
\boxed{2)}&100000&4)&000101&6)&110010&8)&100001&10)&110000&12)&110001\\
13)&111001&15)&101001&17)&011000&19)&101100&21)&011001&23)&011110\\
14)&001111&16)&001011&18)&001100&20)&011010&22)&001101&24)&111100\\
25)&010100&26)&100010&27)&110110&&&&&&\\
 \end{array}
\end{equation*}

\normalsize

\subsubsection[$\fe(6, 6)$ of $\sdim = 38|40$]{$\boldsymbol{\fe(6, 6)}$ of $\boldsymbol{\sdim = 38|40}$} All inequivalent Cartan
matrices are as follows: \tiny
\begin{equation*}\label{e66}
\begin{array}{llllllllllllllll}
\boxed{1)}&000001&\boxed{2)}&000100&\boxed{3)}&001000
&\boxed{4)}&010000&5)&011011&6)&101110&7)&111110\\
8)&011100&9)&101111&10)&011101&11)&101010&12)&111101&
13)&010110&14)&101011\\
15)&110011&16)&001001&17)&011111&18)&110100&19)&010011&20)&101000&
21)&111011\\
22)&001010&23)&100011&24)&110101&25)&001110&26)&111000&27)&010010&
28)&100111\\
29)&100100&30)&110111&31)&100101&32)&111010&33)&010101&34)&010111&35)&101101\\
36)&111111\\
\end{array}
\end{equation*}

\normalsize

\subsubsection[$\fe(7,1)$ of $\sdim=80/78|54$]{$\boldsymbol{\fe(7,1)}$ of $\boldsymbol{\sdim=80/78|54}$} All inequivalent Cartan
matrices are as follows:\tiny
\begin{equation*}\label{e71}
\begin{array}{llllllllllllll}
\boxed{1)}&1000000&2)&1000010&3)&1000110&
4)&1001100&25)&0110000&26)&0110010&27)&0110110\\
5)&1010001&6)&1011001& 7)&1100000&8)&1100010&21)&0011010&
22)&0011110&23)&0100001\\
9)&1100110&
10)&1101100&11)&1110001&12)&1111001&17)&0001101&18)&0001111&
19)&0010100\\
13)&0000011&14)&0000101&15)&0000111&
16)&0001011&28)&0111100&24)&0101001&20)&0011000\\
\end{array}
\end{equation*}

\normalsize

\subsubsection[$\fe(7,6)$ of $\sdim=70/68|64$]{$\boldsymbol{\fe(7,6)}$ of $\boldsymbol{\sdim=70/68|64}$} All inequivalent Cartan
matrices are as follows:\tiny
\begin{equation*}\label{e76}
\begin{array}{llllllllllll}
\boxed{1)}&0000010&\boxed{2)}&0000100&3)&0000110&
\boxed{4)}&0001000&62)&1111100&63)&1111110\\
5)&0001010&6)&0001100&
7)&0001110&8)&0010001&60&1111000& 61)&1111010\\
9)&0010011&
10)&0010101&11)&0010111&12)&0011001& 58)&1110100&59)&1110110\\
13)&0011011&14)&0011101&15)&0011111&
\boxed{16)}&0100000&56)&1110000& 57)&1110010\\
17)&0100010&18)&0100100&
19)&0100110&20)&0101000&54)&1101101& 55)&1101111\\
21)&0101010&
22)&0101100&23)&0101110&24)&0110001& 52)&1101001& 53)&1101011\\
25)&0110011&26)&0110101&27)&0110111&
28)&0111001&50)&1100101&51)&1100111\\
29)&0111011&30)&0111101&
31)&0111111&32)&1000001&48)&1100001& 49)&1100011\\
33)&1000011&
34)&1000101&35)&1000111&36)&1001001& 46)&1011100&47)&1011110 \\
37)&1001011&38)&1001101&39)&1001111&
40)&1010000&44)&1011000& 45)&1011010\\
41)&1010010&42)&1010100& 43)&1010110&&
\end{array}
\end{equation*}

\normalsize

\subsubsection[$\fe(7,7)$ of $\sdim=64/62|70$]{$\boldsymbol{\fe(7,7)}$ of $\boldsymbol{\sdim=64/62|70}$} All inequivalent Cartan
matrices are as follows:\tiny
\begin{equation*}\label{e77}
\begin{array}{llllllllllll}
\boxed{1)}&0000001&2)& 0001001&\boxed{3)}&0010000&4)&0010010&
34)&1111011&35)&1111101\\
5)&0010110&6)&0011100&
7)&0100011&8)&0100101&32)&1110101& 33)&1110111\\
9)&0100111&
10)&0101011&11)&0101101&12)&0101111&30)&1101110& 31)&1110011\\
13)&0110100&14)&0111000&15)&0111010&
16)&0111110& 28)&1101000& 29)&1101010\\
17)&1000100&18)&1001000&
19)&1001010&20)&1001110&26)&1011111&27)&1100100\\
21)&1010011& 22)&1010101&23)&1010111&24)&1011011& 25)&1011101
\end{array}
\end{equation*}

\normalsize

\newpage

 \subsubsection[$\fe(8,1)$ of $\sdim=136|112$]{$\boldsymbol{\fe(8,1)}$ of $\boldsymbol{\sdim=136|112}$} All inequivalent
Cartan matrices are as follows:\tiny
\begin{equation*}\label{81}
\begin{array}{llllllllll}
\boxed{1)}&10000000&2)&10000010&3)&10000011&4)&10000101&120)&01111110\\
5)&10000110&6)&10000111&7)&10001011&8)&10001100&119)&01111010\\
9)&10001101&10)&10001111&11)&10010001&12)&10010100&118)&01111001\\
13)&10011000&14)&10011001&15)&10011010&16)&10011110&117)&01111000\\
17)&10100000&18)&10100001&19)&10100010&20)&10100110&116)&01110100\\
21)&10101001&22)& 10101100&23)&10110000&24)&10110001&115)&01110001\\
25)&10110010&26)& 10110110&27&10111001&28)&10111100&114)&01101111\\
29)&11000000&30)&11000010&31)&11000011&32)&11000101&113)&01101101\\
33)&11000110&34)&11000111&35)&11001011&36)&11001100&112)&01101100\\
37)&11001101&38)&11001111&39)&11010001&40)&11010100&111)&01101011\\
41)&11011000&42)&11011001&43)&11011010&44)&11011110&110)&01100111\\
45)&11100000&46)&11100001&47)&11100010&48)&11100110&109)&01100110\\
49)&11101001&50)&11101100&51)&11110000&52)&11110001&108)&01100101\\
53)&11110010&54)& 11110110&55)&11111001&56)&11111100&107)&01100011\\
57)&00000011&\boxed{58)}&00000100&59)&00000101&60)&00000111&106)&01100010\\
\boxed{61)}&00001000&62)& 00001010&63)&00001011&64)&00001101&105)&01100000\\
65)&00001110&66)&00001111&67)&00010011&68)&00010100&104)&01011100\\
69)&00010101&70)& 00010111&71)&00011000&72)&00011010&103)&01011001\\
73)&00011011&74)&00011101&75)&00011110&76)&00011111&102)&01010110\\
77)&00100001&78)&00100100&79)&00101000&80)&00101001&101)&01010010\\
81)&00101010&82)&00101110&83)&00110000&84)&00110010&100)&01010001\\
85)&00110011&86)&00110101&87)&00110110&88)&00110111&99)&01010000\\
89)&00111011&90)&00111100&91)&00111101&92)&00111111&98)&01001100\\
\boxed{93)}&01000000&94)&01000001&95)&01000010&96)&01000110&
97)&01001001\\
\end{array}\end{equation*}

\normalsize

\subsubsection[$\fe(8,8)$ of $\sdim=120|128$]{$\boldsymbol{\fe(8,8)}$ of $\boldsymbol{\sdim=120|128}$} All inequivalent Cartan
matrices are as follows: \tiny
\begin{equation*}\label{e88}
\begin{array}{llllllllll}
\boxed{1)})&00000001&\boxed{2)}&00000010&\boxed{12)}&00100000&\boxed{6)}&00010000
&109)&11010101\\
5)&00001100&4)&00001001&7)&00010001&8)&00010010&110)&11010110\\
9)&00010110&10)&00011001&11)&00011100&3)&00000110&111)&11010111\\
13)&00100010&14)&00100011&15)&00100101&16)&00100110&112)&11011011\\
17)&00100111&18)&00101011&19)&00101100&20)&00101101&113)&11011100\\
21)&00101111&22)& 00110001&23)&00110100&24)&00111000&114)&11011101\\
25)&00111001&26)&00111010&27)&00111110&28)&01000011&115)&11011111\\
29)&01000100&30)& 01000101&31)&01000111&32)&01001000&116)&11100011\\
33)&01001010&34)&01001011&35)&01001101&36)&01001110&117)&11100100\\
37)&01001111&38)&01010011&39)&01010100&40)&01010101&118)&11100101\\
41)&01010111&42)&01011000&43)&01011010&44)&01011011&119)&11100111\\
45)&01011101&46)&01011110&47)&01011111&48)&01100001&120)&11101000\\
49)&01100100&50)&01101000&51)&01101001&52)&01101010&121)&11101010\\
53)&01101110&54)&01110000&55)&01110010&56)&01110011&122)&11101011\\
57)&01110101&58)&01110110&59)&01110111&60)&01111011&123)&11101101\\
61)&01111100&62)&01111101&63)&01111111&64)&10000001&124)&11101110\\
65)&10000100&66)&10001000&67)&10001001&68)&10001010&125)&11101111\\
69)&10001110&70)& 10010000&71)&10010010&72)&10010011&126)&11110011\\
73)&10010101&74)&10010110&75)&10010111&76)&10011011&127)&11110100\\
77)&10011100&78)&10011101&79)&10011111&80)&10100011&128)&11110101\\
81)&10100100&82)&10100101&83)&10100111&84)&10101000&129)&11110111\\
85)&10101010&86)&10101011&87)&10101101&88)&10101110&130)& 11111000\\
89)&10101111&90)&10110011&91)&10110100&92)&10110101&131)&11111010\\
93)&10110111&94)& 10111000&95)&10111010&96)&10111011&132)&11111011\\
97)&10111101&98)&10111110&99)&10111111&100)&11000001&133)&11111101\\
101)&11000100&102)&11001000&103)&11001001&104)&11001010&134)&11111110\\
105)&11001110&106)&11010000&107)&11010010&108)&11010011&\underline{135)}&
\underline{11111111}\\
 \end{array}
\end{equation*}

\normalsize

\newpage

\begin{landscape}
\section[Table. Dynkin diagrams for $p=2$]{Table. Dynkin diagrams for $\boldsymbol{p=2}$}\label{tbl}{}~{} \vskip 0.6 cm

\begin{center}
\footnotesize
\extrarowheight=2pt \begin{tabular}%
{|>{\PBS\raggedright\hspace{0pt}}m{36mm}|
 >{\PBS\raggedright\hspace{0pt}}m{40mm}|
 >{\PBS\raggedright\hspace{0pt}}m{13mm}|
 >{\PBS\raggedright\hspace{0pt}}m{17mm}|
 >{\PBS\raggedright\hspace{0pt}}m{13mm}|
 >{\PBS\centering\hspace{0pt}}m{12mm}|
 >{\PBS\raggedright\hspace{0pt}}m{25mm}|}
\hline
 Diagrams\centering & $\fg$\centering & $v$\centering & $ev$ &
 $od$ &$png$ & $ng\le \min(* \ , *)$\\ \hline
 &&&
 $k_\ev-2$ & $k_\od$ &$\ev$ &$2k_\ev-4, 2k_\od$ \\
 &&&
 $k_\od$ &$k_\ev-2$ & $\od$ &$2k_\ev-3, 2k_\od-1$ \\
 &&&
 $k_\od-2$ & $k_\ev$ &$\ev$ &$2k_\ev, 2k_\od-4$ \\
 &&&
 $k_\ev$ &$k_\od-2$ & $\od$ &$2k_\ev-1,2k_\od-3$ \\ \cline{4-7}
 &&&
 $k_\ev-1$ & $k_\od-1$ & &$2k_\ev-2, 2k_\od-1$ \\
 \raisebox{3.5em}[0pt]{$\left.\arraycolsep=0pt\begin{array}{l}
 1)\;\includegraphics{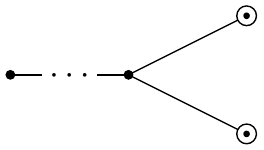}\\
 2)\;\includegraphics{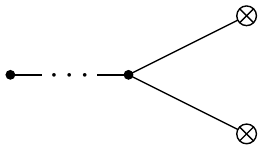}
 \end{array}\;\right\}$} &
 \raisebox{3.5em}[0pt]{$\begin{array}{l}\foo\fc(2;2k_\ev|2k_\od)\subplus\Kee I_0\\
 \text{if $k_\ev+k_\od$ is odd;}\\ \foo\fc(1;2k_\ev|2k_\od)\subplus\Kee I_0\\
 \text{if $k_\ev+k_\od$ is even.}\end{array}$}
 &
 \raisebox{3.5em}[0pt]{$k_{\ev}+k_{\od}$} &
 $k_\od-1$ &$k_\ev-1$ & &$2k_\ev-1, 2k_\od-2$ \\
\hline
 &&&
 $k_\ev-1$ & $k_\od$ &$\ev$ &$2k_\ev-2, 2k_\od$ \\
 &&&
 $k_\od$ &$k_\ev-1$ & $\od$ &$2k_\ev-1,2k_\od-1$ \\ \cline{4-7}
 &&&
 $k_\od-1$ & $k_\ev$ & $\ev$ &$2k_\ev, 2k_\od-2$ \\
 \raisebox{2.2em}[0pt]{$\left.\arraycolsep=0pt\begin{array}{l}
 3)\; \includegraphics{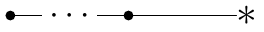}\\
 4)\;\includegraphics{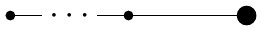}
 \end{array}\;\right\}$} &
 \raisebox{2.2em}[0pt]{$\foo^{(1)}_{\rm I\Pi}(2k_\ev+1|2k_\od)$} &
 \raisebox{2.2em}[0pt]{$k_{\ev}+k_{\od}$} &
 $k_\ev$ &$k_\od-1$ & $\od$ &$2k_\ev-1, 2k_\od-1$ \\
 \hline
 5)\;\includegraphics{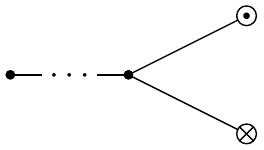}
 & $\begin{array}{l}\fpe\fc(2;m)\subplus\Kee I_0\\ \text{if $m$ is odd;}\\
 \fpe\fc(1;m)\subplus\Kee I_0\\ \text{if $m$ is even.}\end{array}$ & $m$&&&&\\ \hline
\end{tabular}

\end{center}

\subsection{Notation} The Dynkin diagrams in Table~\ref{tbl} correspond
to CM Lie superalgebras close to ortho-orthogonal and periplectic
Lie superalgebras. Each thin black dot may be $\motimes$ or $\odot$;
the last f\/ive columns show conditions on the diagrams; in the last
four columns, it suf\/f\/ices to satisfy conditions in any one row.
Horizontal lines in the last four columns separate the cases
corresponding to dif\/ferent Dynkin diagrams. The notations are:
\begin{enumerate}\itemsep=0pt
\item[] $v$ is the total number of nodes in the diagram;
\item[] $ng$ is the number of ``grey'' nodes $\motimes$'s among the
thin black dots;
\item[] $png$ is the parity of this number;
\item[] $ev$ and $od$ are the number of thin black dots such that the number
of $\otimes$'s to the left from them is even and odd, respectively.
\end{enumerate}
\end{landscape}

\normalsize

\section{Fixed points of symmetries of the Dynkin
diagrams}

\subsection{Recapitulation} For $p=0$, it is well known that the Lie
algebras of series $B$ and $C$ and the exceptions $F$ and $G$ are
obtained as the sets of f\/ixed points of the outer automorphism of an
appropriate Lie algebra of $ADE$ series. All these automorphisms
correspond to the symmetries of the respective Dynkin diagram. Not
all simple f\/inite dimensional Lie superalgebras can be obtained as
the sets of f\/ixed points of the symmetry of an appropriate Dynkin
diagram, but many of them can, see \cite{FSS}.

Recall Serganova's result \cite{Se} on outer automorphisms (i.e.,
the modulo the connected component of the unity of the automorphism
group) of simple f\/inite dimensional Lie superalgebras for $p=0$. The
symmetry of the Dynkin diagram of $\fsl(n)$ corresponds to the
transposition with respect to the side diagonal, conjugate in the
group of automorphisms of $\fsl(n)$ to the ``minus
transposition'' $X\mapsto -X^t$. In the super case, this
automorphism becomes $X\longmapsto -X^{st}$, where
\[
\mat{A&B\\C&D}^{st}=\mat{A^t&-C^t\\B^t&D^t}.
\]
This automorphism, seemingly of order 4, is actually of order 2
modulo the connected component of the unity of the automorphism
group, and is of order 4 only for $\fsl(2n+1|2m+1)$.

The queer Lie superalgebra $\fq(n)$ is obtained as the set of f\/ixed
points of the automorphism
\[
\Pi:\mat{A&B\\C&D}\longmapsto\mat{D&C\\B&A}\] of $\fgl(n|n)$
corresponding to the symmetry of the Dynkin diagram
\[
\stackrel{\text{\raisebox{2ex}{11}}}{\mcirc}-\dots-\stackrel{\text{\raisebox{2ex}{$1n$}}}{\mcirc}-\motimes-
\stackrel{\text{\raisebox{2ex}{21}}}{\mcirc}-
\dots-\stackrel{\text{\raisebox{2ex}{$2n$}}}{\mcirc}\;\;\longmapsto
\;\;\stackrel{\text{\raisebox{2ex}{21}}}{\mcirc}-
\dots-\stackrel{\text{\raisebox{2ex}{$2n$}}}{\mcirc}-\motimes-\stackrel{\text{\raisebox{2ex}{11}}}{\mcirc}-\dots-
\stackrel{\text{\raisebox{2ex}{$1n$}}}{\mcirc}
\]
which interchanges the identical maximal
parts $\mcirc-\dots-\mcirc$ \textit{preserving} the order of nodes;
whereas $\fpe(n)$ is the set of f\/ixed points of the composition
automorphism $\Pi\circ (-st)$.

\subsection{New result} The modular version of the above statements is
given in the next Theorem in which, speaking about ortho-orthogonal
and periplectic superalgebras, we distinguish the cases where the
fork node is grey or white ($g\fg(A)$ and $w\fg(A)$, respectively);
to squeeze the data in the table, we write~$\widehat{\fg}$ instead
of $\fg\subplus\Kee I_0$. We also need the following decomposable
Cartan matrices ($p=2$):
\[
\cN:=\begin{pmatrix}
\ev & 1 & 0 & 0\\
1 & \ev & 0 & 0\\
0 & 1 & \ev & 1 \\
0 & 0 & 1 &\ev
\end{pmatrix}
,\qquad \cM:=\begin{pmatrix}
\ev & 1 & 0 & 0\\
1 & 0 & 0 & 0\\
0 & 1 & \ev & 1 \\
0 & 0 & 1 &\ev
\end{pmatrix}.
\]

\subsubsection[The Lie algebra $\fg({\cal N})$]{The Lie algebra $\boldsymbol{\fg({\cal N})}$}

It is of $\dim 34$ and not
simple; it contains a simple ideal of $\dim=26$ which is
$\fo(1;8)^{(1)}/\fc$ and the quotient is isomorphic to $\fsl(3)$.

\subsubsection[The Lie superalgebra $\fg({\cal M})$]{The Lie superalgebra $\boldsymbol{\fg({\cal M})}$}

It is of $\sdim 18|16$
and not simple. Its even part is $\fhei(2) \oplus_c \fg(C)$, where
$\fhei(2)=\Span \{X^\pm, c\}$ and $c$ is the center of the Lie
algebra $\fg(C)$, where
\[
C:=\begin{pmatrix}
\ev & 0 & 0\\
0 & \ev& 1\\
0 & 1 & \ev
\end{pmatrix}.
\]

The brackets are as follows:
\[[X^\pm, \fg(C)^{(1)}]=0, \qquad [X^\pm, d]=
X^\pm, \qquad [X^+,X^-]=c,\] where $d$ is the grading operator of the
Lie algebra $\fg (C)$.

Now the Cartan subalgebra of $\fg ({\cal M})$ is generated by $h_3$,
$h_6$, $h_1+h_5$, $h_2+h_4$ and the highest weight vector of the module
$\fg({\cal M})_\od$ is $x_{32}+ x_{33}$, where
\[
x_{32}= [ [ [x_1,x_2], [x_3,x_4]] , [[x_3,x_6], [x_4,x_5]]],\qquad
x_{33}=[ [ [x_1,x_2], [x_3,x_6]] , [[x_2,x_3], [x_4,x_5]]].
\]

Its weight is $(0,0,1,0)$ (according to the ordering of the
generators of the Cartan subalgebra as above).

The restriction of the module to $\fhei(2)$ consists of 8 copies of
the 2-dimensional irreducible Fock module; the restriction to
$\fg(C)$ consists of 2 copies of an irreducible 8-dimensional
module.

The lowest weight vector is $y_{32}+y_{33}$ with weight $(0,0,1,0)$.

The Lie superalgebra $\fg({\cal M})$ has a simple ideal, of
$\sdim=10|16$ which is $\foo(1;4|4)^{(1)}/\fc$ (to be described
separately below) and the quotient is isomorphic to $\fsl(3)$.

\sssbegin{Theorem} If the Dynkin diagram of $i\fg(A)$ is symmetric,
it gives rise to an outer automorphism $\sigma$ whose fixed points
constitute the Lie superalgebra $(i\fg(A))^\sigma$ which occupies
the slot under $i\fg(A)$ in the following tables $(\ref{symm})$,
$(\ref{symm2})$, $(\ref{symm3})$.

$1) \; p=2:$ The order $2$ automorphisms of the $\fsl$ series
corresponding to the symmetries of Dynkin diagrams give the
following fixed points, where $\sigma=t$ is the transposition,
unless otherwise stated:
\begin{equation}\label{symm}
\begin{array}{l}\renewcommand{\arraystretch}{1.4}\begin{tabular}{|l|l|l|}
\hline
$\fsl(2n+1)$&$\fgl(2n)$&$\fsl(2k+1|2m)$\\
$\fo(2n+1)$&$\fo(2n)$&$\foo(2k+1|2m)$\\
\hline
$\fgl(n|n)$&$\fgl(n|n)$&$\fgl(2k|2m)$\\
$\fq(n)$, ~$\sigma=\Pi$&$\fpe(n)$, ~$\sigma=\Pi\circ (t)$&$\foo(2k|2m)$\\
\hline
\end{tabular}\; \renewcommand{\arraystretch}{1.4}\begin{tabular}{|l|}
\hline $\fgl(2k+1|2m+1)$\\
$\foo(2k+1|2m+1)$\\
\hline\end{tabular}\end{array}
\end{equation}

$2)\; p=2:$ The order $2$ automorphisms of the orthogonal and
ortho-orthogonal series give the following fixed points $($recall the
definition of $\hat\fg$ in \eqref{hat}$)$:
\begin{equation}\label{symm2}
\renewcommand{\arraystretch}{1.4}\begin{tabular}{|l|l|}
\hline
$\begin{array}{ll}\widehat{\foo\fc(2;2k_\ev|2k_\od)}&\text{for
$k_\ev+k_\od$
odd}\\
\widehat{\foo\fc(1;2k_\ev|2k_\od)}&\text{for $k_\ev+k_\od$
even}\end{array}$&$\begin{array}{ll}\widehat{\fo\fc(2;2k)}&\text{for
$k$ odd}\\
\widehat{\fo\fc(1;2k)}&\text{for $k$ even}\end{array}$\\ \hline
\end{tabular}
\end{equation}

$3)$ The following are the fixed points of order $2$ automorphisms
of the exceptional Lie (super)algebras for $p=3$, and also,
 for $p=2$, of the periplectic superalgebras, and of order $3$ automorphisms of the
orthogonal algebra and ortho-orthogonal superalgebras.
\begin{equation}\label{symm3}\renewcommand{\arraystretch}{1.4}
\begin{tabular}{|l|l|l|l|l|l|}
\hline
 $1\fg(2,3)$ & $2\fg(2,3)$ & $5\fg(2,3)$ &$5\fg(2,6)$&$2\fg(2,6)$&
 $\widehat{\fo\fc(1;8)}$ \\
 $\fpsl(2|2)$ & $\fsl(1|2)$ & $\fosp(3|2) $ & $\fg(1,6) $&$\fg(1,6) $&$\fgl(4)$ \\
 \hline
$g\widehat{\foo\fc(1;4|4)}$&$w\widehat{\foo\fc(1;4|4)}$&
$g\widehat{\fpe\fc(1;4)}$ &$w\widehat{\fpe\fc(1;4)}$ &
$g\widehat{\foo\fc(2;6|2)}$&
$w\widehat{\foo\fc(2;6|2)}$\\
 $\fgl(2|2) $ & $ \fgl (2|2)$ & $ \fgl (1|3) $&$ \fgl (1|3) $ &$ \fgl (2|2) $& $ \fgl (2|2) $\\
 \hline
\end{tabular}\end{equation}
Besides, $\fe(6)^\sigma=\fg({\cal N})$, whereas
\[\renewcommand{\arraystretch}{1.4}
\begin{array}{l}
25\fe(6,1)^\sigma\simeq
26\fe(6,1)^\sigma\simeq27\fe(6,1)^\sigma\simeq\\
1\fe(6,6)^\sigma\simeq7\fe(6,6)^\sigma\simeq5\fe(6,6)^\sigma\simeq33\fe(6,6)^\sigma\simeq
8\fe(6,6)^\sigma\simeq29\fe(6,6)^\sigma\simeq\\
32\fe(6,6)^\sigma\simeq10\fe(6,6)^\sigma\simeq
14\fe(6,6)^\sigma\simeq18\fe(6,6)^\sigma\simeq28\fe(6,6)^\sigma\simeq
36\fe(6,6)^\sigma\simeq
\fg({\cal M}).
\end{array}
\]
\end{Theorem}

\section[A realization of $\fg=\foo(4|4)^{(1)}/\fc$]{A realization of $\boldsymbol{\fg=\foo(4|4)^{(1)}/\fc}$}

This simple Lie
superalgebra $\fg$ admits a realization in which $\fg_\ev\simeq
\fhei(8)\subplus \Kee E$, where $\fhei(8)=\Span(p, q, c)$ with
$p=(p_1,\dots, p_4)$, $q=(q_1,\dots, q_4)$, and $E:=\sum
(p_i\partial_{p_i}+q_i\partial_{q_i})$ and with $c$ being central in
$\fhei$, and in which $\fg_\od$ is a copy of the Fock space (over
$\fhei(8)$) considered purely odd, i.e., as $\Pi(\Kee[p]/(p_1^2,
\dots, p_4^2))$. (Obviously, the indeterminates $p$ and $q$, as well
as $\xi$ and $\eta$, cf.\ Section~\ref{fock}, are interchangeable.)

Indeed, consider the following isomorphism
\begin{equation}\label{focks}
\begin{array}{l}\varphi:
\Pi(\Kee[p]/(p_1^2, \dots, p_4^2))\tto \Span(\varphi_0, \dots, \varphi_{1234}),\vspace{1mm}\\
\varphi_0:=\Pi(1),\  \varphi_{i}:=\Pi(p_i),\
\varphi_{ij}:=\Pi(p_ip_j), \  \dots, \
\varphi_{1234}:=\Pi(p_1p_2p_3p_4).
\end{array}\end{equation}
Now the multiplication is given by the following two tables, where
$D:=c+E$ to save space: \begin{equation}\label{focks1}
\renewcommand{\arraystretch}{1.4}\mbox{\tiny$
\begin{array}{|c|cccccccccccccccc|}
\hline& \varphi_{1234} & \varphi_{234}& \varphi_{134}
&\varphi_{124}& \varphi_{123}& \varphi_{34}&\varphi_{24} &
\varphi_{23} & \varphi_{14} & \varphi_{13} &\varphi_{12} &
\varphi_{4}& \varphi_{3}& \varphi_{2}& \varphi_{1}& \varphi_{0}
\\
\hline \varphi_{1234} & 0 & 0& 0& 0& 0& 0& 0& 0& 0& 0& 0 &
p_4 & p_3 & p_2 & p_1 &D\\
\hline \varphi_{234} & 0 & 0
 & 0 & 0 & 0 & 0 & 0 & 0 & p_4 & p_3 & p_2 & 0 & 0 & 0 & E
& q_1\\
\hline \varphi_{134} & 0 & 0 & 0 & 0 & 0 & 0& p_4 & p_3 & 0 & 0 &
p_1 & 0 & 0 & E & 0 & q_2
\\
\hline \varphi_{124} &0 & 0 & 0 & 0 & 0 & p_4 & 0 & p_2 & 0 & p_1 &
0 & 0 &
E & 0 & 0 & q_3\\
\hline \varphi_{123} &0 & 0 & 0 & 0 & 0 &p_3 & p_2 & 0 & p_1 & 0 &0&
E & 0 & 0 & 0 & q_4
\\
\hline \varphi_{34} &0 & 0 & 0& p_4 & p_3&0 &0 &0 &0 & 0 & D & 0&
0 & q_1 & q_2 &0\\
\hline
\varphi_{24} &0 & 0 & p_4 & 0 &p_2& 0& 0& 0& 0& D & 0 & 0 & q_1 & 0& q_3 & 0\\
\hline \varphi_{23} &0 & 0 & p_3 &p_2 &0 & 0 & 0 & 0 & D & 0&
0 & q_1 & 0 &0&q_4&0\\
\hline
\varphi_{14} &0 & p_4 & 0 & 0 & p_1 & 0 & 0 & D & 0 &0 &0 &0 & q_2 & q_3 & 0 &0\\
\hline
\varphi_{13}&0 & p_3 &0 & p_1 &0 &0 & D & 0& 0& 0& 0& q_2 & 0 &q_4 & 0 & 0\\
\hline \varphi_{12} & 0& p_2& p_1 &0 & 0 & D & 0 & 0 & 0 & 0 &0& q_3
& q_4 & 0 & 0&0
\\
\hline \varphi_{4} &p_4 & 0 & 0 & 0 & E& 0 & 0 & q_1 &0 & q_2 & q_3&
0 & 0 & 0 &0&0
\\
\hline \varphi_{3}&p_3 & 0 & 0 & E &0 & 0 & q_1 & 0 & q_2
& 0 & q_4 & 0 & 0 & 0 & 0 & 0 \\
\hline \varphi_{2} & p_2 & 0 & E &0 & 0 & q_1 & 0 & 0 &
q_3 & q_4 & 0 & 0 & 0 & 0 & 0 & 0 \\
\hline \varphi_{1} & p_1 & E &0 & 0 & 0 & q_2 & q_3 & q_4
& 0 & 0 & 0 & 0 & 0 & 0 & 0 &0\\
\hline
\varphi_{0} & D & q_1 & q_2 & q_3 & q_4 & 0 & 0 & 0& 0& 0& 0& 0& 0& 0& 0& 0 \\
\hline
\end{array}$}
\end{equation}
\begin{equation}\label{focks2}
\begin{array}{|c|cccccccccc|}
\hline& c & D & p_1 & p_2& p_3& p_4& q_1& q_2& q_3& q_4 \\
\hline \varphi_{0} &\varphi_{0} & \varphi_{0} & \varphi_{1} &
\varphi_{2} & \varphi_{3} & \varphi_{4} &0 &0 &0&0
\\
\hline \varphi_{1} & \varphi_{1} & 0 &0 &\varphi_{12} &\varphi_{13}
&
\varphi_{14} &\varphi_{0} & 0 &0 &0\\
\hline \varphi_{2} &\varphi_{2} & 0 & \varphi_{12} & 0 &
\varphi_{23} & \varphi_{24} &
0 & \varphi_{0} & 0 & 0\\
\hline \varphi_{3}
&\varphi_{3}&0 &\varphi_{13}&\varphi_{23}&0& \varphi_{34}&0&0&\varphi_{0}&0\\
\hline \varphi_{4} &\varphi_{4}& 0 & \varphi_{14} & \varphi_{24} &
\varphi_{34} & 0 & 0 & 0 & 0 & \varphi_{0}
\\
\hline \varphi_{12} & \varphi_{12} & \varphi_{12} &0 &0
&\varphi_{123}&\varphi_{124} &\varphi_{2} &\varphi_{1} &0 &0
\\
\hline \varphi_{13}&\varphi_{13}&
\varphi_{13}&0&\varphi_{123}&0&\varphi_{134}&\varphi_{3}&0&\varphi_{1}&0
\\
\hline
\varphi_{14} &\varphi_{14} & \varphi_{14} & 0 & \varphi_{124} & \varphi_{134} & 0 & \varphi_{4} & 0 & 0 &\varphi_{1}\\
\hline \varphi_{23} &\varphi_{23} & \varphi_{23} & \varphi_{123} & 0
& 0 & \varphi_{234} & 0 & \varphi_{3}& \varphi_{2} &0
\\
\hline \varphi_{24} &\varphi_{24} & \varphi_{24} & \varphi_{124} &
0& \varphi_{234} & 0 & 0 & \varphi_{4} & 0
&\varphi_{2}\\
\hline \varphi_{34}
&\varphi_{34}&\varphi_{34}&\varphi_{134}&\varphi_{234}&0&0&0&0&
\varphi_{4}&\varphi_{3}
\\
\hline \varphi_{123}
&\varphi_{123}&0&0&0&0&\varphi_{1234}&\varphi_{23}&\varphi_{13}&\varphi_{12}&0
\\
\hline
\varphi_{124}&\varphi_{124}&0 &0&0&\varphi_{1234}&0&\varphi_{24}&\varphi_{14}&0&\varphi_{12}\\
\hline \varphi_{134}&\varphi_{134}& 0 & 0 & \varphi_{1234} &
0 & 0 & \varphi_{34} & 0 & \varphi_{14} & \varphi_{13}\\
\hline \varphi_{234} & \varphi_{234} &0 &\varphi_{1234} &0 &0 &0 &0
&\varphi_{34}
&\varphi_{24} &\varphi_{23}\\
\hline \varphi_{1234} &\varphi_{1234}& \varphi_{1234} &0 &0 &0 &0 &
\varphi_{234} & \varphi_{134} &
\varphi_{124}&\varphi_{123}\\\hline\end{array}
\end{equation}

\ssbegin{Remark} If $p=0$, every irreducible module over a solvable
Lie algebra is 1-dimensional. A theorem, based on this fact, states
that {\sl any Lie superalgebra $\fg$ is solvable if and only if
$\fg_\ev$ is solvable}. The example above shows that if $p>0$, life
is much more interesting. \end{Remark}

We were unable to answer: For $2n\neq 8$, is there a simple Lie
superalgebra $\fG(2n)$ with $\fG(2n)_\ev\simeq \fhei(2n)\subplus
\Kee E$ and $\fG(2n)_\od\simeq \Pi(\text{Fock module over
$\fhei(2n)$})$?  In the next subsection we cite Irina
Shchepochkina's answer to this question.

\subsection{Shchepochkina's comments. Other simple Lie superalgebras\\
 with solvable even part}

 We consider $\fg=\fg_\ev\oplus \fg_\od$, where
\[\fg_\ev=\fhei(2n)\subplus \Kee\cdot E = \Span(p_1,\dots, p_n,
q_1,\dots, q_n, c, E), \text{~~with $c$ central},\] i.e.,
\[
[p_i,q_j]=\delta_{ij} c, \qquad [E,p_i]=p_i, \qquad [E,q_i]=q_i,
\]
and  $\fg_\od=\Pi\left(\Lambda(p_1,\dots, p_n)\right)$. The Lie
algebra $\fhei(2n)$ acts in $\fg_\od$ as in the Fock space:
\[
\ad_c|\fg_\od =\id_{\fg_\od},\qquad \ad_{p_i}|\fg_\od=p_i\cdot, \qquad
\ad_{ q_i}|\fg_\od=\partial_{p_i}.
\]

The space $\fg_\od$ is spanned by $\varphi_0:=1$, and
$\varphi_{i_1\dots i_k}:=\Pi(p_{i_1}\dots p_{i_k})$ for all sets $I$
of distinct indices. For any $I$, let $I^*$ denote the complementary
set to $\{1, \dots, n\}$. For clarity, we set
$\varphi^*_I=\varphi_{I^*}$.

How can the operator $E$ act in $\fg_\od$? Let us begin with
$E\varphi_0$:
\[
[q_i,E\varphi_0]=[[q_i,E],\varphi_0]+[E, [q_i,\varphi_0]]=0 \quad
\text{~~for all $i$}.
\]
But there is only one (up to a constant factor) element in $\fg_\od$
annihilated by all the $q_i$, namely $\varphi_0$. Hence
$E\varphi_0=\lambda \cdot \varphi_0$.

Since $\fg_\od$ is generated from $\varphi_0$ under the action of
operators $p_i$ of weight $1$ with respect to $E$, the action of $E$
on the monomial $\varphi$ is of the form: \[E\varphi=(\lambda + \deg
\varphi)\varphi,\quad \text{~~i.e., \ \ $\ad_E|\fg_\od=\lambda\cdot
\id|\fg_\od +\deg$}.\] Replacing $E$ by $E+\lambda\cdot c$ (this
does not af\/fect the commutation relations in $\fg_\ev$), we may
assume that $\ad_E|\fg_\od=\deg$.

Let us try to def\/ine the bracket on $\fg_\od$. I claim that it
suf\/f\/ices to determine the only bracket (determine $x$)
\[x=[\varphi_0,\varphi^*_0].\] Everything else follows from
the Jacobi identity. Indeed,
\[
 [\varphi_0,\varphi^*_i]=[\varphi_0,[q_i,\varphi^*_0]]=[[\varphi_0,q_i],\varphi^*_0]+
 [q_i,[\varphi_0,\varphi^*_0]]=[q_i,x],
\]
and the inverse induction on the degree of monomials $\varphi$
yields all the brackets $[\varphi_0,\varphi]$:
\begin{gather}\label{eq1}
[\varphi_0,\varphi^*_{I\cup
i}]=[\varphi_0,[q_i,\varphi^*_I]]=[[\varphi_0,q_i],\varphi^*_I]+
 [q_i,[\varphi_0,\varphi^*_I]]=[q_i,[\varphi_0,\varphi^*_I]].
\end{gather}
If we know the brackets of a monomial $\psi$ with all monomials of
the form $\varphi$, we can recover the brackets of the form
$[p_i\cdot \psi, \varphi]$:
\begin{gather}\label{eq2}
[p_i\cdot \psi, \varphi]=[[p_i, \psi], \varphi]=[p_i, [\psi,
\varphi]]+ [\psi, [p_i,\varphi]].
\end{gather}
Thus, by the induction on the degree of $\psi$ we recover all the
brackets from the brackets of the form $[\varphi_0,\varphi]$.
Equations~\eqref{eq1}--\eqref{eq2} imply that $[\fg_\od,\fg_\od]\subset
\fg'_\ev+\Kee\cdot x$, where $\fg'_\ev:=[\fg_\ev, \fg_\ev]$.

For $\fg$ to be simple, we should have $[\fg_\od,\fg_\od]= \fg_\ev$,
i.e., $x\notin \fg'_\ev$. But $\varphi_0$ and $\varphi^*_0$ are
eigenvectors of the operator $E$ of weight $0$ and $n$,
respectively. Hence their bracket, $x$, is an eigenvector of weight
$n$ and, since $x$ should have a non-zero projection to $E$, the
number $n$ must be even: \[x=\alpha \cdot E +\beta \cdot c,\qquad \text{where $\alpha\ne 0$.}\]

Note that we have one more degree of freedom: We can multiply all
elements of $\fg_\od$ by the same non-zero scalar. This helps us to
f\/ix $\alpha=1$. Thus, $x= E +\beta \cdot c$.

Now, using equations~\eqref{eq1}--\eqref{eq2} we can recover the general
formula for the bracket in $\fg_\od$. I claim that it is of the
following beautiful form
\begin{gather}\label{eq3}
[f,g]=c\!\int\!(\ad_{E}(f)g)\! + (E+\beta
c)\!\int\! (\ad_{c}(f)g)\!+\sum\!\left(p_i \!\int\!(\ad_{q_i}(f)g) \!+ q_i
\int(\ad_{p_i}(f)g\right)\!,\!\!\!
\end{gather}
where $\int$ is the Berezin integral = the coef\/f\/icient of the
highest term (once the basis of the Grassmann algebra is chosen).
To
see this, it suf\/f\/ices to verify that \eqref{eq3} is invariant with
respect to $\ad_{q_i}$ and $\ad_{p_i}$ (the invariance with respect
to $\ad_c$ and $\ad_E$ is obvious).

{\bf Here comes an incomplete argument.} It remains to verify the
Jacobi identity only for triples of odd elements, moreover, it
suf\/f\/ices to check it only for triples of the form $\varphi_0$,
$\varphi^*_0$, $\varphi$. We have
\begin{gather}\label{eq4}
[[\varphi_0, \varphi^*_0], \varphi]=(E+\beta\cdot c)\varphi=(\deg
\varphi +\beta)\varphi.
\end{gather}
What can one say about the sum
\begin{gather}\label{eq5}
[[\varphi_0, \varphi], \varphi^*_0]+[\varphi_0, [\varphi^*_0,
\varphi]]?
\end{gather}
Observe that the f\/irst summand can only be non-zero if
$\varphi=\varphi^*_0$ or $\varphi^*_i$, the second summand can only
be non-zero if $\varphi=\varphi_0$ or $\varphi_i$.

For $\varphi=\varphi^*_0$ and $\varphi=\varphi_0$, the Jacobi
identity holds for any $\beta$, whereas for $\varphi=\varphi^*_i$
and $\varphi=\varphi_i$ only if $\beta=1$ for $n=2$ and $\beta=0$
for $n>2$.

For $n=2$, the above-listed possibilities exhaust all possible
values of $\varphi$.

For $n=4$, there are also elements $\varphi=\varphi_{ij}$ yielding
$0$ in both formulas. However, for $n>4$ and
$\varphi=\varphi_{123}$, equation~\eqref{eq4} yields $\varphi_{123}$,
whereas equation~\eqref{eq5} yields $0$, so there is no Lie superalgebra.

Thus, for $n=2$, one may have a Lie superalgebra $\fG(4)$ with
brackets of its odd elements (the squares of each odd element being~0)
\begin{gather}
\label{Ira1}{}[\varphi_1, \varphi_2] = E, \qquad [\varphi_0,
\varphi_{12}] = E+c,
\end{gather}  where in the right hand sides
stand the elements that act on the odd part of the hypothetical Lie
superalgebra $\fG(4)$ as in~\eqref{focks2}, i.e., $E$ counts the
degree of the element of the Grassmann algebra; the other bracket
being def\/ined similar to \eqref{focks1}:
\begin{gather}
\label{Ira2}{}[\varphi_i, \varphi_{12}] = p_i, \qquad
[\varphi_0,\varphi_1] = q_{2}, \qquad [\varphi_0,\varphi_2] = q_{1},
\end{gather}
the element $p_i$ acts on $\fG(4)_\od$ as the multiplication by
$p_i$ (i.e., $[p_i, \varphi_0]=\varphi_i$, $[p_i,
\varphi_j]=\varphi_{ij}$ and so on), and $q_i$ as~$\partial_{p_i}$.

{\bf Here comes the complete argument}: In the above argument we
forgot that the Jacobi
identity for $p=2$ and odd elements is of the form \eqref{JI}, not of the usual form $[x,[x,x]]=0$.
 And taking $x=\varphi_0$, $y=p_1$ we
fail to satisfy the Jacobi identity  although it is so tempting to
set \[ x_1:=q_1,\qquad x_2:=\varphi_1,\qquad y_1:=p_1,\qquad
y_2:=\varphi_2
\] and (correctly) deduce that the relations between these
$x$'s and $y$'s yield the same Cartan matrix of $\fG(4)$ as that of
$\fo\fo^{(1)}_{\rm I\Pi}(1|4)$. However, due to squaring the space of
$\fG(4)$ is not a Lie superalgebra.

There are, however, simple Lie superalgebras with solvable even part
other than $\fo\fo^{(1)}_{\Pi\Pi}(4|4)/\fc$. Indeed, we know that
the Lie algebras $\fo^{(1)}_{\rm I}(n)$ are solvable (only) for $n=1$ and
2, and $\fo_\Pi(n)$ are solvable (only) for $n=2$ and~4. Therefore
the Lie superalgebras
\begin{gather}
\label{Ira3}
\begin{array}{llll}\fo\fo^{(1)}_{\rm II}(1|2), &
\fo\fo^{(1)}_{\rm I\Pi}(1|2),& \fo\fo^{(1)}_{\rm II}(2|2), &
\fo\fo^{(1)}_{\rm I\Pi}(2|2),\\
\fo\fo^{(1)}_{\rm I\Pi}(2|4),& \fo\fo^{(2)}_{\Pi\Pi}(2|4),&
\fo\fo^{(1)}_{\Pi\Pi}(4|4),&
\end{array}
\end{gather}
though simple (perhaps, modulo center), have solvable even parts.
One can not get a simple Lie superalgebra from
$\fo\fo^{(1)}_{\Pi\Pi}(2|2)$ passing to derived and factorizing.

\subsection*{Acknowledgements}

We are very thankful to A.~Lebedev for help (he not only clarif\/ied
the notion of $\fg(A)$ and roots, but also helped us to f\/igure out
the structure of $\fg(A)$ in the most complicated cases and
elucidate the notion of $p$-structure, he also listed inequivalent
Cartan matrices for the $\fe$-cases) and to I.~Shchepochkina for her
contribution. We thank A.~Protopopov for his help with our graphics,
see \cite{Pro}.
Constructive criticism of referees is most thankfully acknowledged;
the paper is much clearer now.

\addcontentsline{toc}{section}{References}
\LastPageEnding

\end{document}